\documentclass[12pt]{article}

\usepackage{amsmath, amsfonts, amssymb, amsthm}

\usepackage{graphicx, epsfig, xcolor}

\usepackage{xr}
\usepackage{enumerate, enumitem, array, multirow, multicol, arydshln, ctable}

\usepackage{natbib, url, hyperref}
\hypersetup{
    colorlinks=true,
    linkcolor=blue,
    anchorcolor=blue,
    citecolor=blue
}

\usepackage{algorithm, algpseudocode}

\newtheorem{theorem}{Theorem}
\newtheorem{remark}{Remark}
\newtheorem{lemma}{Lemma}
\newtheorem{corollary}{Corollary}
\newtheorem{proposition}{Proposition}

\algnewcommand\algorithmicinput{\textbf{Input:}}
\algnewcommand\Input{\item[\algorithmicinput]}
\algnewcommand\algorithmicoutput{\textbf{Output:}}
\algnewcommand\Output{\item[\algorithmicoutput]}

\makeatletter
\newcommand{\fdsy@scale}{1.0}
\newcommand\fdsy@mweight@normal{Book}
\newcommand\fdsy@mweight@small{Book}
\newcommand\fdsy@bweight@normal{Medium}
\newcommand\fdsy@bweight@small{Medium}

\DeclareFontFamily{U}{FdSymbolC}{}
\DeclareFontShape{U}{FdSymbolC}{m}{n}{
    <-7.1> s * [\fdsy@scale] FdSymbolC-\fdsy@mweight@small
    <7.1-> s * [\fdsy@scale] FdSymbolC-\fdsy@mweight@normal
}{}
\DeclareFontShape{U}{FdSymbolC}{b}{n}{
    <-7.1> s * [\fdsy@scale] FdSymbolC-\fdsy@bweight@small
    <7.1-> s * [\fdsy@scale] FdSymbolC-\fdsy@bweight@normal
}{}
\DeclareSymbolFont{arrows}{U}{FdSymbolC}{m}{n}
\SetSymbolFont{arrows}{bold}{U}{FdSymbolC}{b}{n}
\DeclareMathSymbol{\upvDash}{\mathrel}{arrows}{233}

\DeclareMathSymbol{\upmodels}{\mathrel}{arrows}{237}
\makeatother

\newcommand{\blind}{0}

\begin{document}
\renewcommand{\qedsymbol}{}

\def\spacingset#1{\renewcommand{\baselinestretch}%
{#1}\small\normalsize} \spacingset{1}


\if0\blind
{
  \title{\bf Variational Bernstein--von Mises theorem with increasing parameter dimension }
  \author
  {Jiawei Yan$^{1}$, Peirong Xu$^{1,\ast}$, and Tao Wang$^{1,2\ast}$\\
  	\\
  	\normalsize{$^{1}$Department of Statistics, Shanghai Jiao Tong University, Shanghai 200240, China.}\\
  	\normalsize{$^{2}$SJTU-Yale Joint Center of Biostatistics and Data Science, Shanghai Jiao Tong University,} \\
   \normalsize{Shanghai 200240, China.  }\\
  	\\
  	\normalsize{$^\ast$E-mail: neowangtao@sjtu.edu.cn, prxu@sjtu.edu.cn}
  }
  \date{}
  \maketitle
  
} \fi

\bigskip
\begin{abstract}

Variational Bayes (VB) provides a computationally efficient alternative to Markov Chain Monte Carlo, especially for high-dimensional and large-scale inference. However, existing theory on VB primarily focuses on fixed-dimensional settings or specific models. To address this limitation, this paper develops a finite-sample theory for VB in a broad class of parametric models with latent variables. We establish theoretical properties of the VB posterior, including a non-asymptotic variational Bernstein--von Mises theorem. Furthermore, we derive consistency and asymptotic normality of the VB estimator. An application to multivariate Gaussian mixture models is presented for illustration. 

\end{abstract}

\noindent%
\vfill

\newpage
\spacingset{1.9} 
\section{Introduction}
Variational inference, which approximates Bayesian posterior inference as an optimization problem, represents a significant advancement in statistical inference. Compared to Markov Chain Monte Carlo (MCMC) sampling, it offers a more computationally efficient alternative, particularly for high-dimensional models and large datasets \citep{Blei_2017, zhang2018advances}. This scalability has led to successful applications in fields requiring rapid solutions, such as commercial recommendation systems \citep{9826663}, credit risk control \citep{wang2023study}, computational biology \citep{NEURIPS2022_5e956fef}, and artificial intelligence \citep{chen2024auto}. For an overview of variational inference applied to specific models, see Section 5.2 of \cite{Blei_2017} and Section 1.1 of \cite{wang2019frequentist}. However, unlike MCMC methods, which are theoretically grounded, variational inference methods have few theoretical results supporting their use.

Most research on the statistical properties of variational inference or variational approximations has focused on the classical low-dimensional framework, where the number of parameters remains fixed as the sample size increases. \cite{2006Convergence} derived the convergence rate of the variational Bayesian (VB) estimator for Gaussian mixture models; \cite{2014On} obtained asymptotically valid standard errors for the VB estimator in linear models; \cite{doi:10.1137/16M1105384} established the existence and developed the asymptotic theory of the best Gaussian (or Gaussian mixture) variational approximation for a general posterior measure; and \cite{wang2019frequentist} established the Bernstein--von Mises theorem for the VB posterior in a broad class of Bayesian latent variable models, with \cite{NEURIPS2019_f7ae58c7} further extending the results under model misspecification. In a frequentist context, \cite{hall2011theory,2011ASYMPTOTIC} derived the consistency and asymptotic normality of the Gaussian variational approximate estimator for Poisson mixed models; \cite{10.1214/13-AOS1124} established these properties for the variational approximate estimator in stochastic block models; and \cite{westling2019beyond} studied the consistency and asymptotic normality of estimators based on variational approximations through a connection to M-estimation.

In contrast, more recent research has addressed the high-dimensional or non-asymptotic framework. \cite{ray2022variational} studied the optimal convergence rate of the VB posterior in high-dimensional sparse linear regression, while \cite{ray2020spike} extended this work to high-dimensional logistic regression, providing non-asymptotic convergence bounds. \cite{mukherjee2022variational} investigated the asymptotic accuracy of the VB approximation in high-dimensional Bayesian linear regression with product priors. \cite{10.1214/18-EJS1475} proved that VB is consistent for estimation in general mixtures, derived rates of convergence, and showed that the evidence lower bound maximization strategy is consistent for model selection. \cite{pmlr-v84-pati18a} and \cite{han2019statistical} derived finite-sample risk bounds for VB estimators in a class of Bayesian latent variable models, while \cite{yang2020alpha} extended these results to a family of VB methods known as $\alpha$-VB, with standard VB being a special case with $\alpha = 1$. \cite{alquier2016properties} demonstrated that, in the probably approximatively correct (PAC)-Bayesian context, the variational approximation has the same rate of convergence as the Gibbs posterior it approximates. Similarly, in the context of non-parametric and high-dimensional inference, \cite{zhang2020convergence} derived the same rate of convergence for the VB posterior as for the true posterior. 

As noted by \cite{wang2019frequentist}, many established properties of VB posteriors and estimators are restricted to specific models and priors. They addressed this problem by extending these properties to a broader class of parametric and semi-parametric models with latent variables, though still within the classical low-dimensional framework. This class includes Bayesian mixture models, Bayesian generalized linear mixed models, and Bayesian stochastic block models. To bridge the gap between theory and practice, this paper examines the finite-sample, or non-asymptotic, properties of VB posteriors and estimators for the same class of models.

Given a dataset $x=\{x_i\}^n_{i=1}$, where  each sample $x_i$ is associated with a local latent variable $z_i$ and a $p$-dimensional global latent variable $\theta$. Let $z=\{z_i\}^n_{i=1}$. We assume that the joint distribution of $\left(\theta, x, z\right)$ has the form
\begin{equation*}
p(\theta, x, z)=p(\theta) \prod_{i=1}^n p\left(z_i \mid \theta\right) p\left(x_i \mid z_i, \theta\right),
\end{equation*}
where $p(\theta)$ is the prior distribution of $\theta$, $p(z_i \mid \theta)$ is the conditional distribution of $z_i$ given $\theta$, and $p(x_i \mid z_i, \theta)$ is the likelihood of $x_i$ given $z_i$ and $\theta$. We consider the posterior inference problem for the global latent variable $\theta$. VB approximates the posterior distribution $p(\theta, z \mid x)$, which involves intractable integration or summation over the latent variables $\theta$ and $z$, by solving an optimization problem. Specifically, it employs a variational family of distributions and seeks the one that minimizes the Kullback--Leibler (KL) divergence from $p(\theta, z \mid x)$. In this paper, we focus on the mean-field variational family, where the latent variables are assumed to be mutually independent. Mathematically, the variational family for the joint distribution of $\theta$ and $z$ is given by
\begin{equation*}
\mathcal{Q}_{(\theta,z)} = \left\{q(\theta, z) : q(\theta, z) = q(\theta) q(z),\;q(\theta) \in \mathcal{Q}_\theta,\;q(z) \in \mathcal{Q}_z \right\},
\end{equation*}
where the variational families for $\theta$ and $z$ are defined as
\begin{equation*}
\mathcal{Q}_\theta = \left\{q(\theta) : q(\theta) = \prod\limits_{i=1}^p q(\theta_i) \right\}
\end{equation*}
and
\begin{equation*}
\mathcal{Q}_z = \left\{ q(z) : q(z) = \prod\limits_{j=1}^n q(z_j) \right\}.
\end{equation*}
The VB posterior $q^*(\theta, z) = q^*(\theta)q^*(z)$ is then obtained by solving the optimization problem
\begin{equation*}
q^*(\theta, z)=\underset{q(\theta, z) \in \mathcal{Q}_{(\theta,z)}}{\arg \min } \operatorname{KL}(q(\theta, z) \| p(\theta, z \mid x)),
\end{equation*}
which is equivalent to maximizing the evidence lower bound (ELBO) defined as
\begin{equation*}
\operatorname{ELBO}(q(\theta, z))=-\int q(\theta, z) \log \frac{q(\theta, z)}{p(\theta, x, z)}\mathrm{d}\theta\mathrm{d}z.
\end{equation*}
It can be easily shown that
\begin{equation*}
q^*(\theta) = \underset{q(\theta) \in \mathcal{Q}_\theta}{\arg\max}\sup \limits_{q(z) \in \mathcal{Q}_z} \operatorname{ELBO}(q(\theta, z)).
\end{equation*}
This paper aims to investigate the theoretical properties of the VB posterior $q^*(\theta)$ and the VB estimator $\hat{\theta}_n^* = \int\theta  q^*\left(\theta\right)\mathrm{d}\theta$. Specifically, under certain conditions:
\begin{itemize}
    \item[(1)] We provide a non-asymptotic upper bound on the mass of the VB posterior outside a small neighborhood of the population-optimal parameter, establishing the consistency of the VB posterior.
    \item[(2)] We derive a non-asymptotic upper bound on the total variation distance between the VB posterior and the KL minimizer of a normal distribution centered at the population-optimal parameter, thereby establishing the asymptotic normality of the VB posterior.
    \item[(3)] We establish a non-asymptotic upper bound on the $L_2$ error for the VB estimator around the population-optimal parameter, demonstrating its consistency as the number of parameters diverges.
    \item[(4)] Finally, we leverage these non-asymptotic bounds to establish the asymptotic normality of the VB estimator as the number of parameters diverges.
\end{itemize}

The rest of the paper is organized as follows. In Section \ref{sec: preliminaries}, we introduce the VB ideal posterior, derive the local asymptotic normality expansion of the variational log-likelihood, and establish a finite-sample bound for the normalizing constant of the VB ideal posterior. Section \ref{sec: non-asymptotic properties of the VB bridge} introduces the VB bridge, characterizes its connection with the VB posterior, and derives its non-asymptotic properties. In Section \ref{sec: non-asymptotic properties of the VB posterior}, we establish the convergence rate of the VB estimator and present a non-asymptotic variational Bernstein--von Mises theorem for the VB posterior. In Section \ref{sec: applications}, we present an application to multivariate Gaussian mixture models. Finally, Section \ref{sec: conclusion} provides a brief discussion. All proofs of the theoretical results are provided in the Appendix.

\section{Preliminaries}
\label{sec: preliminaries}

\subsection{Notation}
For each integer $d \geq 1$, let $\mathbb{R}^d$ denote the $d$-dimensional Euclidean space. For any vector $u=$ $\left(u_1, \ldots, u_d\right)^{\top} \in \mathbb{R}^d$, let $\|u\|=\sqrt{\sum_{i=1}^d u_i^2}$ represent its $L_2$ norm. Additionally, let $\operatorname{diag}(u)$ denote the diagonal matrix with diagonal elements equal to the components of $u$. For any symmetric matrix $A=\left(a_{i j}\right) \in \mathbb{R}^{d \times d}$, denote its largest eigenvalue by $\lambda_{\max }(A)$, its smallest eigenvalue by $\lambda_{\min }(A)$, and the maximum absolute value of its eigenvalues by $\rho(A)$. The determinant and trace of $A$ are denoted by $\det\left(A\right)$ and $\operatorname{tr}(A)$, respectively. Also, define $\operatorname{diag}(A)$ as the diagonal matrix with the same diagonal elements as $A$. For a differentiable function $f: \mathbb{R}^d \rightarrow \mathbb{R}$, the gradient vector and Hessian matrix of $f$ at $x$ are denoted by $\nabla f(x)$ and $\nabla^2 f(x)$, respectively. For any set $B$, define $I_B(x)$ as the indicator function, which equals 1 if $x \in B$ and 0 otherwise. For two sequences of non-negative numbers $\{a_n\}$ and $\{b_n\}$, the notation $a_n \lesssim b_n$ means there exists a constant $C>0$, independent of $n$, such that $a_n \leq Cb_n$ for sufficiently large $n$. The notation $a_n \asymp b_n$ denotes that both $a_n \lesssim b_n$ and $b_n \lesssim a_n$ hold. For two sequences of random variables, $\{R_{1n}\}$ and $\{R_{2n}\}$, the notation $R_{1n} = O_P(R_{2n})$ means that for any $\epsilon>0$, there exists a constant $C>0$ such that $P\left(\left|R_{1n}/R_{2n}\right|>C\right)<\epsilon$ for sufficiently large $n$.

\subsection{The VB ideal posterior}

In this subsection, we introduce the concept of the VB ideal posterior, following the approach of \cite{wang2019frequentist}. Assume, for the moment, that $\theta$ is unknown but fixed. The data log-likelihood is given by
\begin{equation*}
L_n(\theta; x)= \log \int p(x, z \mid \theta)\mathrm{d}z,
\end{equation*}
which is challenging to evaluate directly. To address this, the variational log-likelihood provides a convenient approximation, defined as
\begin{equation*}
M_n(\theta; x) = \sup\limits_{q(z) \in \mathcal{Q}_z} \int q(z) \log \frac{p(x, z \mid \theta)}{q(z)}\mathrm{d}z=\sum_{i=1}^n m(\theta; x_i),
\end{equation*}
where
\begin{equation*}
m(\theta; x_i)=\sup\limits_{q(z_i)} \int q(z_i) \log \frac{p(x_i, z_i \mid \theta)}{q(z_i)}\mathrm{d}z_i.
\end{equation*}
By replacing $L_n(\theta; x)$ with $M_n(\theta; x)$, we obtain a specific example of the Gibbs posterior:
\begin{equation*}
\pi^*(\theta \mid x) = \frac{p(\theta) \exp \left\{M_n(\theta ; x)\right\}}{\int p(\theta) \exp \left\{M_n(\theta ; x)\right\}\mathrm{d}\theta}.
\end{equation*}
We refer to this as the VB ideal posterior.

Denote $\theta^*$ as the value of $\theta$ that maximizes the expected
variational log-likelihood:
\[\theta^* = \underset{\theta \in \Theta}{\operatorname{argmax}}\;\mathbb{E}_{\theta_0} \left\{M_n(\theta; x)\right\}.\] 
Here, $\Theta$ is an open subset of $\mathbb{R}^p$, $\theta_0$ denotes the true value of $\theta$, and the expectation is taken with respect to $x$ under the distribution
\begin{equation*}
\int p(x, z \mid \theta = \theta_0)\mathrm{d}z.
\end{equation*}
As we will see in the next section, $\theta^*$ plays a key role in the theoretical properties of the VB posterior. Before proceeding, in the following two subsections, we study the theoretical properties of the variational log-likelihood $M_n(\theta; x)$ as well as the normalizing constant
\begin{equation*}
\int p(\theta)\exp\left\{M_n(\theta; x)\right\}\mathrm{d}\theta
\end{equation*}
in the VB ideal posterior.

\subsection{The local asymptotic normality expansion}
\label{subsection: the non-asymptotic LAN expansion of the variational log-likelihood}

Let $D_0^2=-\nabla^2 \mathbb{E}_{\theta_0} M_n\left(\theta^*; x\right)$ and $\Sigma_0^2=\operatorname{Var}_{\theta_0}\left\{\nabla M_n\left(\theta^* ; x\right)\right\}$. For any $r>0$, define 
$\Theta_0(r)=\left\{\theta \in \Theta:\left\|D_0\left(\theta-\theta^*\right)\right\| \leq r\right\}$. 
The following conditions are needed to establish the local asymptotic normality (LAN) expansion of $M_n\left(\theta; x\right)$.

\begin{enumerate}[label=(C\arabic*),series=C]

\item \label{condition: local, twice continuously differentiable}
Let $D^2(\theta) =-\nabla^2 \mathbb{E}_{\theta_0}M_n\left(\theta; x\right)$. For each $r\leq r_0$, there exists a constant $\delta_1(r)\in\left(0,1/2\right]$ such that
\begin{equation*}
\sup _{\theta\in \Theta_0(r)}\rho\left(D_0^{-1} D^2(\theta) D_0^{-1}-I_p\right)\leq \delta_1(r),
\end{equation*}
which implies that
\begin{equation*}
\sup _{\theta\in \Theta_0(r)}\left|\frac{2 \mathbb{E}_{\theta_0}\left\{M_n\left(\theta^*;x\right)-M_n\left(\theta;x\right)\right\}}{\left\|D_0\left(\theta-\theta^*\right)\right\|^2}-1\right| \leq \delta_1(r).
\end{equation*}
The constant $r_0$ will be specified later.

\item \label{condition: local, exponential moment}
    There exist constants $\tau_1\geq1$ and $\tau_2>0$ such that for all $|\lambda| \leq \tau_2$,
    \begin{equation*}
\sup _{\gamma \in \mathbb{R}^p} \log \mathbb{E}_{\theta_0} \exp \left[\lambda \frac{\gamma^{\top} \nabla \left\{M_n\left(\theta^* ; x\right)-\mathbb{E}_{\theta_0}M_n\left(\theta^* ; x\right)\right\}}{\left\|\Sigma_0 \gamma\right\|}\right] \leq \frac{\tau_1^2 \lambda^2}{2}.
    \end{equation*}

\item \label{condition: global, exponential moment, variational} There exists a constant $\tau_3>0$, and, for each $r>0$, a constant $\delta_2(r)>0$ such that $\sup_r \delta_2(r) / p$ is bounded from below and 
\begin{equation*}
\sup _{\theta\in \Theta_0(r)}\sup _{\gamma_1, \gamma_2 \in \mathbb{R}^p} \log \mathbb{E}_{\theta_0} \exp \left[\frac{\lambda}{\tau_3} \frac{\gamma_1^{\top} \nabla^2\left\{M_n\left(\theta; x\right)-\mathbb{E}_{\theta_0} M_n\left(\theta; x\right)\right\}\gamma_2}{\left\|D_0 \gamma_1\right\| \cdot\left\|D_0 \gamma_2\right\|}\right] \leq \frac{\tau_1^2 \lambda^2}{2},
\end{equation*}
for all $\left|\lambda\right|\leq\delta_2(r).$ 

\end{enumerate}

Conditions \ref{condition: local, twice continuously differentiable}--\ref{condition: global, exponential moment, variational} on the variational log-likelihood $M_n(\theta; x)$ resemble those in \cite{Spokoiny2012} and \cite{spokoiny2014bernsteinvonmises}, where the non-asymptotic properties of the data log-likelihood $L_n(\theta; x)$ were studied under possible model misspecification. While \ref{condition: local, twice continuously differentiable} and \ref{condition: local, exponential moment} are local conditions focusing on the behavior of $M_n\left(\theta; x\right)$ at or in the vicinity of $\theta^*$, \ref{condition: global, exponential moment, variational} is a global condition concerning the overall behavior of $M_n\left(\theta; x\right)$. Together, these conditions imply that $M_n(\theta; x)$ behaves similarly to $L_n(\theta; x)$.

\begin{lemma}[Non-asymptotic LAN expansion]
\label{LAN}
Let $h=\sqrt{n}\left(\theta-\theta^*\right)$, $V_{\theta_0}=D_0^2/n$, and $\Delta_{n,\theta_0}=V_{\theta_0}^{-1}\nabla M_n(\theta^*;x) / \sqrt{n}$. Define
\begin{equation*}
z_{\mathbb{H}}(y)= \begin{cases}\sqrt{4p+2y}, & \text { if } 4p+2y \leq \tau_2^2, \\ \tau_2^{-1} y+\frac{1}{2}\left(\tau_2^{-1} 4p+\tau_2\right), & \text { if } 4p+2y>\tau_2^2,\end{cases}
\end{equation*}
and
\begin{equation*}
\Delta(r_0,y)=r_0^2\left\{\delta_1(r_0)+6\tau_1 z_{\mathbb{H}}(y)\tau_3\right\}.
\end{equation*}
Under Conditions \ref{condition: local, twice continuously differentiable}--\ref{condition: global, exponential moment, variational}, with probability at least $ 1 - e^{-y}$, we have
\begin{equation*}
\sup _{h: h^\top V_{\theta_0}h\leq r_0^2} \left| M_n\left(\theta^*+\frac{h}{\sqrt{n}};x\right)-M_n(\theta^*;x)-h^{\top}V_{\theta_0}\Delta_{n, \theta_0}+\frac{1}{2}h^{\top}V_{\theta_0}h\right|\leq\Delta(r_0,y).
\end{equation*}
\end{lemma}

\textit{Proof sketch of Lemma \ref{LAN}.} This lemma is directly derived from Proposition 3.2 of \cite{spokoiny2014bernsteinvonmises}.

\begin{remark}
\label{Remark 1}
Lemma \ref{LAN} is particularly useful when $\Delta(r_0,y)$ is relatively small. For simplicity, consider moderate high-dimensional cases where $n / p$ is sufficiently large. As discussed in Section 5.1.2 of \cite{Spokoiny2012} and Section 2.4.3 of \cite{spokoiny2014bernsteinvonmises}, in regular cases, we have $\delta_1(r_0)\asymp r_0/ \sqrt{n}$ and $\tau_3\asymp 1 / \sqrt{n}$. Consequently, if $y = kp$ for some $k>0$, then with probability at least $1-e^{-kp}$,
\begin{equation*}
\Delta(r_0,y)=\Delta(r_0,kp)\lesssim\frac{r_0^2\left(r_0+\sqrt{p}\right)}{\sqrt{n}}.
\end{equation*}
If we further assume that $r_0^2\lesssim p$, it follows that 
\begin{equation*}
\sup _{h: h^\top V_{\theta_0}h\leq r_0^2} \left| M_n\left(\theta^*+\frac{h}{\sqrt{n}};x\right)-M_n(\theta^*;x)-h^{\top}V_{\theta_0}\Delta_{n, \theta_0}+\frac{1}{2}h^{\top}V_{\theta_0}h\right|\lesssim\sqrt{\frac{p^3}{n}},
\end{equation*}
with probability at least $1-e^{-kp}$.
\end{remark}

\subsection{The bracketing bound for the normalizing constant}

In this subsection, we use the LAN expansion of $M_n(\theta;x)$ to derive a finite-sample bound on the normalizing constant of the VB ideal posterior. Along with Conditions \ref{condition: local, twice continuously differentiable}--\ref{condition: global, exponential moment, variational}, the following additional conditions are required.

\begin{enumerate}[label=(C\arabic*),resume=C]

\item \label{condition: global identification property} 
    For each $r\geq r_0$, there exists a constant $\delta_3(r)>0$ such that $r\delta_3(r)$ is non-decreasing and 
    \begin{equation*}
\inf_{\{\theta:\left\|D_0\left(\theta-\theta^*\right)\right\|=r\}}\mathbb{E}_{\theta_0}\left\{M_n\left(\theta^*;x\right)-M_n\left(\theta;x\right)\right\}\geq\delta_3(r)r^2.
\end{equation*}

\item \label{condition: identifiability} There exists a constant $\tau_4>0$ such that $\tau_4^2 D_0^2 \geq \Sigma_0^2$. 

\item \label{condition: prior} The prior $p(\theta)$ is continuous and positive on a neighborhood of $\theta^*$, and there exists a constant $\tau_5>0$ such that $\rho\left(\nabla^2 \log p(\theta)\right)\leq \tau_5$.

\end{enumerate}

While \ref{condition: global identification property} is a global condition that concerns the behavior of $M_n\left(\theta;x\right)$ in the complement of $\Theta_0(r_0)$, \ref{condition: identifiability} serves as an identifiability condition. The prior condition \ref{condition: prior}, adapted with minor modifications from Assumption 1 in \cite{wang2019frequentist}, is a mild assumption satisfied by many commonly used priors, such as the standard multivariate Gaussian with zero mean and identity covariance matrix.

\begin{lemma}[Non-asymptotic Laplace approximation]
\label{Lemma: the bracketing result of the logarithm of the denominator of the VB Ideal}
Under Conditions \ref{condition: local, twice continuously differentiable}--\ref{condition: prior}, if $\rho\left(V_{\theta_0}^{-1}\right)\lesssim\sqrt{n}$, then for sufficiently large $r_0^2/p$,
\begin{eqnarray*}
& & \left|\log\int p(\theta) \exp \left\{M_n(\theta ; x)\right\} \mathrm{d} \theta+\frac{1}{2}\log\det\left(V_{\theta_0}\right)+\frac{p}{2}\log n\right.\\
& & \left.-M_n(\theta^*;x)-\log p\left(\theta^*\right)-\frac{p}{2}\log(2\pi)-\frac{1}{2} \Delta_{n,\theta_0}^{\top}V_{\theta_0}\Delta_{n,\theta_0}\right|\\
&\lesssim& \Delta(r_0,y) + \frac{1}{\sqrt{n}} \left\{\left\|\Delta_{n,\theta_0}\right\|^2+\left\|\nabla\log p\left(\theta^*\right)\right\|^2\right\} + \frac{1}{n}\mathrm{tr}\left(V_{\theta_0}^{-1}\right),
\end{eqnarray*}
with probability exceeding $1-5e^{-y}-e^{-cr_0^2}$ for some constant $c>0$.
\end{lemma}

\textit{Proof sketch of Lemma \ref{Lemma: the bracketing result of the logarithm of the denominator of the VB Ideal}.}
We decompose the normalizing constant as follows:
\begin{equation*}
\int p(\theta) \exp \left\{M_n(\theta ; x)\right\} \mathrm{d} \theta = \int_{\Theta_0(r_0)}p(\theta) \exp \left\{M_n(\theta ; x)\right\} \mathrm{d} \theta + \int_{\Theta^c(r_0)}p(\theta) \exp \left\{M_n(\theta ; x)\right\} \mathrm{d} \theta.
\end{equation*}
The first term is bounded using Lemma \ref{LAN}. The upper bound for the second term is established via Proposition 3.6 in \cite{spokoiny2014bernsteinvonmises}, while the lower bound is derived based on the tail behavior of the standard multivariate Gaussian distribution (refer to Lemma \ref{lemma: The Tail Moment Control for the p-Dimensional Normal Distribution} below).

In the classical setting, where the dimension of $\theta$ is fixed and the number of samples tends to infinity, the Laplace method is commonly used to approximate the normalizing constant \citep{wang2019frequentist}. Lemma \ref{Lemma: the bracketing result of the logarithm of the denominator of the VB Ideal} extends this approximation to the finite-sample context, where the classical results do not apply.

\begin{lemma}
\label{lemma: The Tail Moment Control for the p-Dimensional Normal Distribution}
Let $N\left(l;0,I_p\right)$ denote the density function of a $p$-dimensional multivariate normal distribution with zero mean and identity covariance matrix. For any $B>0, K\in\mathbb{N}$, and $p\in\mathbb{N}^+$, if $\left(p-1+K\right)/B^2\leq\alpha$ for an independent constant $ \alpha<1$, we have
\begin{equation*}
\int_{\left\{l:\left\|l\right\|>B\right\}} N\left(l;0,I_p\right)\left\|l\right\|^K\mathrm{d}l\asymp\frac{B^{p-2+K}}{(p-2)!!}\exp\left(-\frac{B^2}{2}\right).
\end{equation*}
\end{lemma}

\section{The VB bridge and its theoretical properties}
\label{sec: non-asymptotic properties of the VB bridge}
In this section, we introduce the concept of the VB bridge, examine its relationship with the VB posterior, and derive its non-asymptotic properties. These results are then utilized in Section \ref{sec: non-asymptotic properties of the VB posterior} to establish the non-asymptotic properties of the VB posterior. In the rest of the paper, we focus on the mean-field Gaussian variational family.

\subsection{The VB bridge}
\label{subsec: The VB bridge}
We define the VB bridge as the KL minimizer of the VB ideal posterior $\pi^*(\theta \mid x)$: 
\begin{equation*}
q^{\ddagger}(\theta)=\underset{q(\theta) \in \mathcal{Q}_\theta}{\arg \min } \operatorname{KL}\left(q(\theta) \| \pi^*(\theta \mid x)\right).
\end{equation*}
A straightforward calculation shows that $q^{\ddagger}(\theta)$ is the maximizer over $\mathcal{Q}_\theta$ of the functional
\begin{eqnarray*}
& & \log \int p(\theta) \exp \left\{M_n(\theta ; x)\right\} \mathrm{d} \theta-\mathrm{KL}\left(q(\theta) \| \pi^*(\theta \mid x)\right)\\
&=& \int q(\theta) \left(\log \left[p(\theta) \exp\left\{\sup_{q(z) \in \mathcal{Q}_z} \int q(z) \log \frac{p(x,z | \theta)}{q(z)}\mathrm{d}z \right\}\right] - \log q(\theta)\right) \mathrm{d}\theta.
\end{eqnarray*}

Recall that the VB posterior $q^*(\theta)$ is the maximizer over $\mathcal{Q}_\theta$ of the functional
\begin{eqnarray*}
\mathrm{ELBO}_p(q(\theta)) &\triangleq& \sup_{q(z) \in \mathcal{Q}_z} \operatorname{ELBO}(q(\theta, z))\\
&=& \sup_{q(z) \in \mathcal{Q}_z}\int q(\theta) \left(\log \left[p(\theta) \exp\left\{\int q(z) \log \frac{p(x,z | \theta)}{q(z)}\mathrm{d}z \right\}\right] - \log q(\theta)\right) \mathrm{d}\theta.
\end{eqnarray*}
The functionals for the VB bridge and the VB posterior are thus similar in form, differing only in the placement of the supremum over $\mathcal{Q}_z$. To explicitly characterize this difference, we need the following  conditions.

\begin{enumerate}[label=(C\arabic*),resume=C]
\item \label{condition: global, exponential moment, complete}
There exists a constant $\tau_6>0$, and, for each $r>0$, a constant $\delta_4(r)>0$ such that $\sup_r \delta_4(r) / p$ is bounded from below and 
\begin{equation*}
\sup _{\theta\in \Theta_0(r)}\sup _{\gamma_1, \gamma_2 \in \mathbb{R}^p} \log \mathbb{E}_{\theta_0} \exp \left[\frac{\lambda}{\tau_6} \frac{\gamma_1^{\top} \nabla^2\left\{\log p(x, z \mid\theta)-\mathbb{E}_{\theta_0}\log p(x, z \mid\theta)\right\}\gamma_2}{\left\|D_0 \gamma_1\right\| \cdot\left\|D_0 \gamma_2\right\|}\right] \leq \frac{\tau_1^2 \lambda^2}{2},
\end{equation*}
for all $\left|\lambda\right|\leq\delta_4(r).$

\item \label{condition: lipschitz}
There exist positive constants $C_1$, $C_2$, $s_1$, and $s_2$, such that for any $\theta_1, \theta_2\in \Theta$,
\begin{equation*}
\rho\left(\nabla^2\mathbb{E}_{\theta_0}M_n\left(\theta_1; x\right)-\nabla^2\mathbb{E}_{\theta_0}M_n\left(\theta_2; x\right)\right) \leq nC_1\left\|\theta_1-\theta_2\right\|^{s_1},
\end{equation*}
and
\begin{equation*}
\rho\left(\nabla^2\mathbb{E}_{\theta_0}\log p(x, z \mid\theta_1)-\nabla^2\mathbb{E}_{\theta_0}\log p(x, z \mid\theta_2)\right)\leq nC_2\left\|\theta_1-\theta_2\right\|^{s_2}.
\end{equation*}
\end{enumerate}

\begin{lemma}
\label{lemma: the gap between two functionals}
Let $c(x)=\log \int p(\theta) \exp \left\{M_n(\theta ; x)\right\}\mathrm{d} \theta$. Assume Conditions \ref{condition: global, exponential moment, variational}, \ref{condition: global, exponential moment, complete}, and \ref{condition: lipschitz} hold. Then, for any $q(\theta)=N\left(\theta;\mu,V/n\right)$, we have
\begin{equation*}
0 \leq  c(x)-\mathrm{KL}\left(q(\theta) \| \pi^*(\theta \mid x)\right)- \mathrm{ELBO}_p(q(\theta)) \lesssim p\rho\left(V\right)\left[\left(\frac{p\rho\left(V\right)}{n}\right)^{\frac{s_1}{2}}+\left(\frac{p\rho\left(V\right)}{n}\right)^{\frac{s_2}{2}}\right],
\end{equation*}
with probability at least $1-e^{-cp}$ for some constant $c>0$.
\end{lemma}

\textit{Proof sketch of Lemma \ref{lemma: the gap between two functionals}.} We employ a strategy similar to the proof of Lemma 4 in \cite{wang2019frequentist}, utilizing Taylor expansions around the mean vector $\mu$ to analyze the functional discrepancy.

\begin{remark}
\label{remark: the gap between two functionals}
In deriving the subsequent theoretical results, we only need to focus on the normal distribution
\begin{equation*}
q_0(\theta)=N\left(\theta;\theta^*+ \frac{\Delta_{n, \theta_0}}{\sqrt{n}}, \frac{\mathrm{diag}^{-1}\left(V_{\theta_0}\right)}{n}\right),
\end{equation*}
which will be referenced later. Accordingly, Condition \ref{condition: lipschitz} can be relaxed to
\begin{equation*}
\rho\left(\nabla^2\mathbb{E}_{\theta_0}M_n\left(\theta^*+\frac{\Delta_{n, \theta_0}}{\sqrt{n}}; x\right)-\nabla^2\mathbb{E}_{\theta_0}M_n\left(\theta; x\right)\right) \leq nC_1\left\|\theta^*+\frac{\Delta_{n, \theta_0}}{\sqrt{n}} -\theta\right\|^{s_1}
\end{equation*}
and
\begin{equation*}
\rho\left(\nabla^2\mathbb{E}_{\theta_0}\log p\left(x, z \mid\theta^*+\frac{\Delta_{n, \theta_0}}{n}\right)-\nabla^2\mathbb{E}_{\theta_0}\log p(x, z \mid\theta)\right)\leq nC_2\left\|\theta^*+\frac{\Delta_{n, \theta_0}}{\sqrt{n}}-\theta\right\|^{s_2}.
\end{equation*}
Lemma \ref{lemma: the gap between two functionals} can then be restricted to $q(\theta)$ with mean $\theta^*+\Delta_{n, \theta_0}/\sqrt{n}$, including  the specific case $q_0(\theta)$. If $\rho\left(V\right)\lesssim1$ and $\min\left(s_1,s_2\right)=1$, then we have 
\begin{equation*}
  \log \int p(\theta) \exp \left\{M_n(\theta ; x)\right\} \mathrm{d} \theta-\mathrm{KL}\left(q_0(\theta)\| \pi^*(\theta \mid x)\right)-\mathrm{ELBO}_p(q_0(\theta))=O_p\left(\sqrt{\frac{p^3}{n}}\right).
\end{equation*}
Please refer to Appendix \ref{Some Verifications of Design Conditions} for details. 
\end{remark}

\subsection{Non-asymptotic properties of the VB bridge}

Now, define the VB bridge $q^{\ddagger}(\theta)= N\left(\theta;\hat{\theta}_n^\ddagger, \operatorname{diag}\left(\hat{v}_n^\ddagger\right)/n\right)$, where $\hat{v}_n^\ddagger=\left(\hat{v}_{n1}^{\ddagger}, \ldots, \hat{v}_{np}^{\ddagger}\right)^\top$. Let $v_{\theta_0, ii}$ denote the $i$-th diagonal element of $V_{\theta_0}$, and define 
\begin{equation*}
R_n=\Delta(r_0,y)+\frac{1}{\sqrt{n}} \left\{\left\|\Delta_{n,\theta_0}\right\|^2+\left\|\nabla\log p\left(\theta^*\right)\right\|^2\right\} +\frac{1}{n} \left\{\mathrm{tr}\left(V_{\theta_0}^{-1}\right)+\mathrm{tr}\left(\mathrm{diag}^{-1}\left(V_{\theta_0}\right)\right)\right\}.
\end{equation*}
The following theorem presents non-asymptotic convergence rates for $\hat{\theta}_n^\ddagger$ and $\hat{v}_n^\ddagger$.

\begin{theorem}
\label{theorem: the convergence rates of the VBBE and the covariance matrix of the VB bridge}
Assume Conditions \ref{condition: local, twice continuously differentiable}--\ref{condition: prior} and the first smoothness condition in \ref{condition: lipschitz} for the variational log-likelihood hold. If $\rho\left(V_{\theta_0}^{-1}\right)\lesssim\sqrt{n}$, then for sufficiently large $r_0^2/p$, we have
\begin{equation*}
n\left(\hat{\theta}_n^\ddagger-\theta^*-\frac{\Delta_{n, \theta_0}}{\sqrt{n}}\right)^{\top}V_{\theta_0}\left(\hat{\theta}_n^\ddagger-\theta^*-\frac{\Delta_{n, \theta_0}}{\sqrt{n}}\right)\lesssim R_n,
\end{equation*}
and
\begin{equation*}
\sum\limits_{i=1}^p\left(\hat{v}_{ni}^{\ddagger} v_{\theta_0,ii}-\log \hat{v}_{ni}^{\ddagger} 
 v_{\theta_0,ii}-1\right)\lesssim R_n,
\end{equation*}
with probability at least $1-5e^{-y}-e^{-cr_0^2}$ for some constant $c>0$.
\end{theorem}

\textit{Proof sketch of Theorem \ref{theorem: the convergence rates of the VBBE and the covariance matrix of the VB bridge}.}
The key is to use Lemma \ref{Lemma: the bracketing result of the logarithm of the denominator of the VB Ideal} to quantify the difference between $\mathrm{KL}\left(q(\theta) \| \pi^*(\theta \mid x)\right)$ and $\mathrm{KL}\left(q(\theta) \| N\left(\theta;\theta^*+\Delta_{n, \theta_0}/\sqrt{n},\mathrm{diag}^{-1}\left(V_{\theta_0}\right)/n\right)\right)$, both considered as functionals of $q(\theta)$.

Utilizing the convergence rates of Theorem \ref{theorem: the convergence rates of the VBBE and the covariance matrix of the VB bridge}, we derive a concentration inequality for the VB bridge $q^{\ddagger}(\theta)$ around the target parameter $\theta^*$.

\begin{theorem}
\label{theorem: the consistency of the VB bridge}
Let $B^c\left(\theta^*, \eta\right)$ denote the complement of the ball of radius $\eta$ centered at $\theta^*$. Under the assumptions of Theorem \ref{theorem: the convergence rates of the VBBE and the covariance matrix of the VB bridge}, it follows that
\begin{equation*}
\int_{B^c\left(\theta^*, \eta\right)} q^{\ddagger}(\theta) \mathrm{d} \theta\lesssim\frac{p+R_n}{\lambda_{\min}\left(V_{\theta_0}\right)n\eta^2+\Delta_{n,\theta_0}^{\top}V_{\theta_0}\Delta_{n,\theta_0}}
\end{equation*}
holds with probability at least $1-5e^{-y}-e^{-cr_0^2}$ for some constant $c>0$.
\end{theorem}

\textit{Proof sketch of Theorem \ref{theorem: the consistency of the VB bridge}.}
The proof consists of three steps: (1) renormalizing $q^{\ddagger}(\theta)$ over $\Theta_0\left(r_0\right)$ as 
\begin{equation*}
q^{\ddagger,\Theta_0(r_0)}(\theta) = \frac{q^{\ddagger}(\theta) I_{\Theta_0(r_0)}(\theta)}{\int q^{\ddagger}(\theta) I_{\Theta_0(r_0)}(\theta) \mathrm{d} \theta},
\end{equation*}
and establishing both upper and lower bounds for $\int q^{\ddagger,\Theta_0\left(r_0\right)}(\theta) M_n(\theta ; x) \mathrm{d} \theta$; (2) deriving a non-asymptotic upper bound for $\int_{B^c\left(\theta^*, \eta\right)}q^{\ddagger,\Theta_0\left(r_0\right)}(\theta)\mathrm{d}\theta$ using Lemma \ref{Lemma: the bracketing result of the logarithm of the denominator of the VB Ideal} and Theorem \ref{theorem: the convergence rates of the VBBE and the covariance matrix of the VB bridge}; and (3) extending the result from $q^{\ddagger,\Theta_0\left(r_0\right)}(\theta)$ to $q^{\ddagger}(\theta)$ by leveraging the thin-tail property of the Gaussian distribution, as described in Lemma \ref{lemma: The Tail Moment Control for the p-Dimensional Normal Distribution}.

We now present a non-asymptotic upper bound for the total variation (TV) distance between the KL minimizer of the rescaled VB bridge and the rescaled reference Gaussian distribution $N\left(\tilde{\theta};\Delta_{n, \theta_0},\mathrm{diag}^{-1}(V_{\theta_0})\right)$, where the rescaled parameter is defined as $\tilde{\theta}=\sqrt{n}\left(\theta-\theta^*\right)$, and the rescaled VB ideal is given by
\begin{equation*}
\pi_{\tilde{\theta}}^*(\tilde{\theta} \mid x)=n^{-\frac{p}{2}}\pi^*\left(\theta^*+\frac{\tilde{\theta}}{\sqrt{n}} \mid x\right).
\end{equation*}

\begin{theorem}
\label{theorem: the asymptotic normality of the VB bridge}
Under the same conditions as Theorem \ref{theorem: the convergence rates of the VBBE and the covariance matrix of the VB bridge}, we have
\begin{equation*}
\left\| \underset{q(\tilde{\theta})\in \mathcal{Q}_\theta}{\arg \min } \operatorname{KL}\left(q(\tilde{\theta}) \| \pi_{\tilde{\theta}}^*(\tilde{\theta}\mid x)\right)-N\left(\tilde{\theta};\Delta_{n, \theta_0}, \mathrm{diag}^{-1}(V_{\theta_0})\right)\right\|^2_{\mathrm{TV}} \lesssim \frac{\max_{i}(v_{\theta_0,ii})R_n}{\lambda_{\min}\left(V_{\theta_0}\right)},
\end{equation*}
with probability at least $1-5e^{-y}-e^{-cr_0^2}$ for some constant $c>0$.
\end{theorem}

\textit{Proof sketch of Theorem \ref{theorem: the asymptotic normality of the VB bridge}.} By applying the Pinsker inequality and leveraging the invariance of the KL divergence under linear transformations within Gaussian distributions, we have
\begin{eqnarray*}
& & \left\| \underset{q(\tilde{\theta})\in \mathcal{Q}_\theta}{\arg \min } \operatorname{KL}\left(q(\tilde{\theta}) \| \pi_{\tilde{\theta}}^*(\tilde{\theta}\mid x)\right)-N\left(\tilde{\theta};\Delta_{n, \theta_0},\mathrm{diag}^{-1}(V_{\theta_0})\right)\right\|^2_{\mathrm{TV}}\\
&\leq& \frac{1}{2}\operatorname{KL}\left(\underset{q(\tilde{\theta})\in \mathcal{Q}_\theta}{\arg \min } \operatorname{KL}\left(q(\tilde{\theta}) \| \pi_{\tilde{\theta}}^*(\tilde{\theta}\mid x)\right)\|N\left(\tilde{\theta};\Delta_{n, \theta_0},\mathrm{diag}^{-1}(V_{\theta_0})\right)\right)\\
&=& \frac{1}{2}\operatorname{KL}\left(\underset{q(\theta) \in \mathcal{Q}_\theta}{\arg \min } \operatorname{KL}\left(q(\theta) \| \pi^*(\theta\mid x)\right)\|N\left(\tilde{\theta};\theta^*+\frac{\Delta_{n, \theta_0}}{\sqrt{n}},\frac{\mathrm{diag}^{-1}(V_{\theta_0})}{n}\right)\right).
\end{eqnarray*}
We then utilize the convergence rates established in Theorem \ref{theorem: the convergence rates of the VBBE and the covariance matrix of the VB bridge} to obtain the bound.

\section{Non-asymptotic properties of the VB posterior}
\label{sec: non-asymptotic properties of the VB posterior}
Lemma \ref{lemma: the gap between two functionals} indicates that the VB posterior inherits similar non-asymptotic properties to the VB bridge under mild conditions. In this section, we establish the non-asymptotic properties of the VB posterior, including the non-asymptotic variational Bernstein--von Mises theorem.

Define the VB posterior $q^*(\theta) \sim N\left(\theta;\hat{\theta}_n^*, \operatorname{diag}\left(\hat{v}^*_n\right)/n\right)$, where  $\hat{v}_n^*=\left(\hat{v}_{n1}^*, \ldots, \hat{v}_{np}^*\right)^\top$. The following theorem presents non-asymptotic convergence rates for $\hat{\theta}_n^*$ and $\hat{v}_n^*$.
 
\begin{theorem}
\label{theorem: the convergence rates of the distribution parameters of the VB posterior}
Assume Conditions \ref{condition: local, twice continuously differentiable}–\ref{condition: lipschitz} hold. If $\rho\left(V_{\theta_0}^{-1}\right)\lesssim\sqrt{n}$, then for sufficiently large $r_0^2/p$, we have
\begin{equation*}
\begin{aligned}
&n\left(\hat{\theta}_n^*-\theta^*-\frac{\Delta_{n, \theta_0}}{\sqrt{n}}\right)^{\top}V_{\theta_0}\left(\hat{\theta}_n^*-\theta^*-\frac{\Delta_{n, \theta_0}}{\sqrt{n}}\right)\\\lesssim&R_n+\frac{p}{\min v_{\theta_0,ii}}\left\{\left(\frac{p}{n\min v_{\theta_0,ii}}\right)^{\frac{s_1}{2}}+\left(\frac{p}{n\min v_{\theta_0,ii}}\right)^{\frac{s_2}{2}}\right\},
\end{aligned}
\end{equation*}
and
\begin{equation*}
\sum\limits_{i=1}^p\left(\hat{v}_{ni}^*v_{\theta_0,ii}-\log\hat{v}_{ni}^*v_{\theta_0,ii}-1\right)\\
\lesssim R_n+\frac{p}{\min v_{\theta_0,ii}}\left\{\left(\frac{p}{n\min v_{\theta_0,ii}}\right)^{\frac{s_1}{2}}+\left(\frac{p}{n\min v_{\theta_0,ii}}\right)^{\frac{s_2}{2}}\right\},
\end{equation*}
with probability at least $1-5e^{-y}-e^{-cp}-e^{-cr_0^2}$ for some constant $c>0$.
\end{theorem}

\textit{Proof sketch of Theorem \ref{theorem: the convergence rates of the distribution parameters of the VB posterior}.} We use the techniques from the proof of Theorem \ref{theorem: the convergence rates of the VBBE and the covariance matrix of the VB bridge} in Appendix \ref{Proof of Theorem 1} to evaluate the following term
\begin{equation*}
\mathrm{KL}\left(q^*(\theta) \| \pi^*(\theta \mid x)\right)-\mathrm{KL}\left(N\left(\theta;\theta^*+\frac{\Delta_{n, \theta_0}}{\sqrt{n}},\frac{\operatorname{diag}^{-1}\left(V_{\theta_0}\right)}{n}\right) \middle\| \pi^*(\theta \mid x)\right).
\end{equation*}
We then combine Lemma \ref{lemma: the gap between two functionals} and Theorem \ref{theorem: the convergence rates of the VBBE and the covariance matrix of the VB bridge} to derive the convergence rates.

\begin{corollary}
\label{corollary of theorem 4}
Assume Conditions \ref{condition: local, twice continuously differentiable}--\ref{condition: lipschitz} hold. Additionally, suppose that (i) $r_0^2\asymp p$, and $r_0^2 / p$ can be sufficiently large, (ii) $\rho\left(V_{\theta_0}^{-1}\right)\lesssim 1$, (iii) $1 / \min v_{\theta_0,ii} \lesssim 1$, and (iv) $\left\|\nabla\log p\left(\theta^*\right)\right\|^4\lesssim p^3$. Then,
\begin{equation*}
\|\hat{\theta}_n^*-\theta^*\| = O_p\left(\sqrt{\frac{p}{n}}\right).
\end{equation*}
\end{corollary}

\textit{Proof Sketch of Corollary \ref{corollary of theorem 4}.} The result follows from Theorem \ref{theorem: the convergence rates of the distribution parameters of the VB posterior} by analyzing the order of the upper bound.

Theorem \ref{theorem: the convergence rates of the distribution parameters of the VB posterior} extends the consistency of the VB estimator, as established by \cite{wang2019frequentist}, to non-asymptotic settings. Corollary \ref{corollary of theorem 4} further shows that, under mild conditions, the VB estimator remains consistent as long as $p/n$ approaches zero. The convergence rate $\sqrt{p/n}$ aligns with those derived for maximum likelihood estimators and the corresponding Bayes estimators in general parametric models \citep{spokoiny2014bernsteinvonmises,he2000parameters}.

\begin{corollary}
\label{Corollary 2}
Under the same conditions as Corollary \ref{corollary of theorem 4}, if $\min\left(s_1,s_2\right) \geq 1$ and $p^3=o(n)$, then
\begin{equation*}
\sqrt{n}\alpha^{\top}\left(\hat{\theta}_n^*-\theta^*\right)/\sigma_\alpha \overset{\mathcal{L}}{\longrightarrow}  N(0,1),
\end{equation*}
for any $\alpha \in S_p=\left\{\alpha \in \mathbb{R}^p:\|\alpha\|=1\right\}$, where $\sigma_\alpha^2=\alpha^{\top}V_{\theta_0}^{-1}\operatorname{Var}_{\theta_0}\left(\nabla m\left(\theta^*; x\right)\right)V_{\theta_0}^{-1}\alpha$.
\end{corollary}

\textit{Proof Sketch of Corollary \ref{Corollary 2}.} This result is directly implied by the proof of Corollary \ref{corollary of theorem 4}.

Corollary \ref{Corollary 2} demonstrates that, under mild conditions, the VB estimator is asymptotically normal if $p^3/n$ converges to zero. Again, the result is consistent with those derived for maximum likelihood estimators and corresponding Bayes estimators in general parametric models \citep{he2000parameters,spokoiny2014bernsteinvonmises}.

We now present the central results of this paper: the non-asymptotic variational Bernstein--von Mises theorem. This theorem characterizes the concentration of the VB posterior around the target parameter $\theta^*$ and bounds the TV distance between the rescaled VB posterior and the rescaled reference Gaussian distribution.

\begin{theorem}[Non-asymptotic Variational Bernstein--von Mises Theorem]
\label{theorem: the non-asymptotic variational Bernstein--von Mises theorem}
Under the conditions of Theorem \ref{theorem: the convergence rates of the distribution parameters of the VB posterior}, with probability at least $1-5e^{-y}-e^{-cp}-e^{-cr_0^2}$ for some constant $c>0$, the following results hold:
\begin{itemize}
    \item[(1)]
    \begin{equation*}
\int_{B^c\left(\theta^*, \eta\right)} q^*(\theta)\mathrm{d} \theta\lesssim\frac{p+R_n+\frac{p}{\min v_{\theta_0,ii}}\left\{\left(\frac{p}{n\min v_{\theta_0,ii}}\right)^{\frac{s_1}{2}}+\left(\frac{p}{n\min v_{\theta_0,ii}}\right)^{\frac{s_2}{2}}\right\}}{\lambda_{\min}\left(V_{\theta_0}\right)n\eta^2+\Delta_{n, \theta_0}^{\top}V_{\theta_0}\Delta_{n, \theta_0}},
\end{equation*}
\item[(2)]
\begin{equation*}
\begin{aligned}
&\left\|q_{\tilde{\theta}}^*(\tilde{\theta})-N\left(\tilde{\theta};\Delta_{n, \theta_0},\mathrm{diag}^{-1}(V_{\theta_0})\right)\right\|^2_{\mathrm{TV}} 
\\\lesssim&\frac{\max\nolimits_{i}v_{\theta_0,ii}}{\lambda_{\min}\left(V_{\theta_0}\right)}\left[R_n+\frac{p}{\min v_{\theta_0,ii}}\left\{\left(\frac{p}{n\min v_{\theta_0,ii}}\right)^{\frac{s_1}{2}}+\left(\frac{p}{n\min v_{\theta_0,ii}}\right)^{\frac{s_2}{2}}\right\}\right].
\end{aligned}
\end{equation*}
\end{itemize}
\end{theorem}

\textit{Proof Sketch of Theorem \ref{theorem: the non-asymptotic variational Bernstein--von Mises theorem}.} 
We employ techniques similar to those used in the proofs of Theorems \ref{theorem: the consistency of the VB bridge} and \ref{theorem: the asymptotic normality of the VB bridge}, and leverage the results from Theorems \ref{theorem: the convergence rates of the VBBE and the covariance matrix of the VB bridge}-\ref{theorem: the convergence rates of the distribution parameters of the VB posterior}.

\begin{corollary}
\label{corollary of theorem 5}
Under the conditions of Corollary \ref{corollary of theorem 4},
\begin{equation*}
\int_{B^c\left(\theta^*, \eta\right)} q^*(\theta) \mathrm{d} \theta=O_p\left(\frac{p}{n\eta^2}\right).
\end{equation*}
Moreover, if $\max_i(v_{\theta_0, ii})\lesssim 1$ and $\min\left(s_1,s_2\right) \geq 1$, we have
\begin{equation*}
\left\|q_{\tilde{\theta}}^*(\tilde{\theta})-N\left(\tilde{\theta};\Delta_{n, \theta_0},\mathrm{diag}^{-1}(V_{\theta_0})\right)\right\|^2_{\mathrm{TV}}=O_p\left(\sqrt{\frac{p^3}{n}}\right).
\end{equation*}
\end{corollary}

Theorem \ref{theorem: the non-asymptotic variational Bernstein--von Mises theorem} extends the variational Bernstein--von Mises theorem of \cite{wang2019frequentist} to non-asymptotic settings. The critical dimensions $p = o\left(n\right)$ and $p^3 = o\left(n\right)$, as established in Corollary \ref{corollary of theorem 5}, align with the conditions for the traditional Bernstein--von Mises Theorem with increasing dimensions \citep{spokoiny2014bernsteinvonmises}.

\begin{remark}
\label{remark: the recovery}
We have derived the non-asymptotic properties of the VB posterior and its mean. When the expected variational log-likelihood is uniquely maximized at $\theta_0$, Theorem \ref{theorem: the non-asymptotic variational Bernstein--von Mises theorem} and Corollary \ref{corollary of theorem 5} suggest that the VB posterior can accurately recover the true parameter and marginal variance within the Gaussian mean-field variational family, consistent with the fixed-dimensional results in \cite{wang2019frequentist}. Additionally, since
\begin{equation*}
\mathbb{H}\left(N\left(\tilde{\theta};\Delta_{n, \theta_0},\mathrm{diag}^{-1}(V_{\theta_0})\right)\right) \leq \mathbb{H}\left(N\left(\tilde{\theta};\Delta_{n, \theta_0},V_{\theta_0}^{-1}\right)\right),
\end{equation*}
where $\mathbb{H}$ denotes the entropy functional, Theorem \ref{theorem: the non-asymptotic variational Bernstein--von Mises theorem} and Corollary \ref{corollary of theorem 5} corroborate the findings in \cite{wang2019frequentist} that the VB posterior derived from the Gaussian mean-field family tends to be under-dispersed, as discussed in their Lemma 8. This limitation may be alleviated by expanding the variational family to better capture dependencies among latent variables.
\end{remark}

\section{An application to multivariate Gaussian mixture models}
\label{sec: applications}
\cite{wang2019frequentist} studied the VB posterior and the VB estimator for univariate Gaussian mixture models. In this section, we extend this framework to Bayesian mixtures of multivariate Gaussians with increasing dimensions. Specifically, the model takes the following form:
\begin{equation}
\label{abstract GMM}
\begin{aligned}
\mu_k & \sim p_0, & & k=1, \ldots, K, \\
c_i & \sim \text{Categorical}\left(1/K,\ldots,1/K\right), & & i=1, \ldots, n, \\
x_i \mid c_i, \mu & \sim N\left(\mu c_i,I_p\right), & & i=1, \ldots, n,
\end{aligned}
\end{equation}
where $\mu_k\in\mathbb{R}^p$ for $k=1,\ldots, K$ are the $p$-dimensional cluster means drawn independently from a prior $p_0$; $c_i$ for $i=1,\ldots,n$ are the cluster assignments, drawn independently from a categorical distribution with uniform weights; $\mu=\left(\mu_1,\ldots,\mu_K\right)$ denotes the matrix of cluster means; and $x_i$, for $i=1,\ldots, n$, are the observed data points, each drawn from a normal distribution with mean vector $\mu c_i$ and identity covariance matrix. Here $\mu$ is a $p \times K$ global latent matrix and $\{c_1, \ldots, c_n\}$ are local latent variables. We aim to infer $\mu$, assuming that the data are generated from $p(x; \mu^0)=\int p(x, c; \mu^0)\mathrm{d}c$, where $\mu^0$ represents the true value of $\mu$.

The following proposition provides the expression of the variational log-likelihood $m(\mu;x_i)$ for the $i$-th observation $x_i$. 

\begin{proposition}
\label{m=logp}
Assume that model \eqref{abstract GMM} holds. Then, the variational log-likelihood coincides with the data log-likelihood:
\begin{equation*}
m(\mu;x_i)=-\log K-\frac{p}{2} \log (2 \pi)+\log\left(\sum\limits_{k=1}^K \exp \left(-\frac{1}{2}\left\|x_i-\mu_k\right\|^2\right)\right)=\log p(x_i;\mu).
\end{equation*}
\end{proposition}

Building on Proposition \ref{m=logp}, the next proposition establishes that the expected variational log-likelihood is maximized precisely at the true parameter $\mu^0$.

\begin{proposition}
\label{mustar=mu0}
Define $\mu^* = \underset{\mu \in \Theta}{\operatorname{argmax}}\;\mathbb{E}_{\mu^0}\{\sum^n_{i=1}m(\mu;x_i)\}$. Assume that the covariance matrix $V_{\mu^0}=\operatorname{Var}_{\mu^0}\left(\nabla \log p(x_i;\mu^0)\right)$ is positive definite, then $\mu^* = \mu^0$ up to a permutation. 
\end{proposition}

\begin{proposition}
\label{Proposition: the (variational) information matrix}
Suppose the following condition holds:
\begin{equation*}
\sum_{k_1< k_2}\frac{\left\|\mu_{k_1}^0-\mu_{k_2}^0\right\|^2+4p}{4\exp\left(\frac{1}{8}\left\|\mu_{k_1}^0-\mu_{k_2}^0\right\|^2\right)}\leq C,
\end{equation*}
for some independent constant $C\in\left(0,1\right)$. Then, for any $1\leq r\leq Kp$, we have
\begin{equation*}
\frac{1-C}{K}\leq\lambda_{\min}\left(V_{\mu^0}\right)\leq \left(V_{\mu^0}\right)_{r,r}\leq\lambda_{\max}\left(V_{\mu^0}\right)\leq\frac{1}{K},
\end{equation*}
where $\left(V_{\mu^0}\right)_{r,r}$ represents the $r$-th diagonal element of $V_{\mu^0}$.
\end{proposition}

Proposition \ref{Proposition: the (variational) information matrix} establishes an important property of $V_{\mu^0}$. It shows that when the cluster means are sufficiently separated, the variational information matrix can be approximated by $I_{Kp}/K$. In more complex scenarios, such as those involving unequal component weights or non-identity covariance matrices, we conjecture that $V_{\mu^0}$ behaves similarly, provided the weights are not excessively imbalanced and the covariance matrices are well-conditioned.

To establish the theoretical properties of the VB posterior and the VB estimator under model \eqref{abstract GMM}, it suffices to verify Conditions \ref{condition: local, twice continuously differentiable}–\ref{condition: lipschitz}. However, most of these conditions are difficult to verify due to the presence of the log-sum term in the variational log-likelihood. For a discussion of some of these conditions, refer to Appendix \ref{Some Verifications of Design Conditions}.

In the remainder of this section, we use simulated data to evaluate the performance of the VB estimator. Specifically, we consider the Gaussian mixture model in \eqref{abstract GMM} with $K=2$ and $p_0$ being $N(0,\sigma^2I_p)$, where $\sigma^2$ denotes the prior variance. Furthermore, we set $\mu_1^0=-w1_p$ and $\mu_2^0=w1_p$, where $w>0$ quantifies the separation between the two true means, and serves as a measure of the component dispersion in this model. 

Let $m_k$ and $D_k$ be the variational parameters for the global latent variable $\mu_k, k = 1, 2$, where $m_k$ is a $p$-dimensional vector and $D_k$ is $p$ × $p$ positive definite diagonal matrix. Additionally, let $\varphi_i$ be the variational parameters for the local latent variable $c_i$, $i = 1,\ldots,n$. The corresponding ELBO, denoted as $\operatorname{ELBO}(m,D,\varphi)$, is given by
\begin{equation*}
\begin{aligned}
&\operatorname{ELBO}(m,D,\varphi)\\
=&\frac{1}{2} \left\{\log\det\left(D_1\right)+\log\det\left(D_2\right)\right\} - \frac{1}{2\sigma^2} \left\{\mathrm{tr}(D_1)+\mathrm{tr}(D_2)+\left\|m_1\right\|^2+\left\|m_2\right\|^2\right\}\\& +p-p\log\sigma^2-\frac{1}{2}\sum\limits_{i=1}^n\left[\varphi_{i1}\left\{\left\|x_i-m_1\right\|^2+\mathrm{tr}(D_1)\right\}+\varphi_{i2}\left\{\left\|x_i-m_2\right\|^2+\mathrm{tr}(D_2)\right\}\right]\\&-\sum\limits_{i=1}^n\left(\varphi_{i1}\log\varphi_{i1}+\varphi_{i2}\log\varphi_{i2}\right)-n\log2-\frac{np}{2}\log(2\pi).
\end{aligned}
\end{equation*}
Maximizing $\operatorname{ELBO}\left(m, D,\varphi\right)$ with respect to the variational parameters can be performed using a coordinate ascent variational inference (CAVI) algorithm, as detailed in Algorithm \ref{CAVI for the p-dimensional 2-component GMM}.

\begin{algorithm}
\caption{CAVI for the $p$-dimensional 2-component Gaussian Mixture Model}
\label{CAVI for the p-dimensional 2-component GMM}
\begin{algorithmic}[ht]
\Input $\left\{x_1,\ldots,x_n\right\}$ and $\sigma^2$.
\State Initialize $\left(m^{(0)},D^{(0)},\varphi^{(0)}\right)$, and $t=0$;
\While{$\operatorname{ELBO}(m,D,\varphi)$ not converged}
    \State $t=t+1$;
    \For{$i \gets 1$ \textbf{to} $n$}
    \State $\varphi_{i1}^{(t)} \gets \frac{\exp\left(x_i^{\top}m_1^{(t-1)}-\frac{1}{2}\left[\left\|m_1^{(t-1)}\right\|^2+\mathrm{tr}\left(D_1^{(t-1)}\right)\right]\right)}{\exp\left(x_i^{\top}m_1^{(t-1)} -\frac{1}{2} \left[\left\|m_1^{(t-1)}\right\|^2+\mathrm{tr}\left(D_1^{(t-1)}\right)\right] \right) + \exp\left(x_i^{\top}m_2^{(t-1)}-\frac{1}{2} \left[\left\|m_2^{(t-1)}\right\|^2+\mathrm{tr}\left(D_1^{(t-1)}\right)\right]\right)}$;
    \State $\varphi_{i2}^{(t)}\gets1-\varphi_{i1}^{(t)}$;
    \EndFor
        \State $m_1^{(t)}\gets\frac{\sum\limits_{i=1}^n\varphi_{i1}^{(t)}x_i}{\sum\limits_{i=1}^n\varphi_{i1}^{(t)}+\frac{1}{\sigma^2}},\quad m_2^{(t)}\gets\frac{\sum\limits_{i=1}^n\varphi_{i2}^{(t)}x_i}{\sum\limits_{i=1}^n\varphi_{i2}^{(t)}+\frac{1}{\sigma^2}}$;
    \State $D_1^{(t)}\gets\frac{I_p}{\sum\limits_{i=1}^n\varphi_{1i}^{(t)}+\frac{1}{\sigma^2}},\quad D_2^{(t)}\gets\frac{I_p}{\sum\limits_{i=1}^n\varphi_{2i}^{(t)}+\frac{1}{\sigma^2}}$;
\EndWhile
\Output $\left(m^*,D^*,\varphi^*\right)=\left(m^{(t)},D^{(t)},\varphi^{(t)}\right)$.
\end{algorithmic}
\end{algorithm}

To evaluate the estimation accuracy, we employ the mean squared error (MSE), defined as $\left(\left\|m_1^*-\mu_1^0\right\|^2+\left\|m_2^*-\mu_2^0\right\|^2\right)/p$. We vary $n\in\left\{50, 200, 1000, 5000\right\}, p\in\left\{2, 5, 10, 50\right\}$, and $\left(\sigma^2,w\right)\in\left\{(1, 10), (1, 50), (25, 10), (25, 50)\right\}$. The results, based on 100 simulation runs, are reported in Table \ref{table: Mean and standard deviation of the MSE}. As expected, the performance improves with increasing sample size or prior variance.

\begin{table}[ht]
\centering
\renewcommand{\arraystretch}{1} 
\setlength{\tabcolsep}{8pt} 
\caption{Mean and standard deviation (in parentheses) of the MSE, computed over 100 simulation runs, for various configurations of $n, p$, and $(\sigma^2, w)$.}
\label{table: Mean and standard deviation of the MSE}
\resizebox{\linewidth}{!}{
\begin{tabular}{cc*{4}{c}}
\toprule
$n$ & $p$ & $(\sigma^2,w)=(1, 10)$ & $(\sigma^2,w)=(1, 50)$ & $(\sigma^2,w)=(25, 10)$ & $(\sigma^2,w)=(25, 50)$ \\
\midrule
50 & 2 & 0.4115(0.1545)&7.8830(1.0260)&0.0905(0.0620)&0.0960(0.0740)
 \\
 &  5   & 0.3894(0.1106)&7.8304(0.6798)&0.0824(0.0372)&0.0990(0.0396)
 \\
 &  10   & 0.3960(0.0726)&7.9920(0.7303)&0.0863(0.0323)&0.0991(0.0328)

 \\
 &  50   & 0.3922(0.0513)&7.9346(0.5741)&0.0834(0.0134)&0.0938(0.0138)

 \\
200 & 2 & 0.0383(0.0213)&0.5230(0.1145)&0.0211(0.0146)&0.0208(0.0128)
 \\
 &  5   & 0.0408(0.0154)&0.5104(0.0580)&0.0198(0.0094)&0.0202(0.0090)
\\
 &  10   & 0.0380(0.0099)&0.5109(0.0441)&0.0190(0.0060)&0.0217(0.0068)

 \\
 &  50   & 0.0392(0.0052)&0.5138(0.0235)&0.0199(0.0029)&0.0206(0.0026)

 \\
 1000 & 2 & 0.0051(0.0035)&0.0226(0.0084)&0.0037(0.0027)&0.0042(0.0030)
 \\
 &  5   & 0.0046(0.0019)&0.0250(0.0059)&0.0041(0.0019)&0.0043(0.0173)
 \\
 &  10   & 0.0048(0.0017)&0.0240(0.0043)&0.0042(0.0014)&0.0040(0.0013)

 \\
 &  50   & 0.0047(0.0007)&0.0239(0.0019)&0.0041(0.0006)&0.0040(0.0005)

 \\
 5000 & 2 & 0.0009(0.0007)&0.0016(0.0009)&0.0007(0.0006)&0.0009(0.0006)
 \\
 &  5   & 0.0008(0.0003)&0.0016(0.0006)&0.0008(0.0003)&0.0008(0.0004)
 \\
 &  10   & 0.0008(0.0003)&0.0016(0.0004)&0.0008(0.0002)&0.0008(0.0003)

 \\
 &  50   & 0.0008(0.0001)&0.0016(0.0002)&0.0008(0.0001)&0.0008(0.0001)
\\
\bottomrule
\end{tabular}
}
\end{table}

\section{Conclusion}
\label{sec: conclusion}
For a broad class of parametric models with latent variables, we have established the concentration of the VB posterior and provided non-asymptotic bounds on the TV distance between the rescaled VB posterior and a normal distribution. We have also characterized the consistency and asymptotic normality of the corresponding VB estimator. Our results demonstrate that the critical dimensions governing the behavior of the VB posterior and estimator are consistent with those known for maximum likelihood estimators and corresponding Bayes estimators in general parametric models. In doing so, we generalize the results of \cite{wang2019frequentist} to settings with diverging dimensions, thereby resolving a question raised in their discussion. Of independent interest are the theoretical properties of the VB bridge, for which we have connected with the VB posterior.

There are several avenues for future exploration. First, beyond the KL divergence, we are currently investigating alternative divergence measures, such as $\alpha$-divergence and $\gamma$-divergence, to more comprehensively study the non-asymptotic properties of the VB approach. Second, extending our theoretical framework to the $\alpha$-fractional posterior and $\alpha$-variational inference is an interesting direction, with work in this area currently underway. Third, it remains important yet challenging to address scenarios where the parameter dimension diverges more rapidly than, or even exceeds, the sample size. In such cases, additional model restrictions, such as sparsity, are needed.

\bibliographystyle{apalike}

\bibliography{ref1}
\end{document}


\renewcommand{\qedsymbol}{}

\def\spacingset#1{\renewcommand{\baselinestretch}%
	{#1}\small\normalsize} \spacingset{1}

	
	\if\blind
	{
		\title{\bf Factor augmented inverse regression and its application to microbiome data analysis}
		\author
		{Daolin Pang$^{1}$, Hongyu Zhao$ ^{2,3} $, Tao Wang$^{1,3,4,5\ast}$\\
			\\
			\normalsize{$^{1}$Department of Bioinformatics and Biostatistics, Shanghai Jiao Tong University}\\
			\normalsize{$^{2}$Department of Biostatistics, Yale University}\\
			\normalsize{$^{3}$SJTU-Yale Joint Center for Biostatistics and Data Science, Shanghai Jiao Tong University}\\
			\normalsize{$^{4}$Department of Statistics, School of Mathematical Sciences, Shanghai Jiao Tong University}\\
			\normalsize{$^{5}$MoE Key Lab of Artificial Intelligence, AI Institute, Shanghai Jiao Tong University}\\
			\\
			\normalsize{$^\ast$E-mail: neowangtao@sjtu.edu.cn}
		}
		\date{}
		\maketitle
		
	} \fi
	
	\title{\bf Appendix to ``Variational Bernstein–von Mises theorem with increasing parameter dimension"}
	\date{}
	\maketitle
	
	\bigskip
	
	\spacingset{1.9}

\section{Proof of Lemma \ref{LAN}}
\label{Proof of Lemma 1}

Conditions \ref{condition: local, twice continuously differentiable}--\ref{condition: global, exponential moment, variational} meet the assumptions of Proposition 3.2 in \cite{spokoiny2014bernsteinvonmises}. Consequently, we can establish a local quadratic approximation of the variational log-likelihood \( M_n(\theta; x) \) using the following bracketing bound. Specifically, applying Proposition 3.2 and Equation (3.5) from \cite{spokoiny2014bernsteinvonmises} yields that, with probability at least $1-e^{-y}$, the bound holds:
\begin{eqnarray}
\label{LAN expansion}
& & \sup _{\theta \in \Theta_0\left(r_0\right)}\left|M_n(\theta ; x)-M_n\left(\theta^* ; x\right)-\left(\theta-\theta^*\right)^{\top} \nabla M_n\left(\theta^* ; x\right)+ \frac{1}{2}\left\|D_0\left(\theta-\theta^*\right)\right\|^2 \right| \nonumber \\
&\leq& r_0^2\left\{\delta_1(r_0)+6\tau_1 z_{\mathbb{H}}(y)\tau_3\right\}.
\end{eqnarray}

Next, we perform a change of variables by setting $h=\sqrt{n}\left(\theta-\theta^*\right)$ and $\Delta_{n, \theta_0}=V_{\theta_0}^{-1}\nabla M_n\left(\theta^* ; x\right)/\sqrt{n}$. With this substitution, we set $\delta_n$ in \cite{wang2019frequentist} to $\sqrt{n}$, aligning with the framework for i.i.d. cases presented in Chapter 7 of \cite{Vaart_1998}. Recall that $V_{\theta_0}=D_0^2/n$, therefore, the conditions $\left\{\left\|V_{\theta_0}^{\frac{1}{2}} h\right\| \leq r_0\right\}$ and $\left\{\left\|D_0\left(\theta-\theta^*\right)\right\| \leq r_0\right\}$ are equivalent. Substituting these variable transformations into \eqref{LAN expansion} yields the desired result. 

\section{Proof of Lemma \ref{Lemma: the bracketing result of the logarithm of the denominator of the VB Ideal}}
\label{Proof of Lemma 2}

We decompose the normalizing constant $\int p(\theta) \exp \left\{M_n(\theta ; x)\right\} \mathrm{d} \theta$ as follows:
\begin{eqnarray}
\label{The division of the denominator of the VB ideal}
& & \int p(\theta) \exp \left\{M_n(\theta ; x)\right\} \mathrm{d} \theta \nonumber\\
&=& \int_{\Theta_0\left(r_0\right)} p(\theta) \exp \left\{M_n(\theta ; x)\right\} \mathrm{d} \theta+\int_{\Theta_0^c\left(r_0\right)} p(\theta) \exp \left\{M_n(\theta ; x)\right\} \mathrm{d} \theta.
\end{eqnarray}
For each term on the right-hand side of Equation \eqref{The division of the denominator of the VB ideal}, we establish both upper and lower bounds. The overall derivation is organized into four steps.

\textit{Step 1.} We establish an upper bound for $\int_{\Theta_0\left(r_0\right)} p(\theta) \exp \left\{M_n(\theta; x)\right\} \mathrm{d} \theta$. 

Recall that $\Delta(r_0,y)=r_0^2\left\{\delta_1(r_0)+6\tau_1 z_{\mathbb{H}}(y)\tau_3\right\}$. Under Conditions \ref{condition: local, twice continuously differentiable}--\ref{condition: global, exponential moment, variational}, it follows from \eqref{LAN expansion} that
\begin{eqnarray*}
& & \int_{\Theta_0(r_0)} p(\theta) \exp \left\{M_n(\theta ; x)\right\} \mathrm{d} \theta\\
&\leq& \int_{\Theta_0(r_0)} p(\theta) \exp \left\{M_n(\theta^*;x)+\left(\theta-\theta^*\right)^{\top} \nabla M_n\left(\theta^* ; x\right)- \frac{1}{2}\left\|D_0\left(\theta-\theta^*\right)\right\|^2 +\Delta(r_0,y)\right\} \mathrm{d} \theta\\
&=&\exp\left\{M_n(\theta^*;x)+\Delta(r_0,y)\right\} \int_{\Theta_0(r_0)} p(\theta) \exp \left\{\left(\theta-\theta^*\right)^{\top} \nabla M_n\left(\theta^* ; x\right)-\frac{1}{2}\left\|D_0\left(\theta-\theta^*\right)\right\|^2 \right\}\mathrm{d} \theta
\end{eqnarray*}
holds with probability at least $1 - e^{-y}$. Note that 
\begin{equation}
\label{Second order control for the prior}
\left|\log p(\theta)-\log p\left(\theta^*\right)-\left(\theta-\theta^*\right)^{\top}\nabla\log p\left(\theta^*\right)\right|\leq \frac{\tau_5}{2}\left\|\theta-\theta^*\right\|^2
\end{equation}
under Condition \ref{condition: prior}. Therefore,
\begin{eqnarray*}
& & \int_{\Theta_0(r_0)} p(\theta) \exp \left\{\left(\theta-\theta^*\right)^{\top} \nabla M_n\left(\theta^* ; x\right)- \frac{1}{2}\left\|D_0\left(\theta-\theta^*\right)\right\|^2 \right\} \mathrm{d} \theta\\ &\leq& p\left(\theta^*\right)\int_{\Theta_0(r_0)} \exp \left\{\left(\theta-\theta^*\right)^{\top}\nabla\left[ M_n\left(\theta^* ; x\right)+\log p\left(\theta^*\right)\right]-\frac{1}{2}\left\| \left(D_0^2-\tau_5I_p\right)^{\frac{1}{2}} \left(\theta-\theta^*\right)\right\|^2\right\} \mathrm{d} \theta.
\end{eqnarray*}
Additionally, we know that
\begin{eqnarray*}
& & \int_{\Theta_0(r_0)} \exp \left\{\left(\theta-\theta^*\right)^{\top} \nabla\left[ M_n\left(\theta^* ; x\right)+\log p\left(\theta^*\right)\right] - \frac{1}{2}\left\| \left(D_0^2-\tau_5I_p\right)^{\frac{1}{2}} \left(\theta-\theta^*\right)\right\|^2\right\} \mathrm{d} \theta\\
&=&\int \exp \left\{\left(\theta-\theta^*\right)^{\top} \nabla\left[ M_n\left(\theta^* ; x\right)+\log p\left(\theta^*\right)\right] -\frac{1}{2}\left\| \left(D_0^2-\tau_5I_p\right)^{\frac{1}{2}} \left(\theta-\theta^*\right)\right\|^2\right\} \mathrm{d} \theta\\
& & -\int_{\Theta_0^c(r_0)} \exp \left\{\left(\theta-\theta^*\right)^{\top}\nabla\left[ M_n\left(\theta^* ; x\right) + \log p\left(\theta^*\right)\right] -\frac{1}{2}\left\| \left(D_0^2-\tau_5I_p\right)^{\frac{1}{2}} \left(\theta-\theta^*\right)\right\|^2\right\} \mathrm{d} \theta,
\end{eqnarray*}
and
\begin{eqnarray*}
& & \int \exp \left\{\left(\theta-\theta^*\right)^{\top}\nabla\left[ M_n\left(\theta^* ; x\right)+\log p\left(\theta^*\right)\right]-\frac{1}{2} \left\| \left(D_0^2-\tau_5I_p\right)^{\frac{1}{2}} \left(\theta-\theta^*\right)\right\|^2 \right\} \mathrm{d} \theta\\
&=&\frac{(2\pi)^{\frac{p}{2}}}{\det\left(D_0^2-\tau_5I_p\right)^{\frac{1}{2}}} \exp\left\{\frac{1}{2}\left\| \left(D_0^2-\tau_5I_p\right)^{-\frac{1}{2}} \nabla\left[M_n(\theta^* ; x) + \log p(\theta^*)\right]\right\|^2\right\}.
\end{eqnarray*}
Consequently,
\begin{eqnarray}
\label{Upper bound of the central part of the denominator of the VB ideal}
& & \int_{\Theta_0(r_0)} p(\theta) \exp \left\{M_n(\theta ; x)\right\} \mathrm{d} \theta \nonumber\\
&\leq& \exp \left\{M_n(\theta^*;x)+\Delta(r_0,y)+\log p\left(\theta^*\right)\right\} \nonumber\\
& & \times \left[\frac{(2\pi)^{\frac{p}{2}}}{\det\left(D_0^2-\tau_5I_p\right)^{\frac{1}{2}}} \exp\left\{\frac{1}{2}\left\| \left(D_0^2-\tau_5I_p\right)^{-\frac{1}{2}} \nabla\left[M_n(\theta^* ; x) + \log p(\theta^*)\right]\right\|^2\right\} \right. \\
& & \left.-\int_{\Theta_0^c(r_0)} \exp \left\{\left(\theta-\theta^*\right)^{\top}\nabla\left[ M_n\left(\theta^* ; x\right)+\log p\left(\theta^*\right)\right]-\frac{1}{2} \left\| \left(D_0^2-\tau_5I_p\right)^{\frac{1}{2}} \left(\theta-\theta^*\right)\right\|^2
\right\} \mathrm{d} \theta\right]\nonumber
\end{eqnarray}
holds with probability at least $1 - e^{-y}$.

\textit{Step 2.} We derive an upper bound for $\int_{\Theta_0^c(r_0)} p(\theta) \exp \left\{M_n(\theta; x)\right\} \mathrm{d} \theta$. 

To aid the derivation, we impose the following mild assumption on the prior distribution:
\begin{eqnarray*}
& & \frac{\int_{\Theta_0^c(r_0)} p(\theta) \exp \left\{M_n(\theta; x)\right\} \mathrm{d} \theta}{\int_{\Theta_0(r_0)}p\left(\theta^*\right)\exp \left\{M_n(\theta ; x)+\left(\theta-\theta^*\right)^{\top}\nabla\log p\left(\theta^*\right)+\frac{\tau_5}{2}\left\|\theta-\theta^*\right\|^2\right\} \mathrm{d} \theta}\\
&\lesssim&\frac{\int_{\Theta_0^c(r_0)} \exp \left\{M_n(\theta; x)\right\} \mathrm{d} \theta}{\int_{\Theta_0(r_0)} \exp \left\{M_n(\theta; x)\right\} \mathrm{d} \theta}.
\end{eqnarray*}
This inequality holds, for instance, when
\begin{equation*}
\max_{\theta \in \Theta_0^c(r_0)} p(\theta) \lesssim \min_{\theta \in \Theta_0(r_0)} p\left(\theta^*\right)\exp \left\{\left(\theta-\theta^*\right)^{\top}\nabla\log p\left(\theta^*\right)+\frac{\tau_5}{2}\left\|\theta-\theta^*\right\|^2\right\},
\end{equation*}
suggesting that the prior places sufficient mass within the set $\Theta_0(r_0)$. For sufficiently large $r_0^2/p$, condition (3.18) in \cite{spokoiny2014bernsteinvonmises} is satisfied. In conjunction with Conditions \ref{condition: local, exponential moment}--\ref{condition: global identification property}, this enables the application of Proposition 3.6 in \cite{spokoiny2014bernsteinvonmises} to conclude that
\begin{equation}
\label{Tail posterior probability}
\frac{\int_{\Theta_0^c(r_0)} \exp \left\{M_n(\theta; x)\right\} \mathrm{d} \theta}{\int_{\Theta_0(r_0)} \exp \left\{M_n(\theta; x)\right\} \mathrm{d} \theta}\leq 2 \exp \left\{2 \Delta(r_0,y)+2e^{-y}-y\right\}
\end{equation}
holds with probability at least $1-4e^{-y}$. For a detailed derivation of the right-hand side of \eqref{Tail posterior probability}, we refer to (3.12) and (3.19) in \cite{spokoiny2014bernsteinvonmises}. We obtain from \eqref{Tail posterior probability} that 
\begin{eqnarray}
\label{eq1_xpr}
& & \int_{\Theta_0^c(r_0)} p(\theta) \exp \left\{M_n(\theta; x)\right\} \mathrm{d} \theta \nonumber\\
& & \times \int_{\Theta_0(r_0)}p\left(\theta^*\right)\exp \left\{M_n(\theta ; x)+\left(\theta-\theta^*\right)^{\top}\nabla\log p\left(\theta^*\right)+\frac{\tau_5}{2}\left\|\theta-\theta^*\right\|^2\right\} \mathrm{d} \theta \nonumber\\
&\lesssim&\frac{\int_{\Theta_0^c(r_0)} \exp \left\{M_n(\theta; x)\right\} \mathrm{d} \theta}{\int_{\Theta_0(r_0)} \exp \left\{M_n(\theta; x)\right\} \mathrm{d} \theta} \nonumber\\
& & \times \int_{\Theta_0(r_0)}p\left(\theta^*\right)\exp \left\{M_n(\theta ; x)+\left(\theta-\theta^*\right)^{\top}\nabla\log p\left(\theta^*\right)+\frac{\tau_5}{2}\left\|\theta-\theta^*\right\|^2\right\} \mathrm{d} \theta\nonumber\\
&\leq&2 \exp \left\{2 \Delta(r_0,y)+2e^{-y}-y\right\} \nonumber\\
& & \times \int_{\Theta_0(r_0)}p\left(\theta^*\right)\exp \left\{M_n(\theta ; x)+\left(\theta-\theta^*\right)^{\top}\nabla\log p\left(\theta^*\right)+\frac{\tau_5}{2}\left\|\theta-\theta^*\right\|^2\right\} \mathrm{d} \theta
\end{eqnarray}
holds with probability at least $1-4e^{-y}$. 

Based on Steps 1 and 2, we can obtain an upper bound for the normalizing constant $\log\int p(\theta) \exp \left\{M_n(\theta ; x)\right\} \mathrm{d} \theta$. Specifically, similar to the proof in Step 1, we know that the second term in \eqref{eq1_xpr} is bounded as follows:
\begin{eqnarray}
\label{Upper bound of the edge part of the denominator of the VB ideal}
   & & \int_{\Theta_0(r_0)} p\left(\theta^*\right)\exp \left\{M_n(\theta ; x)+\left(\theta-\theta^*\right)^{\top}\nabla\log p\left(\theta^*\right)+\frac{\tau_5}{2}\left\|\theta-\theta^*\right\|^2\right\} \mathrm{d} \theta \nonumber\\
   &\leq& \int_{\Theta_0(r_0)} p\left(\theta^*\right)\exp \left\{\left(\theta-\theta^*\right)^{\top}\nabla\left[ M_n\left(\theta^* ; x\right)+\log p\left(\theta^*\right)\right]-\frac{1}{2} \left\| \left(D_0^2-\tau_5I_p\right)^{\frac{1}{2}} \left(\theta-\theta^*\right)\right\|^2\right\} \mathrm{d} \theta \nonumber\\
&\leq& \exp \left\{M_n(\theta^*;x)+\Delta(r_0,y)+\log p\left(\theta^*\right)\right\} \nonumber\\
& & \times \left[\frac{(2\pi)^{\frac{p}{2}}}{\det\left(D_0^2-\tau_5I_p\right)^{\frac{1}{2}}} \exp\left\{\frac{1}{2}\left\| \left(D_0^2-\tau_5I_p\right)^{-\frac{1}{2}} \nabla\left[M_n(\theta^* ; x) + \log p(\theta^*)\right]\right\|^2\right\} \right. \\
& & \left.-\int_{\Theta_0^c(r_0)} \exp \left\{\left(\theta-\theta^*\right)^{\top}\nabla\left[ M_n\left(\theta^* ; x\right)+\log p\left(\theta^*\right)\right] - \frac{1}{2} \left\| \left(D_0^2-\tau_5I_p\right)^{\frac{1}{2}} \left(\theta-\theta^*\right)\right\|^2 \right\} \mathrm{d} \theta\right]  \nonumber 
\end{eqnarray}
holds with probability at least $1-e^{-y}$ by Lemma \ref{LAN}. Therefore, there exists a constant $C_1>0$, such that 
\begin{eqnarray*}
& & \int p(\theta) \exp \left\{M_n(\theta ; x)\right\} \mathrm{d} \theta\\
&\leq&\left[1+C_1\exp\left\{2 \Delta(r_0,y)+2e^{-y}-y\right\}\right]\exp \left\{M_n(\theta^*;x)+\Delta(r_0,y)+\log p\left(\theta^*\right)\right\}\\
& & \times \left[\frac{(2\pi)^{\frac{p}{2}}}{\det\left(D_0^2-\tau_5I_p\right)^{\frac{1}{2}}} \exp\left\{\frac{1}{2}\left\| \left(D_0^2-\tau_5I_p\right)^{-\frac{1}{2}} \nabla\left[M_n(\theta^* ; x) + \log p(\theta^*)\right]\right\|^2\right\} \right. \\
& & \left.-\int_{\Theta_0^c(r_0)} \exp \left\{\left(\theta-\theta^*\right)^{\top}\nabla\left[ M_n\left(\theta^* ; x\right)+\log p\left(\theta^*\right)\right] - \frac{1}{2} \left\| \left(D_0^2-\tau_5I_p\right)^{\frac{1}{2}} \left(\theta-\theta^*\right)\right\|^2 \right\} \mathrm{d} \theta\right]
\end{eqnarray*}
holds with probability at least $1-5e^{-y}$ by combining the results \eqref{The division of the denominator of the VB ideal}, \eqref{Upper bound of the central part of the denominator of the VB ideal}, \eqref{eq1_xpr} and \eqref{Upper bound of the edge part of the denominator of the VB ideal}.
Recall that $V_{\theta_0}=D_0^2/n$, and $\Delta_{n,\theta_0}=V_{\theta_0}^{-1}\nabla M_n(\theta^*;x) / \sqrt{n}$. Take the logarithm on both sides of the above inequality, we can show that there exists a constant $C_2>0$, such that
\begin{eqnarray}
\label{upper bound of log tail integral}
& & \log\int p(\theta) \exp \left\{M_n(\theta ; x)\right\} \mathrm{d} \theta \nonumber\\
&\leq& M_n(\theta^*;x)+C_2\Delta(r_0,y)+\log p\left(\theta^*\right)+\frac{p}{2}\log(2\pi)-\frac{1}{2}\log\det\left(V_{\theta_0}-\frac{\tau_5}{n}I_p\right) \nonumber \\
& & -\frac{p}{2}\log n + \frac{1}{2} \left\| \left(V_{\theta_0}-\frac{\tau_5}{n}I_p\right)^{-\frac{1}{2}} \left(V_{\theta_0}\Delta_{n,\theta_0}+\frac{1}{\sqrt{n}}\nabla\log p\left(\theta^*\right)\right) \right\|^2
\end{eqnarray}
holds with probability at least $1-5e^{-y}$. Now, we are going to further analysis the fifth term and the last quadratic term in the right-hand side of \eqref{upper bound of log tail integral}. For the fifth term, it can be shown that 
\begin{eqnarray}
\label{bound of log det}
\log\det\left(V_{\theta_0}-\frac{\tau_5}{n}I_p\right) &=& \log\det\left(V_{\theta_0}\right)+\log\det\left(I_p-\frac{\tau_5}{n}V_{\theta_0}^{-1}\right)\nonumber \\
&=&\log\det\left(V_{\theta_0}\right)-C_3\frac{\mathrm{tr}\left(V_{\theta_0}^{-1}\right)}{n}
\end{eqnarray}
for some constant $C_3>0$ given $\rho\left(V_{\theta_0}^{-1}\right)\lesssim\sqrt{n}$.
For the quadratic term, it can be reexpressed as
\begin{eqnarray*}
& & \left\| \left(V_{\theta_0}-\frac{\tau_5}{n}I_p\right)^{-\frac{1}{2}} \left(V_{\theta_0}\Delta_{n,\theta_0}+\frac{1}{\sqrt{n}}\nabla\log p\left(\theta^*\right)\right) \right\|^2\\
&=& \Delta_{n,\theta_0}^{\top}V_{\theta_0}\Delta_{n,\theta_0}+\frac{\tau_5}{n}\left\|\Delta_{n,\theta_0}\right\|^2+\frac{2}{\sqrt{n}}\Delta_{n,\theta_0}^{\top}\nabla\log p\left(\theta^*\right)\\
& & +\frac{1}{n} \left\| \left(V_{\theta_0}-\frac{\tau_5}{n}I_p\right)^{-\frac{1}{2}}\left(\nabla\log p\left(\theta^*\right)+\frac{\tau_5}{\sqrt{n}}\Delta_{n,\theta_0}\right)\right\|^2
\end{eqnarray*}
under condition $\rho\left(V_{\theta_0}^{-1}\right)\lesssim\sqrt{n}$.
By Cauchy-Schwarz inequality, we have
\begin{equation*}
\left|\frac{2}{\sqrt{n}}\Delta_{n,\theta_0}^{\top}\nabla\log p\left(\theta^*\right)\right| \leq \frac{1}{\sqrt{n}}\left(\left\|\Delta_{n,\theta_0}\right\|^2+\left\|\nabla\log p\left(\theta^*\right)\right\|^2\right),
\end{equation*}
and
\begin{equation*}
\left\| \left(V_{\theta_0}-\frac{\tau_5}{n}I_p\right)^{-\frac{1}{2}}\left(\nabla\log p\left(\theta^*\right)+\frac{\tau_5}{\sqrt{n}}\Delta_{n,\theta_0}\right)\right\|^2
\lesssim \sqrt{n}\left(\left\|\nabla\log p\left(\theta^*\right)\right\|^2+\left\|\frac{\tau_5}{\sqrt{n}}\Delta_{n,\theta_0}\right\|^2\right).
\end{equation*}
Therefore, there exists a constant $C_4>0$ such that
\begin{eqnarray}
\label{bound of quadratic term}
& & \left\| \left(V_{\theta_0}-\frac{\tau_5}{n}I_p\right)^{-\frac{1}{2}} \left(V_{\theta_0}\Delta_{n,\theta_0}+\frac{1}{\sqrt{n}}\nabla\log p\left(\theta^*\right)\right) \right\|^2\nonumber\\
&\leq&\Delta_{n,\theta_0}^{\top}V_{\theta_0}\Delta_{n,\theta_0}+\frac{C_4}{\sqrt{n}}\left(\left\|\Delta_{n,\theta_0}\right\|^2+\left\|\nabla\log p\left(\theta^*\right)\right\|^2\right).
\end{eqnarray}
Finally, combining \eqref{upper bound of log tail integral} - \eqref{bound of quadratic term}, we obtain an upper bound for the normalizing constant $\log\int p(\theta) \exp \left\{M_n(\theta ; x)\right\} \mathrm{d} \theta$. That is, 
\begin{eqnarray}
\label{Upper bound of the logarithm of the denominator of the VB ideal omit}
& & \log\int p(\theta) \exp \left\{M_n(\theta ; x)\right\} \mathrm{d} \theta+\frac{1}{2}\log\det\left(V_{\theta_0}\right)+\frac{p}{2}\log n \nonumber\\
& & -M_n(\theta^*;x)-\log p\left(\theta^*\right)-\frac{p}{2}\log(2\pi)-\frac{1}{2}\Delta_{n,\theta_0}^{\top}V_{\theta_0}\Delta_{n,\theta_0}\nonumber\\
&\lesssim&\Delta(r_0,y)+\frac{1}{\sqrt{n}} \left(\left\|\Delta_{n,\theta_0}\right\|^2+\left\|\nabla\log p\left(\theta^*\right)\right\|^2\right) + \frac{1}{n}\mathrm{tr}\left(V_{\theta_0}^{-1}\right)
\end{eqnarray}
holds with probability at least $1-5e^{-y}$. In more general cases of the prior, the variational log-likelihood dominates the prior due to the implicit sample size $n$. By employing similar techniques, we can also derive the relaxed results presented in \eqref{Upper bound of the logarithm of the denominator of the VB ideal omit}. For further details, please refer to Propositions 3.6 and 3.7 in \cite{spokoiny2014bernsteinvonmises}.

\textit{Step 3.} We derive a lower bound for $\int_{\Theta_0\left(r_0\right)} p(\theta) \exp \left\{M_n(\theta; x)\right\} \mathrm{d} \theta$. 

Analogous to the derivation of \eqref{Upper bound of the central part of the denominator of the VB ideal}, we can similarly use \eqref{LAN expansion} and \eqref{Second order control for the prior} to establish the lower bound
\begin{eqnarray*}
& & \int_{\Theta_0(r_0)} p(\theta) \exp \left\{M_n(\theta ; x)\right\} \mathrm{d} \theta\\
&\geq& \int_{\Theta_0(r_0)} p(\theta) \exp \left\{M_n(\theta^*;x)+\left(\theta-\theta^*\right)^{\top} \nabla M_n\left(\theta^* ; x\right)-\frac{1}{2} \left\|D_0\left(\theta-\theta^*\right)\right\|^2 -\Delta(r_0,y)\right\} \mathrm{d} \theta\\
&\geq& \exp \left\{M_n(\theta^*;x)+\log p\left(\theta^*\right)-\Delta(r_0,y)\right\}\\
& & \times\int_{\Theta_0(r_0)} \exp \left\{\left(\theta-\theta^*\right)^{\top}\nabla\left[ M_n\left(\theta^* ; x\right)+\log p\left(\theta^*\right)\right] - 
 \frac{1}{2} \left\| \left(D_0^2+\tau_5I_p\right)^{\frac{1}{2}} \left(\theta-\theta^*\right)\right\|^2\right\} \mathrm{d} \theta\\
&=& \exp\left\{M_n(\theta^*;x)+\log p\left(\theta^*\right)-\Delta(r_0,y)\right\}\\
& & \times \left[\frac{(2\pi)^{\frac{p}{2}}}{\det\left(D_0^2+\tau_5I_p\right)^{\frac{1}{2}}} \exp\left\{\frac{1}{2}\left\| \left(D_0^2+\tau_5I_p\right)^{-\frac{1}{2}} \nabla\left[M_n(\theta^* ; x) + \log p(\theta^*)\right]\right\|^2\right\} \right.\\
& & \left.-\int_{\Theta_0^c(r_0)} \exp \left\{\left(\theta-\theta^*\right)^{\top}\nabla\left[ M_n\left(\theta^* ; x\right)+\log p\left(\theta^*\right)\right] - \frac{1}{2} \left\| \left(D_0^2+\tau_5I_p\right)^{\frac{1}{2}} \left(\theta-\theta^*\right)\right\|^2\right\} \mathrm{d} \theta\right]
\end{eqnarray*}
which holds with probability at least $1 - e^{-y}$. To simplify the notation, we denote 
\begin{equation*}
T \triangleq \int_{\Theta_0^c(r_0)} \exp \left\{\left(\theta-\theta^*\right)^{\top}\nabla\left[ M_n\left(\theta^* ; x\right)+\log p\left(\theta^*\right)\right]-\frac{1}{2} \left\| \left(D_0^2+\tau_5I_p\right)^{\frac{1}{2}} \left(\theta-\theta^*\right)\right\|^2\right\} \mathrm{d} \theta.
\end{equation*}
We then further analysis this integral. By the change of variable $s=D_0\left(\theta-\theta^*\right)$, the integral $T$ can be rewritten as:
\begin{eqnarray*}
T&=& \frac{1}{\mathrm{det}\left(D_0\right)}\int_{\left\{\left\|s\right\|\geq r_0\right\}} \exp\left\{s^{\top}V_{\theta_0}^{\frac{1}{2}}\Delta_{n,\theta_0} + \frac{1}{\sqrt{n}}s^{\top}V_{\theta_0}^{-\frac{1}{2}}\nabla\log p\left(\theta^*\right) -\frac{1}{2} \left\|\left(I_p+\frac{\tau_5}{n}V_{\theta_0}^{-1}\right)^{-\frac{1}{2}}s\right\|^2\right\}\mathrm{d}s\\
&=& \frac{(2\pi)^{\frac{p}{2}}}{\det\left(D_0^2+\tau_5I_p\right)^{\frac{1}{2}}}
\exp\left\{ \frac{1}{2} \left\| \left(D_0^2+\tau_5I_p\right)^{-\frac{1}{2}} \nabla\left[M_n(\theta^* ; x) + \log p(\theta^*)\right]\right\|^2 \right\}\\
 & & \times \det\left(I_p+\frac{\tau_5}{n}V_{\theta_0}^{-1}\right)^{\frac{1}{2}} \exp\left\{ -\frac{1}{2} \left\| \left(D_0^2+\tau_5I_p\right)^{-\frac{1}{2}} \nabla\left[M_n(\theta^* ; x) + \log p(\theta^*)\right]\right\|^2 \right\}\\
& & \times \int_{\left\{\left\|s\right\|\geq r_0\right\}} \frac{1}{(2\pi)^{\frac{p}{2}}} \exp\left\{s^{\top}V_{\theta_0}^{\frac{1}{2}}\Delta_{n,\theta_0} + \frac{1}{\sqrt{n}}s^{\top}V_{\theta_0}^{-\frac{1}{2}}\nabla\log p\left(\theta^*\right) -\frac{1}{2} 
\left\|\left(I_p+\frac{\tau_5}{n}V_{\theta_0}^{-1}\right)^{\frac{1}{2}}s\right\|^2
\right\}\mathrm{d}s.
\end{eqnarray*}
Therefore, 
\begin{eqnarray}
\label{eq2-xpr}
& & \int_{\Theta_0(r_0)} p(\theta) \exp \left\{M_n(\theta ; x)\right\} \mathrm{d} \theta\nonumber\\
&\geq& \exp \left\{M_n(\theta^*;x)+\log p\left(\theta^*\right)-\Delta(r_0,y)\right\}\nonumber\\
& &\times  \frac{(2\pi)^{\frac{p}{2}}}{\det\left(D_0^2+\tau_5I_p\right)^{\frac{1}{2}}} \exp\left\{\frac{1}{2}\left\| \left(D_0^2+\tau_5I_p\right)^{-\frac{1}{2}} \nabla\left[M_n(\theta^* ; x) + \log p(\theta^*)\right]\right\|^2\right\}\nonumber\\
& & \times\left[1 -  
\det\left(I_p+\frac{\tau_5}{n}V_{\theta_0}^{-1}\right)^{\frac{1}{2}} \exp\left\{ -\frac{1}{2} \left\| \left(D_0^2+\tau_5I_p\right)^{-\frac{1}{2}} \nabla\left[M_n(\theta^* ; x) + \log p(\theta^*)\right]\right\|^2 \right\}\right. \\
& & \left.\times 
\int_{\left\{\left\|s\right\|\geq r_0\right\}} \frac{1}{(2\pi)^{\frac{p}{2}}} \exp\left\{s^{\top}V_{\theta_0}^{\frac{1}{2}}\Delta_{n,\theta_0} + \frac{1}{\sqrt{n}}s^{\top}V_{\theta_0}^{-\frac{1}{2}}\nabla\log p\left(\theta^*\right) -\frac{1}{2} 
\left\|\left(I_p+\frac{\tau_5}{n}V_{\theta_0}^{-1}\right)^{\frac{1}{2}}s\right\|^2
\right\}\mathrm{d}s \right]\nonumber
\end{eqnarray}
holds with probability at least $1 - e^{-y}$.  

To further refine the analysis of the above lower bound \eqref{eq2-xpr}, we now derive an upper bound for the term 
\begin{equation}
\label{third step quasi-normal tail}
\int_{\left\{\left\|s\right\|\geq r_0\right\}} \frac{1}{(2\pi)^{\frac{p}{2}}} \exp\left\{s^{\top}V_{\theta_0}^{\frac{1}{2}}\Delta_{n,\theta_0} + \frac{1}{\sqrt{n}}s^{\top}V_{\theta_0}^{-\frac{1}{2}}\nabla\log p\left(\theta^*\right) -\frac{1}{2} 
\left\|\left(I_p+\frac{\tau_5}{n}V_{\theta_0}^{-1}\right)^{\frac{1}{2}}s\right\|^2
\right\}\mathrm{d}s,
\end{equation}
and show it is exponentially negligible in the form $\exp\left\{-\xi p\right\}$ for some constant $\xi>0$. Note that 
\begin{equation*}
\mathbb{E}_{\theta_0}\left\|V_{\theta_0}^{\frac{1}{2}}\Delta_{n,\theta_0}\right\|^2=\mathrm{tr}\left(V_{\theta_0}^{-1}\operatorname{Var}_{\theta_0}\left\{\nabla m\left(\theta^*; x_i\right)\right\}\right)\leq\tau_4^2p
\end{equation*}
under Condition \ref{condition: identifiability}. Then, employing the Chernoff bound technique, we have 
\begin{eqnarray*}
P_{\theta_0}\left(\left\|V_{\theta_0}^{\frac{1}{2}}\Delta_{n,\theta_0}\right\|^2>\beta^2 r_0^2\right) &=& P_{\theta_0}\left(\left\|V_{\theta_0}^{\frac{1}{2}}\Delta_{n,\theta_0}\right\|>\beta r_0\right) 
\leq\frac{\mathbb{E}_{\theta_0}\exp\left\{\frac{\lambda}{\tau_4}\left\|V_{\theta_0}^{\frac{1}{2}}\Delta_{n,\theta_0}\right\|\right\}}{\exp\left\{\frac{\lambda}{\tau_4}\beta r_0\right\}}\\
& \leq& \frac{\mathbb{E}_{\theta_0} \exp \left\{\lambda \frac{\Delta_{n,\theta_0}^{\top} \nabla \left[M_n\left(\theta^* ; x\right)-\mathbb{E}_{\theta_0}M_n\left(\theta^* ; x\right)\right]}{\left\|\Sigma_0 \Delta_{n,\theta_0}\right\|}\right\}}{\exp\left\{\frac{\lambda}{\tau_4}\beta r_0\right\}}\leq \exp\left\{\frac{\tau_1^2}{2}\lambda^2-\frac{\beta r_0}{\tau_4}\lambda\right\}
\end{eqnarray*}
for any constant $\beta\in\left(0,1/\sqrt{2}\right)$ under Conditions \ref{condition: local, exponential moment} and \ref{condition: identifiability}. Assume that $\tau_2$ satisfies $\beta r_0/(\tau_1^2\tau_4)\leq\tau_2$, as indicated by Remark 2.2 in \cite{spokoiny2014bernsteinvonmises}. Consequently, 
\begin{equation*}
P_{\theta_0}\left(\left\|V_{\theta_0}^{\frac{1}{2}}\Delta_{n,\theta_0}\right\|^2>\beta^2 r_0^2\right)\leq\exp\left\{-\frac{\beta^2r_0^2}{2\tau_1^2\tau_4^2}\right\}.
\end{equation*}
On the other hand, the triangle inequality states that 
\begin{equation*}
\left\|s\right\|^2\leq2\left(\left\|s-V_{\theta_0}^{\frac{1}{2}}\Delta_{n,\theta_0}\right\|^2+\left\|V_{\theta_0}^{\frac{1}{2}}\Delta_{n,\theta_0}\right\|^2\right).
\end{equation*}
Therefore, in the integration region $\{\|s\| \geq r_0\}$, we have 
\begin{equation*}
P_{\theta_0}\left(\left\|s-V_{\theta_0}^{\frac{1}{2}}\Delta_{n,\theta_0}\right\|^2\geq\left(\frac{1}{2}-\beta^2\right)r_0^2\right)\geq 1-\exp\left\{-\frac{\beta^2r_0^2}{2\tau_1^2\tau_4^2}\right\}.
\end{equation*}
Consequently, for sufficiently large $r_0^2/p$, we can show that there exists a constant $C_5\in\left(0,1/2\right)$ such that
\begin{eqnarray}
& & \int_{\left\{\left\|s\right\|\geq r_0\right\}} \frac{1}{(2\pi)^{\frac{p}{2}}} \exp\left\{-\frac{1}{2} \left\|s-V_{\theta_0}^{\frac{1}{2}}\Delta_{n,\theta_0}\right\|^2 +\frac{1}{\sqrt{n}} s^{\top}V_{\theta_0}^{-\frac{1}{2}}\nabla\log p\left(\theta^*\right) -\frac{\tau_5}{2n}s^{\top}V_{\theta_0}^{-1}s\right\}
\mathrm{d}s \nonumber\\
&\leq& \int_{\left\{\left\|s-V_{\theta_0}^{\frac{1}{2}}\Delta_{n,\theta_0}\right\|\geq \sqrt{\frac{1}{2}-\beta^2}r_0\right\}} \frac{1}{(2\pi)^{\frac{p}{2}}} \exp\left\{-C_5\left\|s-V_{\theta_0}^{\frac{1}{2}}\Delta_{n,\theta_0}\right\|^2\right\}\mathrm{d}\left(s-V_{\theta_0}^{\frac{1}{2}}\Delta_{n,\theta_0}\right) \nonumber\\ 
&=&\int_{\left\{\left\|l\right\|\geq \sqrt{\frac{1}{2}-\beta^2}r_0\right\}} \frac{1}{(2\pi)^{\frac{p}{2}}} \exp\left\{-C_5\left\|l\right\|^2\right\}\mathrm{d}l
\lesssim \frac{\left(\sqrt{\frac{1}{2}-\beta^2}r_0\right)^{p-2}}{(p-2)!!\exp \left\{\left(\frac{1}{2}-\beta^2\right)C_5r_0^2\right\}}
\label{upper bound of the normal tail}
\end{eqnarray}
holds with probability at least $1-\exp\left\{-\beta^2r_0^2 / (2\tau_1^2\tau_4^2) \right\}$. To further show this upper bound is exponentially negligible, we apply the Stirling's formula to this upper bound and obtain
\begin{equation}
\label{eq3-xpr}
\frac{\left(\sqrt{\frac{1}{2}-\beta^2}r_0\right)^{p-2}}{(p-2)!!\exp \left\{\left(\frac{1}{2}-\beta^2\right)C_5r_0^2\right\}}\lesssim\frac{\left[\left(\frac{1}{2}-\beta^2\right)r_0^2\right]^{\frac{p-2}{2}}}{\sqrt{p-2}\left(\frac{p-2}{e}\right)^{\frac{p-2}{2}}\exp \left\{\left(\frac{1}{2}-\beta^2\right)C_5r_0^2\right\}}.
\end{equation}
Consider the function
\begin{equation*}
f(x)=\frac{p-2}{2}\log x-C_5x.
\end{equation*}
For sufficiently large $r_0^2/p$, there exists a constant $\alpha$, independent of other parameters, such that
\begin{equation}
\label{assumption for alpha C3}
    \frac{p-1}{\left(1-2\beta^2\right)C_5r_0^2} \leq\alpha<1 \quad \mbox{and} \quad 1+\log\frac{1}{2C_5\alpha}-\frac{1}{\alpha}<0.
\end{equation}
Hence, it is easy to show that $f^\prime(x) <0$ for $x=\left(\frac{1}{2}-\beta^2\right)r_0^2\geq (p-1)/(2C_5\alpha)$. Exploiting the monotonic decreasing property of $f(x)$, we apply the logarithm to the right-hand side of \eqref{eq3-xpr}, which yields
\begin{eqnarray}
\label{logarithmic analysis 1}
& & \log\left(\frac{\left[\left(\frac{1}{2}-\beta^2\right)r_0^2\right]^{\frac{p-2}{2}}}{\sqrt{p-2}\left(\frac{p-2}{e}\right)^{\frac{p-2}{2}}\exp \left\{\left(\frac{1}{2}-\beta^2\right)C_5r_0^2\right\}}\right)\nonumber \\
&=&\frac{p-2}{2}\log\left(\left(\frac{1}{2}-\beta^2\right)r_0^2\right)-C_5\left(\frac{1}{2}-\beta^2\right)r_0^2-\frac{p-2}{2}\log(p-2)+\frac{p-2}{2}-\frac{1}{2}\log(p-2)\nonumber\\
&\leq& \frac{p-2}{2}\left\{1+\log\left(\frac{1}{2C_5\alpha}\right)-\frac{1}{\alpha}\right\}-\frac{1}{2}\log(p-2)+\frac{p-2}{2}\log\left(1+\frac{1}{p-2}\right)-\frac{1}{2\alpha}. 
\end{eqnarray}
Note that the leading term on the right-hand side of \eqref{logarithmic analysis 1} is negative due to the appropriately chosen $\alpha$ in \eqref{assumption for alpha C3}. Combined with \eqref{upper bound of the normal tail} and \eqref{eq3-xpr}, it implies that the term \eqref{third step quasi-normal tail} is exponentially negligible in the form $\exp\left\{-\xi p\right\}$ for some constant $\xi>0$. Therefore, after some tedious but straightforward computations, we can show that
the following inequality holds based on \eqref{eq2-xpr} with  probability at least $1-\exp\left\{-y\right\}-\exp\left\{-\beta^2r_0^2 / (2\tau_1^2\tau_4^2) \right\}$:
\begin{eqnarray}
\label{Lower bound of the center part of the denominator of the VB ideal}
& &\int_{\Theta_0(r_0)} p(\theta) \exp \left\{M_n(\theta ; x)\right\} \mathrm{d} \theta \nonumber \\
&\geq &\exp \left\{M_n(\theta^*;x)+\log p\left(\theta^*\right)-C_2\left( \Delta(r_0,y)+\frac{1}{\sqrt{n}} \left(\left\|\Delta_{n,\theta_0}\right\|^2+\left\|\nabla\log p\left(\theta^*\right)\right\|^2\right)+\frac{1}{n}\mathrm{tr}\left(V_{\theta_0}^{-1}\right) \right)\right\} \nonumber  \\
& & \times \frac{(2\pi)^{\frac{p}{2}}}{\det\left(D_0^2+\tau_5I_p\right)^{\frac{1}{2}}}\exp\left\{\frac{1}{2} \left\|\left(D_0^2+\tau_5I_p\right)^{-\frac{1}{2}}\nabla\left[ M_n\left(\theta^* ; x\right)+\log p\left(\theta^*\right)\right]\right\|^2\right\}.
\end{eqnarray}

\textit{Step 4.} We derive a lower bound for $\int_{\Theta_0^c\left(r_0\right)} p(\theta) \exp \left\{M_n(\theta; x)\right\} \mathrm{d} \theta$. 

It follows directly that 
\begin{equation}
\label{Lower bound of the edge part of the denominator of the VB ideal}
\int_{\Theta_0^c\left(r_0\right)} p(\theta) \exp \left\{M_n(\theta; x)\right\} \mathrm{d} \theta>0.
\end{equation}

Therefore, based on Steps 3 and 4, we can obtain a lower bound for the normalizing constant $\log\int p(\theta) \exp \left\{M_n(\theta ; x)\right\} \mathrm{d} \theta$. Specifically, by combing \eqref{Lower bound of the center part of the denominator of the VB ideal} and \eqref{Lower bound of the edge part of the denominator of the VB ideal}, we know that
\begin{eqnarray*}
   & & \log\int p(\theta) \exp \left\{M_n(\theta ; x)\right\} \mathrm{d} \theta\\ 
  &>& \log\int_{\Theta_0(r_0)} p(\theta) \exp \left\{M_n(\theta ; x)\right\} \mathrm{d} \theta\\
  &\geq& 
  M_n(\theta^*;x)+\log p\left(\theta^*\right) - C_2\left( \Delta(r_0,y)+\frac{1}{\sqrt{n}} \left(\left\|\Delta_{n,\theta_0}\right\|^2+\left\|\nabla\log p\left(\theta^*\right)\right\|^2\right)+\frac{1}{n}\mathrm{tr}\left(V_{\theta_0}^{-1}\right) \right)\\
   & & + \frac{p}{2}\log(2\pi)-\frac{p}{2}\log n-\frac{1}{2}\log\det\left(V_{\theta_0}\right)
-\frac{1}{2}\log\det\left(I_p+\frac{\tau_5}{n}V_{\theta_0}^{-1}\right)\\
& & +\frac{1}{2}\Delta_{n,\theta_0}^{\top}V_{\theta_0}\Delta_{n,\theta_0} -\frac{\tau_5}{2n}\left\|\Delta_{n,\theta_0}\right\|^2 
+\frac{1}{\sqrt{n}} \Delta_{n,\theta_0}^{\top}\nabla\log p\left(\theta^*\right)\\
& & + \frac{1}{n} \left\| \left(V_{\theta_0}+\frac{\tau_5}{n}I_p\right)^{-\frac{1}{2}}\left(\nabla\log p\left(\theta^*\right)-\frac{\tau_5}{\sqrt{n}}\Delta_{n,\theta_0}\right) \right\|^2
\end{eqnarray*}
holds with probability at least  $1-\exp\left\{-y\right\}-\exp\left\{-\beta^2r_0^2 / (2\tau_1^2\tau_4^2) \right\}$ for  $\beta\in\left(0,1/\sqrt{2}\right)$. For the term $\log\det\left(I_p+\frac{\tau_5}{n}V_{\theta_0}^{-1}\right)$ and the quadratic term
$\left\| \left(V_{\theta_0}+\frac{\tau_5}{n}I_p\right)^{-\frac{1}{2}}\left(\nabla\log p\left(\theta^*\right)-\frac{\tau_5}{\sqrt{n}}\Delta_{n,\theta_0}\right) \right\|^2$,
we can do similar analysis as in \eqref{bound of log det} and \eqref{bound of quadratic term}, respectively. Then, we can show that 
\begin{eqnarray}
\label{Lower bound of the logarithm of the denominator of the VB ideal}
& & M_n(\theta^*;x)+\log p\left(\theta^*\right)+\frac{p}{2}\log(2\pi)+\frac{1}{2}\Delta_{n,\theta_0}^{\top}V_{\theta_0}\Delta_{n,\theta_0} \nonumber \\
& & -\log\int p(\theta) \exp \left\{M_n(\theta ; x)\right\} \mathrm{d} \theta-\frac{1}{2}\log\det\left(V_{\theta_0}\right)-\frac{p}{2}\log n \nonumber\\
&\lesssim& \Delta(r_0,y) +\frac{1}{\sqrt{n}} \left(\left\|\Delta_{n,\theta_0}\right\|^2+\left\|\nabla\log p\left(\theta^*\right)\right\|^2\right) +\frac{1}{n}\mathrm{tr}\left(V_{\theta_0}^{-1}\right)  
\end{eqnarray}
holds with probability at least  $1-\exp\left\{-y\right\}-\exp\left\{-\beta^2r_0^2 / (2\tau_1^2\tau_4^2) \right\}$. 

Finally, by integrating \eqref{Upper bound of the logarithm of the denominator of the VB ideal omit} with \eqref{Lower bound of the logarithm of the denominator of the VB ideal}, a finite-sample bound on the logarithm of the normalizing constant of the VB ideal posterior can be derived directly. That is,
\begin{eqnarray*}
& & \left|\log\int p(\theta) \exp \left\{M_n(\theta ; x)\right\} \mathrm{d} \theta+\frac{1}{2}\log\det\left(V_{\theta_0}\right)+\frac{p}{2}\log n\right.\\
& & \left.-M_n(\theta^*;x)-\log p\left(\theta^*\right)-\frac{p}{2}\log(2\pi)-\frac{1}{2}\Delta_{n,\theta_0}^{\top}V_{\theta_0}\Delta_{n,\theta_0} \right|\\
&\lesssim&\Delta(r_0,y)+\frac{1}{\sqrt{n}}\left(\left\|\Delta_{n,\theta_0}\right\|^2+\left\|\nabla\log p\left(\theta^*\right)\right\|^2\right) +\frac{1}{n}\mathrm{tr}\left(V_{\theta_0}^{-1}\right)   
\end{eqnarray*}
holds with probability at least $1-5e^{-y}-e^{-cr_0^2}$ for some constant $c\in\left(0,1/4\tau_1^2\tau_4^2\right)$.

\section{Proof of Lemma \ref{lemma: The Tail Moment Control for the p-Dimensional Normal Distribution}}
\label{Proof of Lemma 3}

Let $l=(l_1,\ldots,l_p)^{\top}\in\mathbb{R}^p$. Using the multivariate polar coordinate transformation, we obtain
\begin{equation*}
\begin{aligned}
l_1 & =r \cos \left(\theta_1\right), \\
l_2 & =r \sin \left(\theta_1\right) \cos \left(\theta_2\right), \\
l_3 & =r \sin \left(\theta_1\right) \sin \left(\theta_2\right) \cos \left(\theta_3\right), \\
\vdots & \\
l_{p-1} & =r \sin \left(\theta_1\right) \sin \left(\theta_2\right) \cdots \sin \left(\theta_{p-2}\right) \cos \left(\theta_{p-1}\right), \\
l_p & =r \sin \left(\theta_1\right) \sin \left(\theta_2\right) \cdots \sin \left(\theta_{p-2}\right)\sin \left(\theta_{p-1}\right),
\end{aligned}
\end{equation*}
where $r \geq 0, \theta_1, \ldots, \theta_{p-2} \in[0, \pi]$, and $\theta_{p-1} \in[0,2 \pi)$. The integration factor becomes
\begin{equation*}
\left|\mathrm{det}\left(\frac{\partial(l_1, \ldots, l_p)}{\partial(r, \theta_1, \ldots, \theta_{p-1})}\right)\right| = r^{p-1} \prod_{k=1}^{p-2} \sin^{p-k-1}(\theta_k).
\end{equation*}
We then have
\begin{eqnarray*}
& & \int_{\left\{\left\|l\right\|>B\right\}} N\left(l;0,I_p\right)\left\|l\right\|^K\mathrm{d}l\\
&=&\int^{+\infty}_B\int_0^\pi\cdots\int_0^\pi\int_0^{2\pi}\frac{1}{(2\pi)^{\frac{p}{2}}} e^{-\frac{r^2}{2}} r^{p-1+K}\sin ^{p-2}\left(\theta_1\right)\cdots \sin \left(\theta_{p-2}\right)\mathrm{d}r\mathrm{d}\theta_1\cdots\mathrm{d}\theta_{p-2}\mathrm{d}\theta_{p-1}\\
&=&\frac{1}{(2\pi)^{\frac{p}{2}-1}} \int^{+\infty}_B e^{-\frac{r^2}{2}}r^{p-1+K}\mathrm{d}r \prod\limits_{k=1}^{p-2}\left(\int_0^\pi\sin ^{p-1-k}\left(\theta_k\right)\mathrm{d}\theta_k\right).
\end{eqnarray*}
Using the Wallis formula, we derive
\begin{equation*}
\int_0^\pi\sin^k x\mathrm{d}x=\begin{cases}\frac{(k-1)!!}{k!!}\times 2;&2\nmid k,\\ \frac{(k-1)!!}{k!!}\times \pi ;& 2\mid k,\end{cases}
\end{equation*}
which implies that
\begin{equation*}
\prod\limits_{k=1}^{p-2}\left(\int_0^\pi\sin ^{p-1-k}\left(\theta_k\right)\mathrm{d}\theta_k\right)=\frac{2^{\lfloor\frac{p-1}{2}\rfloor}\pi^{\lfloor\frac{p-2}{2}\rfloor}}{(p-2)!!}.
\end{equation*}
Therefore, combining this with
\begin{equation*}
    \frac{2^{\lfloor\frac{p-1}{2}\rfloor}\pi^{\lfloor\frac{p-2}{2}\rfloor}}{(2\pi)^{\frac{p}{2}-1}}=\left(\frac{2}{\pi}\right)^{\frac{p}{2}}\leq 1,
\end{equation*}
we have
\begin{equation*}
\int_{\left\{\left\|l\right\|>B\right\}} N\left(l;0,I_p\right)\left\|l\right\|^K\mathrm{d}l
\leq\frac{1}{(p-2)!!}\int^{+\infty}_B e^{-\frac{r^2}{2}} r^{p-1+K}\mathrm{d}r.
\end{equation*}

Now, we are going to show that
\begin{equation*}
\int_B^{+\infty} e^{-\frac{r^2}{2}}r^{p-1+K}\mathrm{d}r \lesssim B^{p-2+K} e^{-\frac{B^2}{2}}.
\end{equation*}
It can be established via the following two steps.

\textit{Case 1.}  If $2\nmid p-1+K$, then
\begin{eqnarray*}
        & & \int_B^{+\infty} e^{-\frac{r^2}{2}} r^{p-1+K}\mathrm{d}r\\ 
        &=& 2^{\frac{p-2+K}{2}}\int_{\frac{B^2}{2}}^{+\infty}
    x^{\frac{p-2+K}{2}}e^{-x}\mathrm{d}x\\
    &=& 2^{\frac{p-2+K}{2}} e^{-\frac{B^2}{2}} \Gamma\left(\frac{p-2+K}{2} + 1\right) \sum\limits_{i=0}^{\frac{p-2+K}{2}}\frac{1}{i!}\left(\frac{B^2}{2}\right)^i.
 \end{eqnarray*}
Since
\begin{equation*}
\frac{p-2+K}{B^2}<\frac{p-1+K}{B^2}\leq\alpha<1,
\end{equation*}
we can show that
\begin{equation*}
\int_B^{+\infty} e^{-\frac{r^2}{2}}r^{p-1+K}\mathrm{d}r \lesssim B^{p-2+K} e^{-\frac{B^2}{2}}.
\end{equation*}

\textit{Case 2.} If $2\mid p-1+K$, then
 \begin{eqnarray*}
& & \int_B^{+\infty} e^{-\frac{r^2}{2}}  r^{p-1+K}\mathrm{d}r\\
&\leq&\int_B^{+\infty} \frac{r}{B}e^{-\frac{r^2}{2}}  r^{p-1+K}\mathrm{d}r= \frac{1}{B}\int_B^{+\infty} e^{-\frac{r^2}{2}} r^{p+K}\mathrm{d}r\\
&=&\frac{1}{B}2^{\frac{p-1+K}{2}}e^{-\frac{B^2}{2}}  \Gamma\left(\frac{p-1+K}{2}+1\right)  \sum\limits_{i=0}^{\frac{p-1+K}{2}} \frac{1}{i!}\left(\frac{B^2}{2}\right)^i.  
 \end{eqnarray*}
Since
\begin{equation*}
\frac{p-1+K}{B^2}\leq\alpha<1,
\end{equation*}
we can show that
\begin{equation*}
\int_B^{+\infty} e^{-\frac{r^2}{2}} r^{p-1+K}\mathrm{d}r \lesssim B^{p-2+K}e^{-\frac{B^2}{2}}.
\end{equation*}
Therefore, according to cases, we can conclude that
\begin{equation}
\label{Lem3-eq1}
\int_{\left\{\left\|l\right\|>B\right\}} N\left(l;0,I_p\right)\left\|l\right\|^K\mathrm{d}l\lesssim \frac{B^{p-2+K}}{(p-2)!!} e^{-\frac{B^2}{2}}.
\end{equation}
Similarly, we have
\begin{eqnarray*}
& & \int_B^{+\infty} e^{-\frac{r^2}{2}}  r^{p-1+K}\mathrm{d}r\\
&\geq&\int_B^{+\infty} \frac{B}{r}e^{-\frac{r^2}{2}}  r^{p-1+K}\mathrm{d}r\\
&=& 2^{\frac{p-3+K}{2}}B e^{-\frac{B^2}{2}} \Gamma\left(\frac{p-3+K}{2}+1\right) \sum\limits_{i=0}^{\frac{p-3+K}{2}} \frac{1}{i!} \left(\frac{B^2}{2}\right)^i\\
&\gtrsim& B^{p-2+K} e^{-\frac{B^2}{2}},
\end{eqnarray*}
which implies that
\begin{equation}
\label{Lem3-eq2}
\int_{\left\{\left\|l\right\|>B\right\}} N\left(l;0,I_p\right)\left\|l\right\|^K\mathrm{d}l \gtrsim \frac{B^{p-2+K}}{(p-2)!!} e^{-\frac{B^2}{2}}.
\end{equation}
By combining \eqref{Lem3-eq1} and \eqref{Lem3-eq2}, the desired result follows, completing the proof of the lemma.

\section{Proof of Lemma \ref{lemma: the gap between two functionals}}
\label{Proof of Lemma 4}

First, for any distribution $q(z)$, we have
\begin{eqnarray*}
& & \int q(\theta)\left(\log \left[p(\theta) \exp \left\{\int q(z) \log \frac{p(x, z \mid \theta)}{q(z)} \mathrm{d} z\right\}\right]-\log q(\theta)\right) \mathrm{d} \theta\\ 
&\leq&\int q(\theta)\left(\log \left[p(\theta) \exp \left\{\sup _{q(z)\in \mathcal{Q}_z}\int q(z) \log \frac{p(x, z \mid \theta)}{q(z)} \mathrm{d} z\right\}\right]-\log q(\theta)\right) \mathrm{d} \theta.
\end{eqnarray*}
Consequently, 
\begin{eqnarray*}
& & \sup _{q(z)\in \mathcal{Q}_z}\int q(\theta)\left(\log \left[p(\theta) \exp \left\{\int q(z) \log \frac{p(x, z \mid \theta)}{q(z)} \mathrm{d} z\right\}\right]-\log q(\theta)\right) \mathrm{d} \theta\\ 
&\leq& \int q(\theta)\left(\log \left[p(\theta) \exp \left\{\sup _{q(z)\in \mathcal{Q}_z}\int q(z) \log \frac{p(x, z \mid \theta)}{q(z)} \mathrm{d} z\right\}\right]-\log q(\theta)\right) \mathrm{d} \theta.
\end{eqnarray*}
Recall that
\begin{eqnarray}
\label{log-kl}
& & \log \int p(\theta) \exp \left\{M_n(\theta ; x)\right\} \mathrm{d} \theta-\mathrm{KL}\left(q(\theta) \| \pi^*(\theta \mid x)\right) \nonumber\\
&=& \int q(\theta) \left\{ \log \left[p(\theta) \exp\left\{\sup_{q(z) \in \mathcal{Q}_z} \int q(z) \log \frac{p(x,z | \theta)}{q(z)}\mathrm{d}z \right\} \right] - \log q(\theta)\right\} \mathrm{d}\theta,
\end{eqnarray}
and
\begin{equation}
\label{profiled elbo}
 \mathrm{ELBO}_p(q(\theta)) = \sup_{q(z) \in \mathcal{Q}_z} \int q(\theta) \left\{ \log \left[p(\theta) \exp\left\{\int q(z) \log \frac{p(x,z | \theta)}{q(z)}\mathrm{d}z \right\} \right] - \log q(\theta)\right\} \mathrm{d}\theta.
\end{equation}
Therefore, we know that
\begin{equation}
\label{upper bound of profiled elbo}
0\leq\log\int p(\theta)\exp\left\{M_n(\theta ; x)\right\}\mathrm{d}\theta-\mathrm{KL}\left(q(\theta) \| \pi^*(\theta \mid x)\right)-\mathrm{ELBO}_p(q(\theta)) .
\end{equation}

Next, we derive an upper bound for
\begin{equation*}
\log \int p(\theta) \exp \left\{M_n(\theta ; x)\right\} \mathrm{d} \theta-\mathrm{KL}\left(q(\theta) \| \pi^*(\theta \mid x)\right)-\mathrm{ELBO}_p(q(\theta)) .
\end{equation*}
A direct expansion of \eqref{log-kl} and \eqref{profiled elbo} yields
\begin{eqnarray*}
& & \log\int p(\theta)\exp\left\{M_n(\theta ; x)\right\}\mathrm{d}\theta-\mathrm{KL}\left(q(\theta) \| \pi^*(\theta \mid x)\right)\\ 
&=&\int q(\theta) \log p(\theta) \mathrm{d} \theta-\int q(\theta) \log q(\theta) \mathrm{d} \theta+ \int q(\theta) \sup_{q(z) \in \mathcal{Q}_z}\int q(z) \log \frac{p(x, z \mid \theta)}{q(z)} \mathrm{d} z \mathrm{d} \theta,
\end{eqnarray*}
and 
\begin{eqnarray*}
\label{details of profiled elbo}
& & \mathrm{ELBO}_p(q(\theta)) \\
&=& \int q(\theta) \log p(\theta) \mathrm{d} \theta-\int q(\theta) \log q(\theta) \mathrm{d} \theta+\sup_{q(z) \in \mathcal{Q}_z} \int q(\theta) \int q(z) \log \frac{p(x, z \mid \theta)}{q(z)} \mathrm{d} z \mathrm{d} \theta,
\end{eqnarray*}
respectively. Therefore, we only need to focus on bounding
\begin{equation}
\label{Deq1-xpr}
\int q(\theta) \sup_{q(z) \in \mathcal{Q}_z}\int q(z) \log \frac{p(x, z \mid \theta)}{q(z)} \mathrm{d} z \mathrm{d} \theta-\sup_{q(z) \in \mathcal{Q}_z} \int q(\theta) \int q(z) \log \frac{p(x, z \mid \theta)}{q(z)} \mathrm{d} z \mathrm{d} \theta.
\end{equation}
Recall the variational log-likelihood
\begin{equation*}
    M_n(\theta ; x) = \sup_{q(z) \in \mathcal{Q}_z} \int q(z) \log \frac{p(x, z \mid \theta)}{q(z)} \mathrm{d} z.
\end{equation*}
By applying a first-order Taylor expansion, we obtain
\begin{eqnarray*}
M_n(\theta ; x) &= & M_n(\mu ; x)+(\theta-\mu)^{\top} \nabla M_n(\mu ; x) \\
& & +(\theta-\mu)^{\top} \int_0^1(1-t) \nabla^2 M_n(\mu+t(\theta-\mu) ; x) \mathrm{d} t(\theta-\mu).
\end{eqnarray*}
Let $h=\sqrt{n}(\theta-\mu)$ and $q_h(h)=N(h; 0, V)$. Then, for the first term in \eqref{Deq1-xpr}, we have
\begin{eqnarray*}
& & \int q(\theta) \sup_{q(z) \in \mathcal{Q}_z}\int q(z) \log \frac{p(x, z \mid \theta)}{q(z)} \mathrm{d} z \mathrm{d} \theta\\
&=& \int q(\theta)M_n(\theta ; x)\mathrm{d} \theta\\
&=& M_n(\mu ; x)+\int q(\theta)\left(\theta-\mu\right)^{\top}\int_0^1 (1 - t) \nabla^2 M_n(\mu + t(\theta - \mu); x) \mathrm{d}t\left(\theta-\mu\right)\mathrm{d} \theta\\
&=& M_n(\mu ; x)+\frac{1}{n}\int q_h(h)h^\top\left[\int_0^1 (1 - t) \nabla^2 M_n\left(\mu+\frac{th}{\sqrt{n}};x\right) \mathrm{d}t\right] h\mathrm{d}h,
\end{eqnarray*}
by noting that $q(\theta)=N(\theta; \mu, V / n) \in \mathcal{Q}_\theta$. Similarly, a first-order Taylor expansion yields that
\begin{eqnarray*}
& & \int q(z) \log \frac{p(x, z \mid \theta)}{q(z)} \mathrm{d} z\\
&=&\int q(z) \log \frac{p(x, z \mid \mu)}{q(z)} \mathrm{d} z+\left(\theta-\mu\right)^{\top}\nabla\left(\int q(z) \log \frac{p(x, z \mid \mu)}{q(z)} \mathrm{d} z\right)\\
& & +\left(\theta-\mu\right)^{\top}\int_0^1 (1 - t) \nabla^2\left(\int q(z) \log \frac{p(x, z \mid \mu + t(\theta - \mu))}{q(z)} \mathrm{d} z\right) \mathrm{d}t\left(\theta-\mu\right).
\end{eqnarray*}
Therefore, for the second term in \eqref{Deq1-xpr}, we have
\begin{eqnarray*}
& & \int q(\theta)\int q(z) \log \frac{p(x, z \mid \theta)}{q(z)} \mathrm{d} z\mathrm{d} \theta\\
&=& \int q(z) \log \frac{p(x, z \mid \mu)}{q(z)} \mathrm{d} z\\
& & +\int q(\theta)\left(\theta-\mu\right)^{\top}\int_0^1 (1 - t) \nabla^2\left(\int q(z) \log \frac{p(x, z \mid \mu + t(\theta - \mu))}{q(z)} \mathrm{d} z\right) \mathrm{d}t\left(\theta-\mu\right)\mathrm{d} \theta\\
&=& \int q(z) \log \frac{p(x, z \mid \mu)}{q(z)} \mathrm{d} z\\
& & +\frac{1}{n}\int q_h(h)h^\top\left[\int_0^1 (1 - t) \nabla^2 \left(\int q(z) \log \frac{p\left(x, z \mid \mu + \frac{th}{\sqrt{n}}\right)}{q(z)} \mathrm{d} z\right) \mathrm{d}t\right] h\mathrm{d}h.
\end{eqnarray*}
Let $q^{\mu}(z)=\underset{q(z)\in\mathcal{Q}_z}{\arg\max}\int q(z) \log \frac{p(x, z \mid \mu)}{q(z)} d z$. Then,
\begin{equation*}
M_n(\mu ; x)=\int q^{\mu}(z) \log \frac{p(x, z \mid \mu)}{q^{\mu}(z)} d z.
\end{equation*}
Therefore, the term \eqref{Deq1-xpr} can be bounded by
\begin{equation*}
 \frac{1}{n}\int q_h(h)h^\top\left\{\int_0^1 (1 - t)\left[\nabla^2 M_n\left(\mu+\frac{th}{\sqrt{n}};x\right)-\nabla^2 \left(\int q^\mu(z) \log \frac{p\left(x, z \mid \mu + \frac{th}{\sqrt{n}}\right)}{q^\mu(z)} \mathrm{d} z\right)\right] \mathrm{d}t\right\}h\mathrm{d}h.   
\end{equation*}
Decompose the term $\nabla^2 M_n\left(\mu+\frac{th}{\sqrt{n}};x\right)-\nabla^2 \left(\int q^\mu(z) \log \frac{p\left(x, z \mid \mu + \frac{th}{\sqrt{n}}\right)}{q^\mu(z)} \mathrm{d} z\right)$ as follows:
\begin{eqnarray*}
& & \nabla^2 M_n\left(\mu+\frac{th}{\sqrt{n}};x\right)-\nabla^2 \left(\int q^\mu(z) \log \frac{p\left(x, z \mid \mu + \frac{th}{\sqrt{n}}\right)}{q^\mu(z)} \mathrm{d} z\right)\\
&=&\nabla^2 M_n\left(\mu+\frac{th}{\sqrt{n}};x\right)-\nabla^2 \left(\int q^\mu(z) \log \frac{p\left(x, z \mid \mu\right)}{q^\mu(z)} \mathrm{d} z\right)\\
& & +\nabla^2 \left(\int q^\mu(z) \log \frac{p\left(x, z \mid \mu\right)}{p\left(x, z \mid \mu + \frac{th}{\sqrt{n}}\right)} \mathrm{d} z\right)\\
&=& \nabla^2 M_n\left(\mu+\frac{th}{\sqrt{n}};x\right)-\nabla^2 M_n\left(\mu;x\right) \\
& & +\int q^\mu(z)\left[\nabla^2\log p\left(x, z \mid \mu\right)-\nabla^2\log p\left(x, z \mid \mu + \frac{th}{\sqrt{n}}\right)\right]\mathrm{d} z.
\end{eqnarray*}
The upper bound of \eqref{Deq1-xpr} can be rewritten as
\begin{eqnarray}
\label{decomposition of VL and CL}
& & \frac{1}{n}\int q_h(h)h^\top\left\{\int_0^1 (1 - t)\left[\nabla^2 M_n\left(\mu+\frac{th}{\sqrt{n}};x\right)-\nabla^2 M_n\left(\mu;x\right)\right] \mathrm{d}t\right\}h\mathrm{d}h \nonumber\\
& & +\frac{1}{n}\int q_h(h)h^\top\left\{\int_0^1 (1 - t)\int q^\mu(z)\left[\nabla^2\log p\left(x, z \mid \mu\right)-\nabla^2\log p\left(x, z \mid \mu + \frac{th}{\sqrt{n}}\right)\right]\mathrm{d} z \mathrm{d}t\right\}h\mathrm{d}h \nonumber\\
&\triangleq& T_1 + T_2.
\end{eqnarray}

Now, we derive the upper bounds for terms $T_1$ and $T_2$, respectively. For the term $T_1$, we decompose it as follows:
\begin{eqnarray*}
 T_1 &=&  \frac{1}{n}\int q_h(h)h^\top\left\{\int_0^1 (1 - t)\left[\nabla^2\mathbb{E}_{\theta_0} M_n\left(\mu+\frac{th}{\sqrt{n}};x\right)-\nabla^2\mathbb{E}_{\theta_0} M_n\left(\mu;x\right)\right] \mathrm{d}t\right\}h\mathrm{d}h\\
 & & + \frac{1}{n} \int q_h(h)h^\top\left\{\int_0^1 (1 - t)\left[\nabla^2\mathbb{E}_{\theta_0} M_n\left(\mu;x\right)-\nabla^2M_n\left(\mu;x\right)\right] \mathrm{d}t\right\}h\mathrm{d}h\\
 & & + \frac{1}{n} \int q_h(h)h^\top\left\{\int_0^1 (1 - t)\left[\nabla^2 M_n\left(\mu+\frac{th}{\sqrt{n}};x\right)-\nabla^2\mathbb{E}_{\theta_0} M_n\left(\mu+\frac{th}{\sqrt{n}};x\right)\right] \mathrm{d}t\right\}h\mathrm{d}h\\
 &\triangleq& T_{11} + T_{12} + T_{13}.
\end{eqnarray*}
We proceed to analyze the terms $T_{11}, T_{12}$ and $T_{13}$ one by one. By Condition \ref{condition: lipschitz}, we have
\begin{equation*}
    T_{11} \lesssim \frac{1}{n} \int q_h(h)h^\top\left\{\int_0^1 (1 - t)n\left\|\frac{th}{\sqrt{n}}\right\|^{s_1} \mathrm{d}t\right\}h\mathrm{d}h \lesssim \frac{1}{n^{\frac{s_1}{2}}} \int q_h(h)\left\|h\right\|^{2+s_1}\mathrm{d}h.
\end{equation*}
Let $\tilde{h}=V^{-\frac{1}{2}}h$, then $q_{\tilde{h}}(\tilde{h}) = N(\tilde{h}; 0, I_p)$ with $\tilde{h}=(\tilde{h}_1,\cdots, \tilde{h}_p)^\top$. Note that
\begin{equation*}
\left(\frac{1}{p} \sum_{i=1}^p \tilde{h}_i^2\right)^{\frac{1}{2}} \leq\left(\frac{1}{p} \sum_{i=1}^p\left|\tilde{h}_i\right|^{2+s_1}\right)^{\frac{1}{2+s_1}}
\end{equation*}
by the power mean inequality. We then have
\begin{equation*}
\left\|\tilde{h}\right\|^{2+s_1}=\left(\tilde{h}_1^2+\cdots+\tilde{h}_p^2\right)^{\frac{2+s_1}{2}}\leq p^{\frac{s_1}{2}} \sum_{i=1}^p\left|\tilde{h}_i\right|^{2+s_1}.
\end{equation*}
Note that
\begin{equation*}
\int q_{\tilde{h}}(\tilde{h})\left(\sum_{i=1}^p\left|\tilde{h}_i\right|^{2+s_1}\right)\mathrm{d}\tilde{h}=p\mathbb{E}_{N(0,1)}\left|X\right|^{2+s_1}\lesssim p.
\end{equation*}
Therefore,
\begin{eqnarray*}
    T_{11} & \lesssim& \frac{1}{n^{\frac{s_1}{2}}} \int q_h(h)\left\|h\right\|^{2+s_1}\mathrm{d}h = \frac{1}{n^{\frac{s_1}{2}}} \int q_{\tilde{h}}(\tilde{h})\left\|V^{\frac{1}{2}}\tilde{h}\right\|^{2+s_1}\mathrm{d}\tilde{h}\\
    &\leq& \frac{1}{n^{\frac{s_1}{2}}} \rho\left(V\right)^{\frac{2+s_1}{2}}\int q_{\tilde{h}}(\tilde{h})\left\|\tilde{h}\right\|^{2+s_1}\mathrm{d}\tilde{h} \lesssim p\rho\left(V\right) \left(\frac{p\rho\left(V\right)}{n}\right)^{\frac{s_1}{2}}. 
\end{eqnarray*}
For the term $T_{12}$, we use the exponential moment condition \ref{condition: global, exponential moment, variational} to construct a Chernoff bound, effectively controlling the deviation probability. Specifically, for any arbitrarily large constant $M$ and a given $h\in\mathbb{R}^p$, we have 
\begin{eqnarray*}
& & P_{\theta_0}\left(\lambda_{\max}\left(\frac{1}{n} \nabla^2\left[\mathbb{E}_{\theta_0}M_n\left(\mu;x\right)- M_n\left(\mu; x\right)\right]\right)>M\left(\frac{\rho(V)p}{n}\right)^{\frac{s_1}{2}}\right)\\
&=& P_{\theta_0}\left(\frac{1}{n\left\|h\right\|^2} h^\top\nabla^2\left[\mathbb{E}_{\theta_0}M_n\left(\mu;x\right)- M_n\left(\mu; x\right)\right]h > M\left(\frac{\rho(V)p}{n}\right)^{\frac{s_1}{2}}\right)\\
&=& P_{\theta_0}\left(\frac{\lambda}{\tau_3}\frac{h^\top\nabla^2\left[M_n\left(\mu;x\right)-\mathbb{E}_{\theta_0} M_n\left(\mu; x\right)\right](-h)}{nh^\top V_{\theta_0}h}>\frac{\lambda M\left\|h\right\|^2}{\tau_3 h^\top V_{\theta_0}h}\left(\frac{\rho(V)p}{n}\right)^{\frac{s_1}{2}}\right)\\
&\leq& \exp\left(\frac{\tau_1^2}{2}\lambda^2-\frac{M\left\|h\right\|^2}{\tau_3 h^\top V_{\theta_0}h}\left(\frac{\rho(V)p}{n}\right)^{\frac{s_1}{2}}\lambda\right),
\end{eqnarray*}
for $\lambda\in\left[0,\delta_2(r)\right]$. Note that
\begin{equation*}
\frac{M\left\|h\right\|^2}{\tau_1^2\tau_3 h^\top V_{\theta_0}h}\left(\frac{\rho(V)p}{n}\right)^{\frac{s_1}{2}}\asymp \frac{\sqrt{n}M\left\|h\right\|^2}{h^\top V_{\theta_0}h} \left(\frac{\rho(V)p}{n}\right)^{\frac{s_1}{2}},
\end{equation*}
where $\tau_3 \asymp 1 / \sqrt{n}$ in typical cases, as indicated by \cite{spokoiny2014bernsteinvonmises}. Therefore, under Condition \ref{condition: global, exponential moment, variational} that ensures $\delta_2(r)\gtrsim p$, we can choose an appropriate $\lambda$ such that 
\begin{equation*}
P_{\theta_0}\left(\lambda_{\max}\left(\frac{1}{n}\nabla^2\left[\mathbb{E}_{\theta_0}M_n\left(\mu;x\right)- M_n\left(\mu; x\right)\right]\right)>M\left(\frac{\rho(V)p}{n}\right)^{\frac{s_1}{2}}\right)\leq e^{-cp}
\end{equation*}
for some $c>0$. Consequently, 
\begin{equation*}
    T_{12} \lesssim \left(\frac{\rho(V)p}{n}\right)^{\frac{s_1}{2}}\int q_h(h)\left\|h\right\|^2\mathrm{d}h \lesssim p\rho\left(V\right) \left(\frac{p\rho\left(V\right)}{n}\right)^{\frac{s_1}{2}}
\end{equation*}
with probability at least $1-e^{-cp}$. For the term $T_{13}$, there exists a fixed vector $\tilde{\mu}$ such that
\begin{equation*}
T_{13} \lesssim \frac{1}{n}\int q_h(h)h^\top\left\{\int_0^1 (1 - t)\left[\nabla^2 M_n\left(\tilde{\mu};x\right)-\nabla^2\mathbb{E}_{\theta_0} M_n\left(\tilde{\mu};x\right)\right] \mathrm{d}t\right\}h\mathrm{d}h
\end{equation*}
by the integral mean value theorem. Accordingly, by an argument similar to that used for term $T_{12}$, we can establish that
\begin{equation*}
T_{13} \lesssim p\rho\left(V\right) \left(\frac{p\rho\left(V\right)}{n}\right)^{\frac{s_1}{2}}
\end{equation*}
with probability at least $1-e^{-cp}$ by Condition \ref{condition: global, exponential moment, variational}.
In summary, we have the following result for term $T_1$:
\begin{equation}
\label{Deq2-xpr}
    T_1  \lesssim p\rho\left(V\right) \left(\frac{p\rho\left(V\right)}{n}\right)^{\frac{s_1}{2}}
\end{equation}
holds with probability at least $1-e^{-cp}$.

For the term $T_2$, we decompose it as follows:
\begin{eqnarray*}
   T_2  
   &=& \frac{1}{n} \int q_h(h)h^\top\left\{\int_0^1 (1 - t)\int q^\mu(z) \right.\\
   & & \times \left.\left[\nabla^2\mathbb{E}_{\theta_0}\log p\left(x, z \mid \mu\right)-\nabla^2\mathbb{E}_{\theta_0}\log p\left(x, z \mid \mu + \frac{th}{\sqrt{n}}\right)\right]\mathrm{d} z \mathrm{d}t\right\}h\mathrm{d}h\\
    & +& \frac{1}{n} \int q_h(h)h^\top\left\{\int_0^1 (1 - t)\int q^\mu(z)\left[\nabla^2\log p\left(x, z \mid \mu\right)-\nabla^2\mathbb{E}_{\theta_0}\log p\left(x, z \mid \mu\right)\right]\mathrm{d} z \mathrm{d}t\right\}h\mathrm{d}h\\
    & +& \frac{1}{n} \int q_h(h)h^\top\left\{\int_0^1 (1 - t)\int q^\mu(z) \right.\\
    & & \times \left.\left[\nabla^2\mathbb{E}_{\theta_0}\log p\left(x, z \mid \mu + \frac{th}{\sqrt{n}}\right)-\nabla^2\log p\left(x, z \mid \mu + \frac{th}{\sqrt{n}}\right)\right]\mathrm{d} z \mathrm{d}t\right\}h\mathrm{d}h\\
    &\triangleq& T_{21} + T_{22} + T_{23}.
\end{eqnarray*}
By an argument similar to that used for term $T_{11}$, we can establish that
\begin{eqnarray*}
  T_{21}  &\lesssim&  \frac{1}{n} \int q_h(h)h^\top\left\{\int_0^1 (1 - t)\int q^\mu(z)n\left\|\frac{th}{\sqrt{n}}\right\|^{s_2}\mathrm{d} z \mathrm{d}t\right\}h\mathrm{d}h\\
  &\lesssim& \frac{1}{n^{\frac{s_2}{2}}} \int q_h(h)\left\|h\right\|^{2+s_2}\mathrm{d}h \lesssim p\rho\left(V\right) \left(\frac{p\rho\left(V\right)}{n}\right)^{\frac{s_2}{2}}
\end{eqnarray*}
by Condition \ref{condition: lipschitz}. For the terms $T_{22}$ and $T_{23}$, we apply the exponential moment condition \ref{condition: global, exponential moment, complete} to construct a Chernoff bound, thereby effectively controlling their deviation probabilities, as in the cases of $T_{12}$ and $T_{13}$. We can show that
\begin{equation*}
T_{22} \lesssim p\rho\left(V\right) \left(\frac{p\rho\left(V\right)}{n}\right)^{\frac{s_2}{2}}
\end{equation*}
and
\begin{equation*}
T_{23} \lesssim p\rho\left(V\right) \left(\frac{p\rho\left(V\right)}{n}\right)^{\frac{s_2}{2}}
\end{equation*}
hold with probability at least $1-e^{-cp}$. In summary, we have the following result for term $T_2$:
\begin{equation}
\label{Deq3-xpr}
    T_2  \lesssim p\rho\left(V\right) \left(\frac{p\rho\left(V\right)}{n}\right)^{\frac{s_2}{2}}
\end{equation}
holds with probability at least $1-e^{-cp}$.

Therefore, combing the results \eqref{decomposition of VL and CL}, \eqref{Deq2-xpr}, and \eqref{Deq3-xpr}, we have
\begin{eqnarray*}
    & & \int q(\theta) \sup_{q(z) \in \mathcal{Q}_z}\int q(z) \log \frac{p(x, z \mid \theta)}{q(z)} \mathrm{d} z \mathrm{d} \theta-\sup_{q(z) \in \mathcal{Q}_z} \int q(\theta) \int q(z) \log \frac{p(x, z \mid \theta)}{q(z)} \mathrm{d} z \mathrm{d} \theta\\
    & \lesssim & p\rho\left(V\right) \left(\frac{p\rho\left(V\right)}{n}\right)^{\frac{s_1}{2}} + p\rho\left(V\right) \left(\frac{p\rho\left(V\right)}{n}\right)^{\frac{s_2}{2}},
\end{eqnarray*}
with probability at least $1-e^{-cp}$. Further, together with \eqref{upper bound of profiled elbo}, we finish the proof.

\section{Proof of Theorem \ref{theorem: the convergence rates of the VBBE and the covariance matrix of the VB bridge}}
\label{Proof of Theorem 1}

We begin by introducing two $p$-dimensional variational parameters:
\begin{equation*}
(\mathbf{m},\mathbf{v})=\left(\left(m_1,\cdots,m_p\right)^{\top},\left(v_1,\cdots,v_p\right)^{\top}\right),
\end{equation*}
which are used to parameterize the variational distribution $q(\theta)\in\mathcal{Q}_\theta$ as follows:
\begin{equation*}
q(\theta)=q(\theta ;\mathbf{m},\mathbf{v})=N\left(\theta;\mathbf{m},\frac{1}{n}\mathrm{diag}\left(\mathbf{v}\right)\right).
\end{equation*}
Next, we define two functionals that are central to our analysis. The first functional, $F_n(\mathbf{m},\mathbf{v})$, quantifies the Kullback–Leibler (KL) divergence between the variational distribution $q(\theta;\mathbf{m},\mathbf{v})$ and the VB ideal posterior $\pi^*(\theta \mid x)$:
\begin{eqnarray}
\label{def of Fn}
 F_n(\mathbf{m},\mathbf{v}) &\triangleq &\operatorname{KL}\left(q(\theta ;\mathbf{m},\mathbf{v}) \| \pi^*(\theta \mid x)\right) \nonumber\\
&=& \int q(\theta ;\mathbf{m},\mathbf{v}) \log q(\theta ;\mathbf{m},\mathbf{v}) \mathrm{d} \theta-\int q(\theta ;\mathbf{m},\mathbf{v}) \log \pi^*(\theta \mid x) \mathrm{d} \theta \nonumber\\
&=& -\frac{1}{2}p(1+\log(2\pi))-\frac{1}{2}\sum\limits_{i=1}^p\log \frac{v_i}{n}-\int q(\theta ;\mathbf{m},\mathbf{v}) \log p(\theta) \mathrm{d} \theta\nonumber\\
& & +\log\int p(\theta)\exp\left\{M_n(\theta ; x)\right\}\mathrm{d}\theta -\int q(\theta ; \mathbf{m},\mathbf{v}) M_n(\theta ; x) \mathrm{d} \theta.
\end{eqnarray}
The second functional, $F_0(\mathbf{m},\mathbf{v})$, measures the KL divergence between $q(\theta;\mathbf{m},\mathbf{v})$ and the reference Gaussian distribution $q_0(\theta)$:
\begin{eqnarray}
\label{def of F0}
F_0(\mathbf{m},\mathbf{v}) & \triangleq & \mathrm{KL}\left(q(\theta;\mathbf{m},\mathbf{v}) \|q_0(\theta)\right)=\mathrm{KL}\left(q(\theta;\mathbf{m},\mathbf{v}) \|N\left(\theta ; \theta^*+\frac{\Delta_{n, \theta_0}}{\sqrt{n}},\frac{V_{\theta_0}^{-1}}{n}\right)\right) \nonumber\\
& = & \int q(\theta ;\mathbf{m},\mathbf{v}) \log q(\theta ;\mathbf{m},\mathbf{v}) \mathrm{d} \theta-\int q(\theta ;\mathbf{m},\mathbf{v}) \log N\left(\theta ; \theta^*+\frac{\Delta_{n, \theta_0}}{\sqrt{n}},\frac{V_{\theta_0}^{-1}}{n}\right) \mathrm{d} \theta \nonumber\\
&=& -\frac{1}{2}p(1+\log(2\pi)) - \frac{1}{2}\sum\limits_{i=1}^p\log \frac{v_i}{n}+\frac{p}{2}\log(2\pi)+\frac{1}{2}\log\det\left(\frac{V_{\theta_0}^{-1}}{n}\right) \nonumber\\
& & +\frac{n}{2}\int q(\theta ; \mathbf{m},\mathbf{v})  \left(\theta-\theta^*-\frac{\Delta_{n, \theta_0}}{\sqrt{n}}\right)^{\top}V_{\theta_0}\left(\theta-\theta^*-\frac{\Delta_{n, \theta_0}}{\sqrt{n}}\right)\mathrm{d} \theta \nonumber\\
&=& -\frac{p}{2}-\frac{1}{2}\sum\limits_{i=1}^p\log v_i-\frac{1}{2}\log\det\left(V_{\theta_0}\right)+ \frac{1}{2} \mathrm{tr}\left(V_{\theta_0}\mathrm{diag}\left(\mathbf{v}\right)\right) + \frac{1}{2}\Delta_{n, \theta_0}^{\top} V_{\theta_0} \Delta_{n, \theta_0} \nonumber\\
& & +\frac{n}{2}\left(\mathbf{m}-\theta^*\right)^{\top}V_{\theta_0}\left(\mathbf{m}-\theta^*\right)-\sqrt{n}\left(\mathbf{m}-\theta^*\right)^{\top}V_{\theta_0}\Delta_{n, \theta_0}.
\end{eqnarray}
The optimal variational parameters are then defined as the minimizers of these functionals:
\begin{equation*}
(\mathbf{m}_n^*,\mathbf{v}_n^*)=\underset{\mathbf{m},\mathbf{v}}{\arg\min}F_n(\mathbf{m},\mathbf{v}),\quad (\mathbf{m}_0^*,\mathbf{v}_0^*)=\underset{\mathbf{m},\mathbf{v}}{\arg\min}F_0(\mathbf{m},\mathbf{v}).
\end{equation*}
By definition, it follows that
\begin{equation}
\label{def mvn}
\left(\mathbf{m}_n^*, \mathbf{v}_n^*\right) = (\hat{\theta}_n^\ddagger, \hat{v}_n^\ddagger ).
\end{equation}
And a straightforward calculation yields
\begin{equation}
\label{def mv0}
(\mathbf{m}_0^*,\mathbf{v}_0^*)=\left(\theta^*+\frac{\Delta_{n, \theta_0}}{\sqrt{n}},\left(\frac{1}{v_{\theta_0,11}},\ldots,\frac{1}{v_{\theta_0,pp}}\right)^{\top}\right).
\end{equation}

We proceed study the relationship between $F_n(\mathbf{m},\mathbf{v})$ and $F_0(\mathbf{m},\mathbf{v})$. To this end, we perform a detailed analysis of the three integral terms on the right-hand side of \eqref{def of Fn}. By Condition \ref{condition: prior}, a second-order Taylor expansion of $\log p(\theta)$ around $\mathbf{m}$ yields
\begin{equation*}
\left|\log p(\theta)-\log p\left(\mathbf{m}\right)-\left(\theta-\mathbf{m}\right)^{\top}\nabla\log p\left(\mathbf{m}\right)\right|\leq\frac{\tau_5}{2} \left\|\theta-\mathbf{m}\right\|^2,
\end{equation*}
which is analogous to the bound derived in \eqref{Second order control for the prior}. Similarly, we can show that 
\begin{equation}
\label{q m v logp}
\int q(\theta ;\mathbf{m},\mathbf{v}) \log p(\theta) \mathrm{d} \theta=\log p(\mathbf{m})+O\left(\frac{1}{n}\mathrm{tr}\left(\mathrm{diag}\left(\mathbf{v}\right)\right)\right).
\end{equation}

For the term $\log\int p(\theta)\exp\left\{M_n(\theta; x)\right\}\mathrm{d}\theta$, Lemma \ref{Lemma: the bracketing result of the logarithm of the denominator of the VB Ideal} ensures that
\begin{eqnarray}
\label{The Bracket Results of log}
& & \log\int p(\theta) \exp \left\{M_n(\theta ; x)\right\} \mathrm{d} \theta \nonumber\\
&=& M_n(\theta^*;x)+\log p\left(\theta^*\right)+\frac{p}{2}\log(2\pi)-\frac{1}{2}\log\det\left(V_{\theta_0}\right)-\frac{p}{2}\log n+\frac{1}{2}\Delta_{n,\theta_0}^{\top}V_{\theta_0}\Delta_{n,\theta_0} \nonumber\\
& & +O\left(\Delta(r_0,y)+\frac{1}{\sqrt{n}} \left\{\left\|\Delta_{n,\theta_0}\right\|^2+\left\|\nabla\log p\left(\theta^*\right)\right\|^2\right\} + \frac{1}{n}\mathrm{tr}\left(V_{\theta_0}^{-1}\right)\right)
\end{eqnarray}
holds with probability at least $1-5e^{-y}-e^{-cr_0^2}$ for some constant $c>0$.

To evaluate $\int q(\theta;\mathbf{m},\mathbf{v}) M_n(\theta; x) \mathrm{d} \theta$, we decompose it as follows:
\begin{equation*}
\int q(\theta ;\mathbf{m},\mathbf{v}) M_n(\theta ; x) \mathrm{d} \theta=\int_{\Theta_0\left(r_0\right)} q(\theta ;\mathbf{m},\mathbf{v}) M_n(\theta ; x) \mathrm{d} \theta+\int_{\Theta_0^c\left(r_0\right)} q(\theta ;\mathbf{m},\mathbf{v}) M_n(\theta ; x) \mathrm{d} \theta.
\end{equation*}
Note that Lemma \ref{LAN} guarantees the existence of some  $R\in\left[-\Delta(r_0,y),\Delta(r_0,y)\right]$ such that
\begin{eqnarray*}
  & & \int_{\Theta_0\left(r_0\right)} q(\theta ;\mathbf{m},\mathbf{v}) M_n(\theta ; x) \mathrm{d} \theta\\
  &=& \int_{\Theta_0\left(r_0\right)} q(\theta ;\mathbf{m},\mathbf{v})\left[M_n(\theta^*;x)+\left(\theta-\theta^*\right)^{\top} \nabla M_n\left(\theta^* ; x\right)- 
  \frac{1}{2} \left\|D_0\left(\theta-\theta^*\right)\right\|^2 + R\right]\mathrm{d} \theta\\
  &=& M_n\left(\theta^* ; x\right)-M_n\left(\theta^* ; x\right)\int_{\Theta_0^c\left(r_0\right)} q(\theta ;\mathbf{m},\mathbf{v})\mathrm{d} \theta+\sqrt{n}\left(\mathbf{m}-\theta^*\right)^{\top}V_{\theta_0}\Delta_{n, \theta_0}\\
  & & -\sqrt{n}\Delta_{n, \theta_0}^{\top}V_{\theta_0}\int_{\Theta_0^c\left(r_0\right)} q(\theta ;\mathbf{m},\mathbf{v})\left(\theta-\theta^*\right)\mathrm{d} \theta-\frac{1}{2}\mathrm{tr}\left(V_{\theta_0}\mathrm{diag}\left(\mathbf{v}\right)\right) -\frac{n}{2}\left(\mathbf{m}-\theta^*\right)^{\top}V_{\theta_0}\left(\mathbf{m}-\theta^*\right)\\
  & & +\frac{n}{2}\int_{\Theta_0^c\left(r_0\right)} q(\theta ;\mathbf{m},\mathbf{v})\left(\theta-\theta^*\right)^{\top}V_{\theta_0}\left(\theta-\theta^*\right)\mathrm{d} \theta+\int_{\Theta_0\left(r_0\right)} q(\theta ;\mathbf{m},\mathbf{v})R\mathrm{d} \theta\\
  &=& M_n\left(\theta^* ; x\right) - \frac{1}{2} \mathrm{tr}\left(V_{\theta_0}\mathrm{diag}\left(\mathbf{v}\right)\right) - \frac{n}{2}\left(\mathbf{m}-\theta^*\right)^{\top}V_{\theta_0}\left(\mathbf{m}-\theta^*\right)+\sqrt{n}\left(\mathbf{m}-\theta^*\right)^{\top}V_{\theta_0}\Delta_{n, \theta_0}+R\\
  & & -R\int_{\Theta_0^c\left(r_0\right)} q(\theta ;\mathbf{m},\mathbf{v})\mathrm{d} \theta-\sqrt{n}\Delta_{n, \theta_0}^{\top}V_{\theta_0}\int_{\Theta_0^c\left(r_0\right)} q(\theta ;\mathbf{m},\mathbf{v})\left(\theta-\theta^*\right)\mathrm{d} \theta\\
  & & +\frac{n}{2}\int_{\Theta_0^c\left(r_0\right)} q(\theta ;\mathbf{m},\mathbf{v})\left(\theta-\theta^*\right)^{\top}V_{\theta_0}\left(\theta-\theta^*\right)\mathrm{d} \theta-M_n\left(\theta^* ; x\right)\int_{\Theta_0^c\left(r_0\right)} q(\theta ;\mathbf{m},\mathbf{v})\mathrm{d} \theta.  
\end{eqnarray*}
Therefore, $\int q(\theta;\mathbf{m},\mathbf{v}) M_n(\theta; x) \mathrm{d} \theta$ can be rewritten as

\begin{eqnarray*}
\int q(\theta ;\mathbf{m},\mathbf{v}) M_n(\theta ; x) \mathrm{d} \theta 
&=&  
\int_{\Theta_0\left(r_0\right)} q(\theta ;\mathbf{m},\mathbf{v}) M_n(\theta ; x) \mathrm{d} \theta + \int_{\Theta_0^c\left(r_0\right)} q(\theta ;\mathbf{m},\mathbf{v}) M_n(\theta ; x) \mathrm{d} \theta\\
 &=& M_n\left(\theta^* ; x\right) -\frac{1}{2}\mathrm{tr}\left(V_{\theta_0}\mathrm{diag}\left(\mathbf{v}\right)\right) -\frac{n}{2}\left(\mathbf{m}-\theta^*\right)^{\top}V_{\theta_0}\left(\mathbf{m}-\theta^*\right) \\
 & & +\sqrt{n}\left(\mathbf{m}-\theta^*\right)^{\top}V_{\theta_0}\Delta_{n, \theta_0}+ R - R\int_{\Theta_0^c\left(r_0\right)} q(\theta ;\mathbf{m},\mathbf{v})\mathrm{d} \theta\\
 & & +  \int_{\Theta_0^c\left(r_0\right)} q(\theta ;\mathbf{m},\mathbf{v}) M_n(\theta ; x) \mathrm{d} \theta - M_n\left(\theta^* ; x\right)\int_{\Theta_0^c\left(r_0\right)} q(\theta ;\mathbf{m},\mathbf{v})\mathrm{d} \theta\\
  & & -\sqrt{n}\Delta_{n, \theta_0}^{\top}V_{\theta_0}\int_{\Theta_0^c\left(r_0\right)} q(\theta ;\mathbf{m},\mathbf{v})\left(\theta-\theta^*\right)\mathrm{d} \theta\\
  & & +\frac{n}{2}\int_{\Theta_0^c\left(r_0\right)} q(\theta ;\mathbf{m},\mathbf{v})\left(\theta-\theta^*\right)^{\top}V_{\theta_0}\left(\theta-\theta^*\right)\mathrm{d} \theta.  
\end{eqnarray*}
We now analyze the remaining four terms. By a Taylor expansion, there exists some $\bar{\theta}$  between $\theta$ and $\theta^*$ such that these terms satisfy: 
\begin{eqnarray}
\label{combination of edge and central}
& & \int_{\Theta_0^c\left(r_0\right)} q(\theta ;\mathbf{m},\mathbf{v}) M_n(\theta ; x) \mathrm{d} \theta-M_n\left(\theta^* ; x\right)\int_{\Theta_0^c\left(r_0\right)} q(\theta ;\mathbf{m},\mathbf{v})\mathrm{d} \theta \nonumber\\
& & -\sqrt{n}\Delta_{n, \theta_0}^{\top}V_{\theta_0}\int_{\Theta_0^c\left(r_0\right)} q(\theta ;\mathbf{m},\mathbf{v})\left(\theta-\theta^*\right)\mathrm{d} \theta+\frac{n}{2}\int_{\Theta_0^c\left(r_0\right)} q(\theta ;\mathbf{m},\mathbf{v})\left(\theta-\theta^*\right)^{\top}V_{\theta_0}\left(\theta-\theta^*\right)\mathrm{d} \theta \nonumber\\
&=& \int_{\Theta_0^c\left(r_0\right)} q(\theta ;\mathbf{m},\mathbf{v})\left[M_n(\theta ; x)-M_n\left(\theta^* ; x\right)-\left(\theta-\theta^*\right)^{\top}\nabla M_n\left(\theta^* ; x\right)\right]\mathrm{d} \theta\nonumber\\
& & -\frac{1}{2}\int_{\Theta_0^c\left(r_0\right)} q(\theta ;\mathbf{m},\mathbf{v})\left(\theta-\theta^*\right)^{\top}\nabla^2 \mathbb{E}_{\theta_0} M_n\left(\theta^*; x\right)\left(\theta-\theta^*\right)\mathrm{d} \theta \nonumber\\
&=& \frac{1}{2}\int_{\Theta_0^c\left(r_0\right)} q(\theta ;\mathbf{m},\mathbf{v})\left(\theta-\theta^*\right)^{\top}\nabla^2M_n\left(\bar{\theta}; x\right)\left(\theta-\theta^*\right)\mathrm{d} \theta \nonumber\\
& & -\frac{1}{2}\int_{\Theta_0^c\left(r_0\right)} q(\theta ;\mathbf{m},\mathbf{v})\left(\theta-\theta^*\right)^{\top}\nabla^2 \mathbb{E}_{\theta_0} M_n\left(\theta^*; x\right)\left(\theta-\theta^*\right)\mathrm{d} \theta \nonumber\\
&=& \frac{1}{2}\int_{\Theta_0^c\left(r_0\right)} q(\theta ;\mathbf{m},\mathbf{v})\left(\theta-\theta^*\right)^{\top} \left[\nabla^2\mathbb{E}_{\theta_0} M_n\left(\bar{\theta}; x\right)-\nabla^2\mathbb{E}_{\theta_0} M_n\left(\theta^*; x\right)\right] \left(\theta-\theta^*\right)\mathrm{d} \theta \nonumber\\
& & +\frac{1}{2} \int_{\Theta_0^c\left(r_0\right)} q(\theta ;\mathbf{m},\mathbf{v})\left(\theta-\theta^*\right)^{\top} \left[\nabla^2M_n\left(\bar{\theta}; x\right)-\nabla^2 \mathbb{E}_{\theta_0} M_n\left(\bar{\theta}; x\right)\right] \left(\theta-\theta^*\right)\mathrm{d} \theta \nonumber\\
&\triangleq& \frac{1}{2} U_1 + \frac{1}{2} U_2.
\end{eqnarray}
We next analyze the terms $U_1$ and $U_2$ one by one.

For the term $U_1$, by the smoothness assumption for the variational log-likelihood in Condition \ref{condition: lipschitz}, we have that
\begin{eqnarray}
\label{bar theta 1}
    |U_1| &\lesssim& n\int_{\Theta_0^c\left(r_0\right)} q(\theta ;\mathbf{m},\mathbf{v})\left\|\theta-\theta^*\right\|^2\left\|\bar{\theta}-\theta^*\right\|^{s_1}\mathrm{d} \theta \nonumber\\
    &\leq&n\int_{\Theta_0^c\left(r_0\right)} q(\theta ;\mathbf{m},\mathbf{v})\left\|\theta-\theta^*\right\|^{s_1+2}\mathrm{d} \theta \nonumber\\
   &=& n^{-\frac{s_1}{2}}\int_{\left\{\left\|V_{\theta_0}^{\frac{1}{2}} h\right\|>r_0\right\}}N\left(h ;\sqrt{n}\left(\mathbf{m}-\theta^*\right),\mathrm{diag}\left(\mathbf{v}\right)\right)\left\|h\right\|^{s_1+2}\mathrm{d}h.
\end{eqnarray}
Furthermore, the triangle inequality states that
\begin{equation*}
\left\|V_{\theta_0}^{\frac{1}{2}} h\right\|^2 \leq 2\left(\left\|V_{\theta_0}^{\frac{1}{2}} \left[h-\sqrt{n}\left(\mathbf{m}-\theta^*\right)\right]\right\|^2 + \left\|V_{\theta_0}^{\frac{1}{2}} \sqrt{n}\left(\mathbf{m}-\theta^*\right)\right\|^2\right).
\end{equation*}
And \cite{spokoiny2014bernsteinvonmises} claimed that $\delta\left(r_0\right) \asymp r_0 / \sqrt{n}$ under regular i.i.d. cases. Therefore, to constrain the search range of $\mathbf{m}$, we can introduce an independent constant $\beta\in\left(0,\frac{1}{\sqrt{2}}\right)$ such that
\begin{equation}
\label{range of m}
\left\|V_{\theta_0}^{\frac{1}{2}} \left(\mathbf{m}-\theta^*\right)\right\|\leq\beta\frac{r_0}{\sqrt{n}}=O\left(\delta\left(r_0\right)\right).
\end{equation}
It implies that there exists a constant $C_1>0$ such that 
\begin{equation*}
\left\{\left\|V_{\theta_0}^{\frac{1}{2}} h\right\|>r_0\right\}\subseteq \left\{\left\|V_{\theta_0}^{\frac{1}{2}} \left[h-\sqrt{n}\left(\mathbf{m}-\theta^*\right)\right]\right\|>C_1r_0\right\}.
\end{equation*}
Based on these results, we have
\begin{eqnarray}
\label{bar theta 2}
& & \int_{\left\{\left\|V_{\theta_0}^{\frac{1}{2}} h\right\|>r_0\right\}}N\left(h ;\sqrt{n}\left(\mathbf{m}-\theta^*\right),\mathrm{diag}\left(\mathbf{v}\right)\right)\left\|h\right\|^{s_1+2}\mathrm{d}h \nonumber\\
&\leq& \int_{\left\{\left\|V_{\theta_0}^{\frac{1}{2}}\left[h-\sqrt{n}\left(\mathbf{m}-\theta^*\right)\right]\right\| >C_1r_0\right\}}N\left(h ;\sqrt{n}\left(\mathbf{m}-\theta^*\right),\mathrm{diag}\left(\mathbf{v}\right)\right)\left\|h\right\|^{s_1+2}\mathrm{d}h \nonumber\\ 
&\lesssim& \int_{\left\{\left\|V_{\theta_0}^{\frac{1}{2}}\left[h-\sqrt{n}\left(\mathbf{m}-\theta^*\right)\right]\right\| > C_1r_0\right\}} N\left(h ;\sqrt{n}\left(\mathbf{m}-\theta^*\right),\mathrm{diag}\left(\mathbf{v}\right)\right)\left\|h-\sqrt{n}\left(\mathbf{m}-\theta^*\right)\right\|^{s_1+2}\mathrm{d}h \nonumber\\
 & & + \int_{\left\{\left\|V_{\theta_0}^{\frac{1}{2}}\left[h-\sqrt{n}\left(\mathbf{m}-\theta^*\right)\right]\right\|>C_1r_0\right\}} N\left(h ;\sqrt{n}\left(\mathbf{m}-\theta^*\right),\mathrm{diag}\left(\mathbf{v}\right)\right)\left\|\sqrt{n}\left(\mathbf{m}-\theta^*\right)\right\|^{s_1+2}\mathrm{d}h \nonumber\\
 &=& \int_{\left\{\left\|V_{\theta_0}^{\frac{1}{2}}\mathrm{diag}\left(\mathbf{v}\right)^{\frac{1}{2}}l\right\|>C_1r_0\right\}}N\left(l ;0,I_p\right)\left\|\mathrm{diag}\left(\mathbf{v}\right)^{\frac{1}{2}}l\right\|^{s_1+2}\mathrm{d}l \nonumber\\
 & & +\left\|\sqrt{n}\left(\mathbf{m}-\theta^*\right)\right\|^{s_1+2}\int_{\left\{\left\|V_{\theta_0}^{\frac{1}{2}}\mathrm{diag}\left(\mathbf{v}\right)^{\frac{1}{2}}l\right\|>C_1r_0\right\}}N\left(l ;0,I_p\right)\mathrm{d}l\nonumber\\
 &\triangleq& U_{11} + U_{12}.
\end{eqnarray}
For sufficiently large $r_0^2/p$, there exists a constant $\alpha$, independent of other parameters, such that
\begin{equation}
\label{assumption for alpha C5}
\frac{p+s_1+1}{C_1^2r_0^2}\rho\left(V_{\theta_0}\mathrm{diag}\left(\mathbf{v}\right)\right)\leq\alpha<1.
\end{equation}
Hence, by suitably restricting the search space of $\mathbf{v}$, as in Assumption 2 of Appendix C in \cite{wang2019frequentist} (which supports the proof of Lemma 3 therein), we have
\begin{eqnarray}
\label{bar theta 3}
 U_{11} &\leq& \rho\left(\mathrm{diag}\left(\mathbf{v}\right)\right)^{\frac{s_1+2}{2}}\int_{\left\{\left\|l\right\|>\frac{C_1r_0}{\sqrt{\rho\left(V_{\theta_0}\mathrm{diag}\left(\mathbf{v}\right)\right)}}\right\}}N\left(l;0,I_p\right)\left\|l\right\|^{s_1+2}\mathrm{d}l \nonumber\\
&\lesssim& \rho\left(\mathrm{diag}\left(\mathbf{v}\right)\right)^{\frac{s_1+2}{2}} \frac{\left[\frac{C_1^2r_0^2}{\rho\left(V_{\theta_0}\mathrm{diag}\left(\mathbf{v}\right)\right)}\right]^{\frac{p+s_1}{2}}}{(p-2)!! \exp\left\{\frac{C_1^2r_0^2}{2\rho\left(V_{\theta_0}\mathrm{diag}\left(\mathbf{v}\right)\right)}\right\}},
\end{eqnarray}
and
\begin{eqnarray}
\label{bar theta 4}
U_{12} &\leq& \left\|\sqrt{n}\left(\mathbf{m}-\theta^*\right)\right\|^{s_1+2}\int_{\left\{\left\|l\right\|>\frac{C_1r_0}{\sqrt{\rho\left(V_{\theta_0}\mathrm{diag}\left(\mathbf{v}\right)\right)}}\right\}}N\left(l;0,I_p\right)\mathrm{d}l \nonumber\\
&\lesssim& \left\|\sqrt{n}\left(\mathbf{m}-\theta^*\right)\right\|^{s_1+2} \frac{\left[\frac{C_1^2r_0^2}{\rho\left(V_{\theta_0}\mathrm{diag}\left(\mathbf{v}\right)\right)}\right]^{\frac{p-2}{2}}}{(p-2)!! \exp\left\{\frac{C_1^2r_0^2}{2\rho\left(V_{\theta_0}\mathrm{diag}\left(\mathbf{v}\right)\right)}\right\}}.
\end{eqnarray}
Therefore, combining \eqref{bar theta 1},  \eqref{bar theta 2}, \eqref{bar theta 3}, and \eqref{bar theta 4}, we derive that
\begin{equation*}
    |U_1| \lesssim n^{-\frac{s_1}{2}}\left[\left\{\frac{\rho\left(\mathrm{diag}\left(\mathbf{v}\right)\right)C_1^2r_0^2}{\rho\left(V_{\theta_0}\mathrm{diag}\left(\mathbf{v}\right)\right)}\right\}^{\frac{s_1+2}{2}} + \left\|\sqrt{n}\left(\mathbf{m}-\theta^*\right)\right\|^{s_1+2}\right] \frac{\left[\frac{C_1^2r_0^2}{\rho\left(V_{\theta_0}\mathrm{diag}\left(\mathbf{v}\right)\right)}\right]^{\frac{p-2}{2}}}{(p-2)!!\exp\left\{\frac{C_1^2r_0^2}{2\rho\left(V_{\theta_0}\mathrm{diag}\left(\mathbf{v}\right)\right)}\right\}}.
\end{equation*}
Note that
\begin{equation*}
\rho\left(\mathrm{diag}\left(\mathbf{v}\right)\right)=\rho\left(V_{\theta_0}^{-1}V_{\theta_0}\mathrm{diag}\left(\mathbf{v}\right)\right)\leq\rho\left(V_{\theta_0}^{-1}\right)\rho\left(V_{\theta_0}\mathrm{diag}\left(\mathbf{v}\right)\right),
\end{equation*}
and the constraint on the search space for $\mathbf{m}$ in \eqref{range of m} ensures that
\begin{equation*}
\left\|\sqrt{n}\left(\mathbf{m}-\theta^*\right)\right\| \leq \frac{\beta r_0}{\sqrt{\lambda_{\min}\left(V_{\theta_0}\right)} }.
\end{equation*}
We further have
\begin{eqnarray*}
|U_1| &\lesssim& n^{-\frac{s_1}{2}}\left[\left(\frac{C_1^2r_0^2}{\lambda_{\min}\left(V_{\theta_0}\right)}\right)^{\frac{s_1+2}{2}}+\left(\frac{\beta^2r_0^2}{\lambda_{\min}\left(V_{\theta_0}\right)}\right)^{\frac{s_1+2}{2}}\right] \frac{\left[\frac{C_1^2r_0^2}{\rho\left(V_{\theta_0}\mathrm{diag}\left(\mathbf{v}\right)\right)}\right]^{\frac{p-2}{2}}}{(p-2)!!\exp\left\{\frac{C_1^2r_0^2}{2\rho\left(V_{\theta_0}\mathrm{diag}\left(\mathbf{v}\right)\right)}\right\}}\\
&\lesssim& n^{-\frac{s_1}{2}} \left\{\frac{\rho\left(V_{\theta_0}\mathrm{diag}\left(\mathbf{v}\right)\right)}{\lambda_{\min}\left(V_{\theta_0}\right)}\right\}^{\frac{s_1+2}{2}} \frac{\left[\frac{C_1^2r_0^2}{\rho\left(V_{\theta_0}\mathrm{diag}\left(\mathbf{v}\right)\right)}\right]^{\frac{p+s_1}{2}}}{(p-2)!!\exp\left\{\frac{C_1^2r_0^2}{2\rho\left(V_{\theta_0}\mathrm{diag}\left(\mathbf{v}\right)\right)}\right\}}\\
&\lesssim& \frac{\rho\left(V_{\theta_0}\mathrm{diag}\left(\mathbf{v}\right)\right)}{\lambda_{\min}\left(V_{\theta_0}\right)}\left\{\frac{\rho\left(V_{\theta_0}\mathrm{diag}\left(\mathbf{v}\right)\right)}{n\lambda_{\min}\left(V_{\theta_0}\right)}\right\}^{\frac{s_1}{2}} \frac{\left[\frac{C_1^2r_0^2}{\rho\left(V_{\theta_0}\mathrm{diag}\left(\mathbf{v}\right)\right)}\right]^{\frac{p+s_1}{2}}}{\sqrt{p-2}\left(\frac{p-2}{e}\right)^{\frac{p-2}{2}}\exp\left\{\frac{C_1^2r_0^2}{2\rho\left(V_{\theta_0}\mathrm{diag}\left(\mathbf{v}\right)\right)}\right\}},    
\end{eqnarray*}
where the Stirling's formula is employed in the last step. Consider the function
\begin{equation*}
f(x)=\frac{p+s_1}{2}\log x-\frac{1}{2}x.
\end{equation*}
It is easy to show that $f'(x) < 0$ when $(p+s_1+1)/x \leq\alpha<1$. Therefore, based on \eqref{assumption for alpha C5}, we derive that
\begin{equation*}
\frac{p+s_1}{2}\log\left(\frac{C_1^2r_0^2}{\rho\left(V_{\theta_0}\mathrm{diag}\left(\mathbf{v}\right)\right)}\right)-\frac{C_1^2r_0^2}{2\rho\left(V_{\theta_0}\mathrm{diag}\left(\mathbf{v}\right)\right)}\leq\frac{p+s_1}{2}\log\left(\frac{p+s_1+1}{\alpha}\right)-\frac{p+s_1+1}{2\alpha}
\end{equation*}
due to the monotonic decreasing property of $f(x)$. It implies that 
\begin{eqnarray*}
& & \log \left\{\frac{\left(\frac{C_1^2r_0^2}{\rho\left(V_{\theta_0}\mathrm{diag}\left(\mathbf{v}\right)\right)}\right)^{\frac{p+s_1}{2}}}{\sqrt{p-2}\left(\frac{p-2}{e}\right)^{\frac{p-2}{2}}\exp\left(\frac{C_1^2r_0^2}{2\rho\left(V_{\theta_0}\mathrm{diag}\left(\mathbf{v}\right)\right)}\right)}\right\} \\
&=& \frac{p+s_1}{2}\log\left(\frac{C_1^2r_0^2}{\rho\left(V_{\theta_0}\mathrm{diag}\left(\mathbf{v}\right)\right)}\right)-\frac{C_1^2r_0^2}{2\rho\left(V_{\theta_0}\mathrm{diag}\left(\mathbf{v}\right)\right)}-\frac{p-2}{2}\log(p-2)+\frac{p-2}{2}-\frac{1}{2}\log(p-2)\\
&\leq& \frac{p+s_1}{2}\log\left(\frac{p+s_1+1}{\alpha}\right)-\frac{p+s_1+1}{2\alpha}-\frac{p-2}{2}\log(p-2)+\frac{p-2}{2}-\frac{1}{2}\log(p-2)\\
&=& \left(1+\log\frac{1}{\alpha}-\frac{1}{\alpha}\right)\frac{p+s_1}{2}+\frac{s_1+1}{2}\log(p+s_1+1)+\frac{p-1}{2}\log\left(1+\frac{s_1+3}{p-2}\right)-\frac{s_1+2}{2}-\frac{1}{2\alpha} \\
&\leq& \left(1+\log\frac{1}{\alpha}-\frac{1}{\alpha}\right)\frac{p}{2}+\frac{s_1+1}{2}\log(p+s_1+1)+\left(1+\log\frac{1}{\alpha}-\frac{1}{\alpha}\right)\frac{s_1}{2}+\frac{s_1+3}{2}\frac{p-1}{p-2}-\frac{s_1+2}{2}-\frac{1}{2\alpha}.
\end{eqnarray*}
For sufficiently large $r_0^2/p$, there exists a constant $\alpha$ satisfies both \eqref{assumption for alpha C5} and 
\begin{equation*}
1+\log\frac{1}{\alpha}-\frac{1}{\alpha}<0.
\end{equation*}
Therefore, by noting that 
$1/\left(n\lambda_{\min}\left(V_{\theta_0}\right)\right)=\rho\left(V_{\theta_0}^{-1}\right)/n\lesssim 1/\sqrt{n}$, we can conclude that $|U_1|$ is exponentially negligible.

For the term $U_2$, we use the exponential moment condition \ref{condition: global, exponential moment, variational} to construct a Chernoff bound as in Appendix \ref{Proof of Lemma 4} for term $T_{12}$. Specifically, for any arbitrarily large constant $M$, any $k>0$, and a given $h\in\mathbb{R}^p$, we have
\begin{eqnarray*}
& & P_{\theta_0}\left(\lambda_{\max}\left(\nabla^2M_n\left(\overline{\theta}; x\right)-\nabla^2 \mathbb{E}_{\theta_0} M_n\left(\overline{\theta}; x\right)\right)>Mnp^k\right) \\
&=& P_{\theta_0}\left(\frac{1}{n\left\|h\right\|^2} h^\top\nabla^2\left[M_n\left(\overline{\theta}; x\right)-\mathbb{E}_{\theta_0}M_n\left(\overline{\theta}; x\right)\right]h > Mp^k\right)\\
&=& P_{\theta_0}\left(\frac{\lambda}{\tau_3}\frac{h^\top\nabla^2\left[M_n\left(\overline{\theta}; x\right)-\mathbb{E}_{\theta_0}M_n\left(\overline{\theta}; x\right)\right]h}{nh^\top V_{\theta_0}h}>\frac{\lambda M\left\|h\right\|^2p^k}{\tau_3 h^\top V_{\theta_0}h}\right)\\
&\leq& \exp\left(\frac{\tau_1^2}{2}\lambda^2-\frac{M\left\|h\right\|^2p^k}{\tau_3 h^\top V_{\theta_0}h}\lambda\right),
\end{eqnarray*}
for $\lambda\in\left[0,\delta_2(r)\right]$. By choosing $k$ sufficiently large such that $p^k$ dominates $\rho\left(V_{\theta_0}\right)$, we obtain that
\begin{equation*}
U_2 \lesssim np^k\int_{\Theta_0^c\left(r_0\right)} q(\theta;\mathbf{m},\mathbf{v})\left\|\theta-\theta^*\right\|^2\mathrm{d} \theta,
\end{equation*}
with probability at least $1-e^{-cp^{k^\prime}}$ for some $c>0$, where $k^\prime$ is also sufficiently large. Similarly, we can also establish that
\begin{equation*}
-U_2 \lesssim np^k\int_{\Theta_0^c\left(r_0\right)} q(\theta;\mathbf{m},\mathbf{v})\left\|\theta-\theta^*\right\|^2\mathrm{d} \theta,
\end{equation*}
with probability at least $1-e^{-cp^{k^\prime}}$. Therefore,
\begin{equation*}
    |U_2| \lesssim np^k\int_{\Theta_0^c\left(r_0\right)} q(\theta;\mathbf{m},\mathbf{v})\left\|\theta-\theta^*\right\|^2\mathrm{d} \theta.
\end{equation*}
Compare the above inequality with that in \eqref{bar theta 1}, we can use the same technique to show that $|U_2|$ also decays exponentially. 

Thus far, we have demonstrated that the terms in \eqref{combination of edge and central} decay exponentially. This implies that these terms can be absorbed into the big-$O$ terms specified in \eqref{q m v logp} and \eqref{The Bracket Results of log}. Incorporating this result into \eqref{def of Fn}--\eqref{The Bracket Results of log}, and noting $R\in\left[-\Delta(r_0,y),\Delta(r_0,y)\right]$, we derive that
\begin{eqnarray}
\label{Fn=F0+error}
F_n(\mathbf{m},\mathbf{v})
&=& -\frac{p}{2}-\frac{1}{2}\sum\limits_{i=1}^p\log v_i-\frac{1}{2}\log\det\left(V_{\theta_0}\right)+\frac{1}{2} \mathrm{tr}\left(V_{\theta_0}\mathrm{diag}\left(\mathbf{v}\right)\right) +\log p(\theta^*)-\log p(\mathbf{m}) \nonumber\\
& & +\frac{n}{2}\left(\mathbf{m}-\theta^*-\frac{\Delta_{n, \theta_0}}{\sqrt{n}}\right)^{\top}V_{\theta_0}\left(\mathbf{m}-\theta^*-\frac{\Delta_{n, \theta_0}}{\sqrt{n}}\right) \nonumber \\
& & + O\left(\Delta(r_0,y)+\frac{1}{\sqrt{n}} \left\{\left\|\Delta_{n,\theta_0}\right\|^2+\left\|\nabla\log p\left(\theta^*\right)\right\|^2\right\} +\frac{1}{n} \left\{\mathrm{tr}\left(V_{\theta_0}^{-1}\right)+\mathrm{tr}\left(\mathrm{diag}\left(\mathbf{v}\right)\right)\right\} \right) \nonumber\\
&=& F_0(\mathbf{m},\mathbf{v})+\log p(\theta^*)-\log p(\mathbf{m}) \\
& & + O\left(\Delta(r_0,y)+\frac{1}{\sqrt{n}} \left\{\left\|\Delta_{n,\theta_0}\right\|^2+\left\|\nabla\log p\left(\theta^*\right)\right\|^2\right\} +\frac{1}{n} \left\{\mathrm{tr}\left(V_{\theta_0}^{-1}\right)+\mathrm{tr}\left(\mathrm{diag}\left(\mathbf{v}\right)\right)\right\} \right) \nonumber
\end{eqnarray}
which holds with probability at least $1-5e^{-y}-e^{-cr_0^2}$ for some constant $c>0$. All subsequent arguments are made on the random set under which  \eqref{Fn=F0+error} holds with probability at least $1-5e^{-y}-e^{-cr_0^2}$; thus, we omit restating this probability in what follows.

Denote $\mathbf{v}_0^* =\left(v_{01},\ldots,v_{0p}\right)^{\top}$ and $\mathbf{v}_n^*=\left(v_{n1},\ldots,v_{np}\right)^{\top}$. By the definitions of $(\mathbf{m}_n^*-\mathbf{m}_0^*)$ and $(\mathbf{m}_0^*,\mathbf{v}_0^*)$, there exists a constant $C_2>0$ such that 
\begin{eqnarray*}
  0 &\leq&  \frac{1}{2}\sum\limits_{i=1}^p\left(\frac{v_{ni}}{v_{0i}}-\log\frac{v_{ni}}{v_{0i}}-1\right)+\frac{n}{2}\left(\mathbf{m}_n^*-\mathbf{m}_0^*\right)^{\top}V_{\theta_0}\left(\mathbf{m}_n^*-\mathbf{m}_0^*\right) \\
  &=& F_0(\mathbf{m}_n^*,\mathbf{v}_n^*)-F_0(\mathbf{m}_0^*,\mathbf{v}_0^*)\\
&\leq& F_n(\mathbf{m}_n^*,\mathbf{v}_n^*)-\log p(\theta^*)+\log p(\mathbf{m}_n^*)-F_0(\mathbf{m}_0^*,\mathbf{v}_0^*)\\
& & + C_2\left(\Delta(r_0,y)+\frac{1}{\sqrt{n}} \left\{\left\|\Delta_{n,\theta_0}\right\|^2+\left\|\nabla\log p\left(\theta^*\right)\right\|^2\right\} +\frac{1}{n} \left\{\mathrm{tr}\left(V_{\theta_0}^{-1}\right)+\mathrm{tr}\left(\mathrm{diag}\left(\mathbf{v}_n^*\right)\right)\right\} \right)\\
&\leq& F_n(\mathbf{m}_0^*,\mathbf{v}_0^*)-\log p(\theta^*)+\log p(\mathbf{m}_n^*)-F_0(\mathbf{m}_0^*,\mathbf{v}_0^*)\\
& & + C_2 \left(\Delta(r_0,y)+\frac{1}{\sqrt{n}} \left\{\left\|\Delta_{n,\theta_0}\right\|^2+\left\|\nabla\log p\left(\theta^*\right)\right\|^2\right\} +\frac{1}{n} \left\{\mathrm{tr}\left(V_{\theta_0}^{-1}\right)+\mathrm{tr}\left(\mathrm{diag}\left(\mathbf{v}_n^*\right)\right)\right\} \right)
\end{eqnarray*}
holds based on \eqref{Fn=F0+error}. Further, recall that 
\[R_n = \Delta(r_0,y)+\frac{1}{\sqrt{n}} \left\{\left\|\Delta_{n,\theta_0}\right\|^2+\left\|\nabla\log p\left(\theta^*\right)\right\|^2\right\} +\frac{1}{n} \left\{\mathrm{tr}\left(V_{\theta_0}^{-1}\right)+\mathrm{tr}\left(\mathrm{diag}\left(\mathbf{v}_0^*\right)\right)\right\}.\]
Then, there exists a constant $C_3>0$ such that 
\begin{eqnarray}
  \label{sandwich bound of v m}
   0 &\leq&  \frac{1}{2}\sum\limits_{i=1}^p\left(\frac{v_{ni}}{v_{0i}}-\log\frac{v_{ni}}{v_{0i}}-1\right)+\frac{n}{2}\left(\mathbf{m}_n^*-\mathbf{m}_0^*\right)^{\top}V_{\theta_0}\left(\mathbf{m}_n^*-\mathbf{m}_0^*\right) \nonumber\\
   &\leq& \log p(\mathbf{m}_n^*)-\log p(\mathbf{m}_0^*) + C_3R_n\\
   & & + C_2 \left(\Delta(r_0,y)+\frac{1}{\sqrt{n}} \left\{\left\|\Delta_{n,\theta_0}\right\|^2+\left\|\nabla\log p\left(\theta^*\right)\right\|^2\right\} +\frac{1}{n} \left\{\mathrm{tr}\left(V_{\theta_0}^{-1}\right)+\mathrm{tr}\left(\mathrm{diag}\left(\mathbf{v}_n^*\right)\right)\right\} \right) \nonumber.
\end{eqnarray}
Under Condition \ref{condition: prior}, we have
\begin{equation*}
\label{log p theta* - log p m}
\left|\log p(\mathbf{m}_n^*)-\log p(\mathbf{m}_0^*)-\left(\mathbf{m}_n^*-\mathbf{m}_0^*\right)^{\top}\nabla\log p(\mathbf{m}_0^*)\right|\leq\frac{\tau_5}{2}\left\|\mathbf{m}_n^*-\mathbf{m}_0^*\right\|^2.
\end{equation*}
Consequently, there exists a constant $C_4>0$ such that
\begin{eqnarray*}
0 &\leq& \frac{1}{2}\sum\limits_{i=1}^p\left(\frac{v_{ni}}{v_{0i}}-\log\frac{v_{ni}}{v_{0i}}-1\right)+\frac{n}{2}\left(\mathbf{m}_n^*-\mathbf{m}_0^*\right)^{\top}V_{\theta_0}\left(\mathbf{m}_n^*-\mathbf{m}_0^*\right)\\
&\leq& \left(\mathbf{m}_n^*-\mathbf{m}_0^*\right)^{\top}\nabla\log p(\mathbf{m}_0^*) +\frac{\tau_5}{2}\left\|\mathbf{m}_n^*-\mathbf{m}_0^*\right\|^2 + C_4 \left( R_n +\frac{1}{n} \mathrm{tr}\left(\mathrm{diag}\left(\mathbf{v}_n^*\right)\right)\right).
\end{eqnarray*}

Now, we are going to show the conclusion in Theorem \ref{theorem: the convergence rates of the VBBE and the covariance matrix of the VB bridge}. The following overall derivation is organized into three steps.

\textit{Step 1.} We want to show that
\begin{equation}
\label{rate of v and m}
 \sum\limits_{i=1}^p\left(\frac{v_{ni}}{v_{0i}}-\log\frac{v_{ni}}{v_{0i}}-1\right)+n\left(\mathbf{m}_n^*-\mathbf{m}_0^*\right)^{\top}V_{\theta_0}\left(\mathbf{m}_n^*-\mathbf{m}_0^*\right)    \lesssim R_n + \frac{1}{n} \mathrm{tr}\left(\mathrm{diag}\left(\mathbf{v}_n^*\right)\right).
\end{equation}
To simplify notation, we denote
\begin{equation*}
    W_1 = \left(\mathbf{m}_n^*-\mathbf{m}_0^*\right)^{\top}\nabla\log p(\mathbf{m}_0^*),\quad
    W_2 = R_n + \frac{1}{n} \mathrm{tr}\left(\mathrm{diag}\left(\mathbf{v}_n^*\right)\right).
\end{equation*}
Note that the term $n\left(\mathbf{m}_n^*-\mathbf{m}_0^*\right)^{\top}V_{\theta_0}\left(\mathbf{m}_n^*-\mathbf{m}_0^*\right)$ dominates the term $\tau_5\left\|\mathbf{m}_n^*-\mathbf{m}_0^*\right\|^2$ due to $\lambda_{\min}\left(V_{\theta_0}\right)\gtrsim 1/\sqrt{n}$, . Therefore, it suffices to show that $W_1$ is dominated by $W_2$. We prove this by contradiction: suppose instead that $W_1$ dominates $W_2$. Then we have
\begin{equation*}
W_2 \lesssim \left\|\mathbf{m}_n^*-\mathbf{m}_0^*\right\|\cdot\left\|\nabla\log p(\mathbf{m}_0^*)\right\|,
\end{equation*}
which implies that
\begin{equation*}
n\lambda_{\min}\left(V_{\theta_0}\right)\left\|\mathbf{m}_n^*-\mathbf{m}_0^*\right\|^2\lesssim\left\|\mathbf{m}_n^*-\mathbf{m}_0^*\right\|\cdot\left\|\nabla\log p(\mathbf{m}_0^*)\right\|.
\end{equation*}
This suggests that
\begin{equation*}
\frac{W_2}{\left\|\nabla\log p(\mathbf{m}_0^*)\right\|} \lesssim\left\|\mathbf{m}_n^*-\mathbf{m}_0^*\right\| \lesssim \frac{\left\|\nabla\log p(\mathbf{m}_0^*)\right\|}{n\lambda_{\min}\left(V_{\theta_0}\right)}.
\end{equation*}
Hence, 
\begin{equation}
\label{proof by contradiction}
n\lambda_{\min}\left(V_{\theta_0}\right)W_2 \lesssim \left\|\nabla\log p(\mathbf{m}_0^*)\right\|^2.
\end{equation}
However, by definition, we have
\begin{equation*}
\left\|\nabla\log p(\mathbf{m}_0^*)\right\|^2=\left\|\nabla\log p\left(\theta^*+\frac{\Delta_{n, \theta_0}}{\sqrt{n}}\right)\right\|^2.
\end{equation*}
It is dominated by $n\lambda_{\min}\left(V_{\theta_0}\right)W_2$ by Condition \ref{condition: prior}. This contradicts conclusion \eqref{proof by contradiction}. Further, note that $n\lambda_{\min}\left(V_{\theta_0}\right)W_2 \gtrsim \sqrt{n}W_2$ due to  $\lambda_{\min}\left(V_{\theta_0}\right)\gtrsim 1/\sqrt{n}$. Therefore, we arrive at conclusion \eqref{rate of v and m}.

{\textit{Step 2.}} We want show that
\begin{equation}
 \label{final rate of ddagger}
 \sum\limits_{i=1}^p\left(\frac{v_{ni}}{v_{0i}}-\log\frac{v_{ni}}{v_{0i}}-1\right)+n\left(\mathbf{m}_n^*-\mathbf{m}_0^*\right)^{\top}V_{\theta_0}\left(\mathbf{m}_n^*-\mathbf{m}_0^*\right) \lesssim R_n
\end{equation}
It suffices to show that the term $\mathrm{tr}\left(\mathrm{diag}\left(\mathbf{v}_n^*\right)\right)/n$ can be absorbed by other terms on the right side of \eqref{rate of v and m}. First, note that the discrepancy between $\mathrm{tr}\left(\mathrm{diag}\left(\mathbf{v}_n^*\right)\right)$ and $\mathrm{tr}\left(\mathrm{diag}\left(\mathbf{v}_0^*\right)\right)$ can be controlled as follows:
\begin{equation*}
\left|\mathrm{tr}\left(\mathrm{diag}\left(\mathbf{v}_n^*\right)\right)-\mathrm{tr}\left(\mathrm{diag}\left(\mathbf{v}_0^*\right)\right)\right|=\left|\sum\limits_{i=1}^p v_{0i}\left(\frac{v_{ni}}{v_{0i}}-1\right)\right|\leq\sqrt{\sum\limits_{i=1}^p v_{0i}^2}\sqrt{\sum\limits_{i=1}^p\left(\frac{v_{ni}}{v_{0i}}-1\right)^2}.
\end{equation*}
Then, we can show that
\begin{equation}
\label{vni/voi-1}
\sum\limits_{i=1}^p\left(\frac{v_{ni}}{v_{0i}}-1\right)^2 \lesssim R_n +\frac{1}{n}\left\|\mathbf{v}_0^*\right\|\sqrt{\sum\limits_{i=1}^p\left(\frac{v_{ni}}{v_{0i}}-1\right)^2}
\end{equation}
by the facts that $n\left(\mathbf{m}_n^*-\mathbf{m}_0^*\right)^{\top}V_{\theta_0}\left(\mathbf{m}_n^*-\mathbf{m}_0^*\right)$ is non-negative, and $x-\log x-1  \asymp (x-1)^2$ when $x \rightarrow 1$. Using a similar contradiction-based argument as above, we can show that $R_n$ dominates $ n^{-1} \left\|\mathbf{v}_0^*\right\| \sqrt{\sum\limits_{i=1}^p\left(\frac{v_{ni}}{v_{0i}}-1\right)^2}$. If it is wrong, that is,
\begin{equation*}
    R_n \lesssim \frac{1}{n}\left\|\mathbf{v}_0^*\right\|\sqrt{\sum\limits_{i=1}^p\left(\frac{v_{ni}}{v_{0i}}-1\right)^2},
\end{equation*}
we can infer that $R_n \lesssim \left\|\mathbf{v}_0^*\right\|^2 / {n^2}$ under \eqref{vni/voi-1}. It leads to a contradiction regarding the orders of $n$ and $p$. Therefore, we arrive at conclusion \eqref{final rate of ddagger}.

{\textit{Step 3.}} Substitute
\begin{equation*}
\left(\theta^*+\frac{\Delta_{n, \theta_0}}{\sqrt{n}},\mathrm{diag}(\mathrm{diag}^{-1}(V_{\theta_0}))\right)=(\mathbf{m}_0^*,\mathrm{diag}(\mathbf{v}_0^*)),\quad(\hat{\theta}_n^\ddagger, \hat{V}_n^\ddagger)=(\mathbf{m}_n^*,\mathrm{diag}(\mathbf{v}_n^*))
\end{equation*}
into \eqref{final rate of ddagger}. It becomes clear that when $r_0^2/p$ is sufficiently large, the convergence rate of $\mathbf{m}_n^*$ matches the expected range in \eqref{range of m}, supporting the conclusions of Theorem \ref{theorem: the convergence rates of the VBBE and the covariance matrix of the VB bridge}.

\section{Proof of Theorem \ref{theorem: the consistency of the VB bridge}}
\label{Proof of Theorem 2}

We utilize the proof techniques from Appendix B in \cite{wang2019frequentist} to bound the integral $\int q^{\ddagger,\Theta_0\left(r_0\right)}(\theta) M_n(\theta ; x) \mathrm{d} \theta$, focusing on the concentration of $q^{\ddagger,\Theta_0\left(r_0\right)}(\theta)$ around $\theta^*$. A tail technique then extends this result to $q^\ddagger(\theta)$. Overall, the proof will be carried out in three main steps.

{\textit{Step 1.}} Renormalize $q^{\ddagger}(\theta)$ over $\Theta_0\left(r_0\right)$ as 
\begin{equation*}
q^{\ddagger,\Theta_0(r_0)}(\theta) = \frac{q^{\ddagger}(\theta) I_{\Theta_0(r_0)}(\theta)}{\int q^{\ddagger}(\theta) I_{\Theta_0(r_0)}(\theta) \mathrm{d} \theta},
\end{equation*}
and establish both upper and lower bounds for $\int q^{\ddagger,\Theta_0\left(r_0\right)}(\theta) M_n(\theta ; x) \mathrm{d} \theta$.

Recall that $h=\sqrt{n}\left(\theta-\theta^*\right)$, by Lemma \ref{LAN}, we have
\begin{eqnarray*}
& & \int q^{\ddagger,\Theta_0\left(r_0\right)}(\theta) M_n(\theta ; x) \mathrm{d} \theta\\
&=& \int q^{\ddagger,\Theta_0\left(r_0\right)}(\theta) M_n\left(\theta^*+\frac{h}{\sqrt{n}}; x\right) \mathrm{d} \theta\\
&\leq& \int q^{\ddagger,\Theta_0\left(r_0\right)}(\theta)\left[M_n(\theta^*;x)+h^{\top}V_{\theta_0}\Delta_{n, \theta_0}-\frac{1}{2}h^{\top}V_{\theta_0}h +\Delta(r_0,y)\right]\mathrm{d} \theta\\
&=& M_n(\theta^*;x)+\Delta(r_0,y)+\int q^{\ddagger,\Theta_0\left(r_0\right)}(\theta)\left(h^{\top}V_{\theta_0}\Delta_{n, \theta_0}-\frac{1}{2}h^{\top}V_{\theta_0}h \right)\mathrm{d} \theta
\end{eqnarray*}
holds with probability at least $ 1 - e^{-y}$. The last term can be split into two parts:
\begin{eqnarray*}
& & \int q^{\ddagger,\Theta_0\left(r_0\right)}(\theta)\left(h^{\top}V_{\theta_0}\Delta_{n, \theta_0}-\frac{1}{2} h^{\top}V_{\theta_0}h \right)\mathrm{d} \theta\\
&=& \int_{B\left(\theta^*, \eta\right)} q^{\ddagger,\Theta_0\left(r_0\right)}(\theta)\left(h^{\top}V_{\theta_0}\Delta_{n, \theta_0}-\frac{1}{2}h^{\top}V_{\theta_0}h \right)\mathrm{d} \theta\\
& & +\int_{B^c\left(\theta^*, \eta\right)} q^{\ddagger,\Theta_0\left(r_0\right)}(\theta)\left(h^{\top}V_{\theta_0}\Delta_{n, \theta_0}-\frac{1}{2} h^{\top}V_{\theta_0}h\right)\mathrm{d} \theta\\
&\triangleq& A_1 + A_2.
\end{eqnarray*}
For the term $A_1$, we have
\begin{eqnarray*}
 A_1 &\leq& \int_{B\left(\theta^*, \eta\right)} \frac{1}{2}\Delta_{n, \theta_0}^{\top}V_{\theta_0}\Delta_{n, \theta_0} q^{\ddagger,\Theta_0\left(r_0\right)}(\theta)\mathrm{d} \theta\\
 &=& \frac{1}{2}\Delta_{n, \theta_0}^{\top}V_{\theta_0}\Delta_{n, \theta_0}\left(1-\int_{B^c\left(\theta^*, \eta\right)} q^{\ddagger,\Theta_0\left(r_0\right)}(\theta)\mathrm{d} \theta\right).   
\end{eqnarray*}
For the term $A_2$, we know that
\begin{eqnarray*}
A_2 &=& \int_{\left\{\left\|h\right\|>\sqrt{n}\eta\right\}} q^{\ddagger,\Theta_0\left(r_0\right)}(\theta) \left(h^{\top}V_{\theta_0}\Delta_{n, \theta_0}-\frac{1}{2}h^{\top}V_{\theta_0}h\right)\mathrm{d} \theta\\
&\leq& - n\eta^2C_1 \lambda_{\min}\left(V_{\theta_0}\right) \int_{B^c\left(\theta^*, \eta\right)} q^{\ddagger,\Theta_0\left(r_0\right)}(\theta)\mathrm{d} \theta,    
\end{eqnarray*}
where $C_1\in\left(0,1/2\right)$ is a constant reflecting the dominance of the quadratic term $h^{\top}V_{\theta_0}h$ when $\left\|h\right\|>\sqrt{n}\eta$. By Condition \eqref{condition: identifiability}, we note that $\mathbb{E}_{\theta_0}\left\|\Delta_{n, \theta_0}\right\|^2\lesssim\mathrm{tr}\left(V_{\theta_0}^{-1}\right)=o\left(\left\|h\right\|^2\right)$, provided $\eta$ approaches zero at a sufficiently slow rate. Thus, we derive the following upper bound for $\int q^{\ddagger,\Theta_0\left(r_0\right)}(\theta) M_n(\theta ; x) \mathrm{d} \theta$:
\begin{eqnarray}
\label{upper bound of renormalized q ddagger M_n final}
& & \int q^{\ddagger,\Theta_0\left(r_0\right)}(\theta) M_n(\theta ; x) \mathrm{d} \theta \nonumber\\
&\leq & M_n(\theta^*;x)+\Delta(r_0,y)+\frac{1}{2}\Delta_{n, \theta_0}^{\top}V_{\theta_0}\Delta_{n, \theta_0} \left(1-\int_{B^c\left(\theta^*, \eta\right)} q^{\ddagger,\Theta_0\left(r_0\right)}(\theta)\mathrm{d} \theta\right) \nonumber\\
& & -C_1\lambda_{\min}\left(V_{\theta_0}\right)n\eta^2\int_{B^c\left(\theta^*, \eta\right)} q^{\ddagger,\Theta_0\left(r_0\right)}(\theta)\mathrm{d} \theta \nonumber\\
&=& M_n(\theta^*;x)+\Delta(r_0,y)+\frac{1}{2}\Delta_{n, \theta_0}^{\top}V_{\theta_0}\Delta_{n, \theta_0}\nonumber\\
& & - \int_{B^c\left(\theta^*, \eta\right)} q^{\ddagger,\Theta_0\left(r_0\right)}(\theta)\mathrm{d} \theta \left(\frac{1}{2}\Delta_{n, \theta_0}^{\top}V_{\theta_0}\Delta_{n, \theta_0} +C_1\lambda_{\min}\left(V_{\theta_0}\right)n\eta^2\right) 
\end{eqnarray}
holds with probability at least $ 1 - e^{-y}$.

We then derive a lower bound for $\int q^{\ddagger,\Theta_0\left(r_0\right)}(\theta) M_n(\theta ; x) \mathrm{d} \theta$. First, according to the definition of KL divergence, $\int q^{\ddagger,\Theta_0\left(r_0\right)}(\theta) M_n(\theta ; x) \mathrm{d} \theta$ can be rewritten as
\begin{eqnarray}
\label{Identity for the normalized q ddagger}
& & \int q^{\ddagger,\Theta_0\left(r_0\right)}(\theta) M_n(\theta ; x) \mathrm{d} \theta \nonumber\\
&=& \int q^{\ddagger,\Theta_0\left(r_0\right)}(\theta)\log q^{\ddagger,\Theta_0\left(r_0\right)}(\theta)\mathrm{d}\theta+\log \int p(\theta) \exp \left\{M_n(\theta ; x)\right\} \mathrm{d} \theta-\int q^{\ddagger,\Theta_0\left(r_0\right)}(\theta)\log p(\theta) \mathrm{d} \theta \nonumber\\
 & & -\operatorname{KL}\left(q^{\ddagger,\Theta_0\left(r_0\right)}(\theta)\| \pi^*(\theta \mid x)\right) \nonumber\\
 &=& \int q^{\ddagger,\Theta_0\left(r_0\right)}(\theta)\log q^{\ddagger,\Theta_0\left(r_0\right)}(\theta)\mathrm{d}\theta+\log \int p(\theta) \exp \left\{M_n(\theta ; x)\right\} \mathrm{d} \theta-\int q^{\ddagger,\Theta_0\left(r_0\right)}(\theta)\log p(\theta) \mathrm{d} \theta \nonumber\\
 & & - \operatorname{KL}\left(q^{\ddagger}(\theta)\| \pi^*(\theta \mid x)\right) + \left\{\operatorname{KL}\left(q^{\ddagger}(\theta)\| \pi^*(\theta \mid x)\right)-\operatorname{KL}\left(q^{\ddagger,\Theta_0\left(r_0\right)}(\theta)\| \pi^*(\theta \mid 
x)\right)\right\} \nonumber\\
&\triangleq& T_1 + T_2 - T_3 - T_4 + T_5.
\end{eqnarray}
Next, we analyze these terms one by one. For the term $T_1$, a straightforward calculation yields
\begin{eqnarray*}
 T_1 &=& \int_{\Theta_0\left(r_0\right)}\frac{q^{\ddagger}(\theta)}{\int_{\Theta_0\left(r_0\right)} q^{\ddagger}(\theta) \mathrm{d} \theta}\log\frac{q^{\ddagger}(\theta)}{\int_{\Theta_0\left(r_0\right)} q^{\ddagger}(\theta) \mathrm{d} \theta}\mathrm{d} \theta \\   &=& \frac{\int_{\Theta_0\left(r_0\right)} q^{\ddagger}(\theta)\log q^{\ddagger}(\theta)\mathrm{d} \theta}{\int_{\Theta_0\left(r_0\right)} q^{\ddagger}(\theta) \mathrm{d} \theta}-\log\int_{\Theta_0\left(r_0\right)} q^{\ddagger}(\theta) \mathrm{d} \theta \\
&=& -\frac{p}{2}-\frac{p}{2}\log(2\pi)+\frac{p}{2}\log n-\frac{1}{2}\sum\limits_{i=1}^p\log \hat{v}_{n,ii}^\ddagger \\
& & +\left(-\frac{p}{2}-\frac{p}{2}\log(2\pi)+\frac{p}{2}\log n-\frac{1}{2}\sum\limits_{i=1}^p\log \hat{v}_{n,ii}^\ddagger\right)\frac{\int_{\Theta_0^c\left(r_0\right)} q^{\ddagger}(\theta) \mathrm{d} \theta}{1-\int_{\Theta_0^c\left(r_0\right)} q^{\ddagger}(\theta) \mathrm{d} \theta} \\
& & -\frac{\int_{\Theta_0^c\left(r_0\right)} q^{\ddagger}(\theta)\log q^{\ddagger}(\theta) \mathrm{d} \theta}{1-\int_{\Theta_0^c\left(r_0\right)} q^{\ddagger}(\theta) \mathrm{d} \theta}-\log\left(1-\int_{\Theta_0^c\left(r_0\right)} q^{\ddagger}(\theta) \mathrm{d} \theta\right).
\end{eqnarray*}
Using the same technique as that applied to the term $U_1$ defined in \eqref{combination of edge and central}, we can similarly show that $\int_{\Theta_0^c\left(r_0\right)} q^{\ddagger}(\theta) \mathrm{d} \theta$ and $\int_{\Theta_0^c\left(r_0\right)} q^{\ddagger}(\theta)\log q^{\ddagger}(\theta) \mathrm{d} \theta$ are all exponentially small by Lemma \ref{lemma: The Tail Moment Control for the p-Dimensional Normal Distribution} and Theorem \ref{theorem: the convergence rates of the VBBE and the covariance matrix of the VB bridge}. Hence, 
\begin{equation}
\label{Identity for the first term}
T_1 =  -\frac{p}{2}-\frac{p}{2}\log(2\pi)+\frac{p}{2}\log n-\frac{1}{2}\sum\limits_{i=1}^p\log \hat{v}_{n,ii}^\ddagger  + T_{1n},
\end{equation}
where $T_{1n}$ is exponentially small. 
For the term $T_2$, by Lemma \ref{Lemma: the bracketing result of the logarithm of the denominator of the VB Ideal}, we have
\begin{eqnarray}
\label{Lower bound for the second term}
T_2 &=& \log\int p(\theta) \exp \left\{M_n(\theta ; x)\right\} \mathrm{d} \theta \nonumber\\
&>& M_n(\theta^*;x)+\log p\left(\theta^*\right)+\frac{p}{2}\log(2\pi)-\frac{1}{2}\log\det\left(V_{\theta_0}\right)-\frac{p}{2}\log n +\frac{1}{2}\Delta_{n,\theta_0}^{\top}V_{\theta_0}\Delta_{n,\theta_0} \nonumber\\ 
& &  -C_3\left(\Delta(r_0,y)+\frac{1}{\sqrt{n}} \left\{\left\|\Delta_{n,\theta_0}\right\|^2+\left\|\nabla\log p\left(\theta^*\right)\right\|^2\right\} +\frac{\mathrm{tr}\left(V_{\theta_0}^{-1}\right)}{n}\right)
\end{eqnarray}
holds for a certain $C_3>0$, with probability at least $1-5e^{-y}-e^{-cr_0^2}$ for some $c>0$. For the term $T_3$, a direct calculation shows that 
\begin{eqnarray*}
T_3 &=& \int_{\Theta_0\left(r_0\right)}\frac{q^{\ddagger}(\theta)}{\int_{\Theta_0\left(r_0\right)} q^{\ddagger}(\theta) \mathrm{d} \theta}\log p(\theta) \mathrm{d} \theta\\
&=& \frac{\int q^{\ddagger}(\theta)\log p(\theta) \mathrm{d} \theta-\int_{\Theta_0^c\left(r_0\right)} q^{\ddagger}(\theta)\log p(\theta) \mathrm{d} \theta}{1-\int_{\Theta_0^c\left(r_0\right)} q^{\ddagger}(\theta) \mathrm{d} \theta}\\
&=& \int q^\ddagger(\theta)\log p(\theta) \mathrm{d} \theta+\int q^\ddagger(\theta)\log p(\theta) \mathrm{d} \theta\frac{\int_{\Theta_0^c\left(r_0\right)} q^{\ddagger}(\theta) \mathrm{d} \theta}{1-\int_{\Theta_0^c\left(r_0\right)} q^{\ddagger}(\theta) \mathrm{d} \theta}-\frac{\int_{\Theta_0^c\left(r_0\right)} q^{\ddagger}(\theta)\log p(\theta) \mathrm{d} \theta}{1-\int_{\Theta_0^c\left(r_0\right)} q^{\ddagger}(\theta) \mathrm{d} \theta}.
\end{eqnarray*}
We first examine the last two terms in the expression above. As shown above, $\int_{\Theta_0^c\left(r_0\right)} q^{\ddagger}(\theta) \mathrm{d} \theta$ and $\int_{\Theta_0^c\left(r_0\right)} q^{\ddagger}(\theta)\log q^{\ddagger}(\theta) \mathrm{d} \theta$ are exponentially small using the technique applied to $U_1$ defined in \eqref{combination of edge and central}. Noting that $p(\theta)$ belongs to a Gaussian family, we can apply the same argument to show that $\int_{\Theta_0^c\left(r_0\right)} q^{\ddagger}(\theta)\log p(\theta) \mathrm{d} \theta$ is also exponentially small. Consequently,
\begin{equation*}
\int q^\ddagger(\theta)\log p(\theta) \mathrm{d} \theta\frac{\int_{\Theta_0^c\left(r_0\right)} q^{\ddagger}(\theta) \mathrm{d} \theta}{1-\int_{\Theta_0^c\left(r_0\right)} q^{\ddagger}(\theta) \mathrm{d} \theta}-\frac{\int_{\Theta_0^c\left(r_0\right)} q^{\ddagger}(\theta)\log p(\theta) \mathrm{d} \theta}{1-\int_{\Theta_0^c\left(r_0\right)} q^{\ddagger}(\theta) \mathrm{d} \theta}
\end{equation*}
is exponentially small. On the other hand, based on Condition \ref{condition: prior}, we obtain that
\begin{equation*}
\left|\int q^\ddagger(\theta)\log p(\theta) \mathrm{d} \theta-\log p\left(\hat{\theta}_n^\ddagger\right)\right|\leq\frac{\tau_5}{2n}\mathrm{tr}\left(\hat{V}_n^\ddagger\right).
\end{equation*}
Therefore,
\begin{equation}
\label{Upper bound for the third term}
T_3  \leq    \log p\left(\hat{\theta}_n^\ddagger\right) +\frac{\tau_5}{2n}\mathrm{tr}\left(\hat{V}_n^\ddagger\right) + T_{3n},
\end{equation}
 where $T_{3n}$ is exponentially small. For the term $T_4$, first note that $T_4 = F_n(\mathbf{m}_n^*,\mathbf{v}_n^*) \leq F_n(\mathbf{m}_0^*,\mathbf{v}_0^*)$ by the definition of function $F_n(\mathbf{m}_n,\mathbf{v}_n)$ defined in \eqref{def of Fn} and its minimizer $(\mathbf{m}_n^*,\mathbf{v}_n^*)$. As in Theorem \ref{theorem: the convergence rates of the VBBE and the covariance matrix of the VB bridge}, let $R_n = \Delta(r_0,y)+  \left\{\left\|\Delta_{n,\theta_0}\right\|^2+\left\|\nabla\log p\left(\theta^*\right)\right\|^2\right\}  / \sqrt{n} + \left\{\mathrm{tr}\left(V_{\theta_0}^{-1}\right)+ \mathrm{tr}\left(\mathrm{diag}\left(\mathbf{v}_0^*\right)\right)\right\} / n$. Then, by the relationship between $F_n(\mathbf{m}_n,\mathbf{v}_n)$
and $F_0(\mathbf{m}_n,\mathbf{v}_n)$ derived in \eqref{Fn=F0+error}, we can further show that, with probability at least $1-5e^{-y}-e^{-cr_0^2}$, the following inequality 
\begin{eqnarray}
  \label{Upper bound for the fourth term}
 T_4 &\leq&  F_0(\mathbf{m}_0^*,\mathbf{v}_0^*)+\log p(\theta^*)-\log p(\mathbf{m}_0^*) +C_3 R_n \nonumber\\
&=& \frac{1}{2}\sum\limits_{i=1}^p\log v_{\theta_0,ii}-\frac{1}{2}\log\det\left(V_{\theta_0}\right)+\log p(\theta^*)-\log p\left(\theta^*+\frac{\Delta_{n, \theta_0}}{\sqrt{n}}\right) + C_3R_n
\end{eqnarray}
holds for some $C_3>0$, where the final equality holds by the definition of $F_0(\mathbf{m}_n,\mathbf{v}_n)$ defined in \eqref{def of F0} and its minimizer $(\mathbf{m}_n^*,\mathbf{v}_n^*)$ defined in \eqref{def mv0}. For the fifth term $T_5$, it is easy to verify by direct computation that 
\begin{eqnarray*}
T_5 &=& \int q^{\ddagger}(\theta)\log\frac{q^{\ddagger}(\theta)}{\pi^*(\theta \mid x)}\mathrm{d}\theta-\int_{\Theta_0\left(r_0\right)}\frac{q^{\ddagger}(\theta)}{\int_{\Theta_0\left(r_0\right)} q^{\ddagger}(\theta) \mathrm{d} \theta}\log \left(\frac{\frac{q^{\ddagger}(\theta)}{\int_{\Theta_0\left(r_0\right)} q^{\ddagger}(\theta) \mathrm{d} \theta}}{\pi^*(\theta \mid x)} \right)\mathrm{d}\theta\\
&=& \int q^\ddagger(\theta)\log q^\ddagger(\theta)\mathrm{d}\theta-\int q^\ddagger(\theta)\log \pi^*(\theta \mid x)\mathrm{d}\theta-\frac{\int_{\Theta_0\left(r_0\right)}q^\ddagger(\theta)\log q^\ddagger(\theta)\mathrm{d}\theta}{\int_{\Theta_0\left(r_0\right)} q^{\ddagger}(\theta) \mathrm{d} \theta}\\
& & +\log\int_{\Theta_0\left(r_0\right)} q^{\ddagger}(\theta) \mathrm{d} \theta+\frac{\int_{\Theta_0\left(r_0\right)} q^\ddagger(\theta)\log \pi^*(\theta \mid x)\mathrm{d}\theta}{\int_{\Theta_0\left(r_0\right)} q^{\ddagger}(\theta) \mathrm{d} \theta}\\
&=& \int q^\ddagger(\theta)\log q^\ddagger(\theta)\mathrm{d}\theta-\int q^\ddagger(\theta)\log \pi^*(\theta \mid x)\mathrm{d}\theta-\frac{\int q^\ddagger(\theta)\log q^\ddagger(\theta)\mathrm{d}\theta-\int_{\Theta_0^c\left(r_0\right)}q^\ddagger(\theta)\log q^\ddagger(\theta)\mathrm{d}\theta}{1-\int_{\Theta_0^c\left(r_0\right)} q^{\ddagger}(\theta) \mathrm{d} \theta}\\
& & +\log\left(1-\int_{\Theta_0^c\left(r_0\right)} q^{\ddagger}(\theta) \mathrm{d} \theta\right)+\frac{\int q^\ddagger(\theta)\log \pi^*(\theta \mid x)\mathrm{d}\theta-\int_{\Theta_0^c\left(r_0\right)} q^\ddagger(\theta)\log \pi^*(\theta \mid x)\mathrm{d}\theta}{1-\int_{\Theta_0^c\left(r_0\right)} q^{\ddagger}(\theta) \mathrm{d} \theta}\\
&=& -\int q^\ddagger(\theta)\log q^\ddagger(\theta)\mathrm{d}\theta\frac{\int_{\Theta_0^c\left(r_0\right)} q^{\ddagger}(\theta) \mathrm{d} \theta}{1-\int_{\Theta_0^c\left(r_0\right)} q^{\ddagger}(\theta) \mathrm{d} \theta}+\int q^\ddagger(\theta)\log \pi^*(\theta \mid x)\mathrm{d}\theta\frac{\int_{\Theta_0^c\left(r_0\right)} q^{\ddagger}(\theta) \mathrm{d} \theta}{1-\int_{\Theta_0^c\left(r_0\right)} q^{\ddagger}(\theta) \mathrm{d} \theta}\\
& & +\frac{\int_{\Theta_0^c\left(r_0\right)}q^\ddagger(\theta)\log q^\ddagger(\theta)\mathrm{d}\theta}{1-\int_{\Theta_0^c\left(r_0\right)} q^{\ddagger}(\theta) \mathrm{d} \theta}+\log\left(1-\int_{\Theta_0^c\left(r_0\right)} q^{\ddagger}(\theta) \mathrm{d} \theta\right)-\frac{\int_{\Theta_0^c\left(r_0\right)} q^\ddagger(\theta)\log \pi^*(\theta \mid x)\mathrm{d}\theta}{1-\int_{\Theta_0^c\left(r_0\right)} q^{\ddagger}(\theta) \mathrm{d} \theta}.
\end{eqnarray*}
By Lemma \ref{LAN} and Condition \ref{condition: global identification property},  $\int_{\Theta_0^c\left(r_0\right)} q^\ddagger(\theta) M_n(\theta; x)\mathrm{d}\theta$ is exponentially small using the technique applied to $U_1$ defined in \eqref{combination of edge and central}. Consequently, 
\begin{eqnarray*}
& & \int_{\Theta_0^c\left(r_0\right)} q^\ddagger(\theta)\log \pi^*(\theta \mid x)\mathrm{d}\theta\\
&=& \int_{\Theta_0^c\left(r_0\right)} q^\ddagger(\theta)\log p(\theta)\mathrm{d}\theta+\int_{\Theta_0^c\left(r_0\right)} q^\ddagger(\theta) M_n(\theta ; x)\mathrm{d}\theta\\
& & - \log\left(\int p(\theta)\exp\left\{M_n(\theta ; x)\right\}\mathrm{d}\theta\right) \int_{\Theta_0^c\left(r_0\right)} q^\ddagger(\theta)\mathrm{d}\theta 
\end{eqnarray*}
is also exponentially small. Thus, we can conclude that the discrepancy $T_5$ is exponentially negligible. Combined with the results \eqref{Identity for the normalized q ddagger}--\eqref{Upper bound for the fourth term}, we establish the following lower bound for $\int q^{\ddagger,\Theta_0\left(r_0\right)}(\theta) M_n(\theta ; x) \mathrm{d} \theta$: with probability at least $1-5e^{-y}-e^{-cr_0^2}$,
\begin{eqnarray}
\label{lower bound of renormalized q ddagger M_n final}
& & \int q^{\ddagger,\Theta_0\left(r_0\right)}(\theta) M_n(\theta ; x) \mathrm{d} \theta \nonumber\\
&\geq& -\frac{p}{2}-\frac{1}{2}\sum\limits_{i=1}^p\log\frac{v_{ni}}{v_{0i}}+M_n(\theta^*;x)+\log p\left(\mathbf{m}_0^*\right)-\log p\left(\mathbf{m}_n^*\right)\nonumber\\
& & +\frac{1}{2}\Delta_{n,\theta_0}^{\top}V_{\theta_0}\Delta_{n,\theta_0} - C_4 \left(R_n + \frac{1}{n} \mathrm{tr}\left(\mathrm{diag}(\mathbf{v}_n^*)\right) \right) 
\end{eqnarray}
holds for a certain $C_4>0$, by the definitions of $\mathbf{m}_0^*$, $\mathbf{m}_n^*$, $\mathbf{v}_0^*$, and $\mathbf{v}_n^*$. 
Moreover, under Condition \ref{condition: prior}, we have
\begin{equation*}
\label{log p theta* - log p m}
\left|\log p(\mathbf{m}_n^*)-\log p(\mathbf{m}_0^*)-\left(\mathbf{m}_n^*-\mathbf{m}_0^*\right)^{\top}\nabla\log p(\mathbf{m}_0^*)\right|\leq\frac{\tau_5}{2}\left\|\mathbf{m}_n^*-\mathbf{m}_0^*\right\|^2.
\end{equation*}
Then, according to the proof of \eqref{rate of v and m}, we can control the difference of the logarithmic terms as follows:
\begin{equation}
\label{control of logm0-logmn}
\left|\log p\left(\mathbf{m}_0^*\right)-\log p\left(\mathbf{m}_n^*\right)\right| \lesssim R_n + \frac{1}{n} \mathrm{tr}\left(\mathrm{diag}(\mathbf{v}_n^*)\right).
\end{equation}
Based on \eqref{vni/voi-1} and the fact that $R_n$ dominates $n^{-1}\left\|\mathbf{v}_0^*\right\| \sqrt{\sum\limits_{i=1}^p\left(\frac{v_{ni}}{v_{0i}}-1\right)^2}$, we have
\begin{equation*}
\left(\sum\limits_{i=1}^p\log\frac{v_{ni}}{v_{0i}}\right)^2 \lesssim p    \sum\limits_{i=1}^p\left(\frac{v_{ni}}{v_{0i}}-1\right)^2 \lesssim p R_n.
\end{equation*}
It implies that
\begin{eqnarray}
\label{control of log vn/v0}
 \left|\sum\limits_{i=1}^p\log\frac{v_{ni}}{v_{0i}}\right| \lesssim \sqrt{p}\sqrt{R_n} \lesssim p+R_n.   
\end{eqnarray}
Therefore, accounting for the absorption of $\mathrm{tr}\left(\mathrm{diag}(\mathbf{v}_n^*)\right)$, there exists a constant $C_5>0$ such that
\begin{equation}
\label{Final lower bound of q ddagger Mn}
\int q^{\ddagger,\Theta_0\left(r_0\right)}(\theta) M_n(\theta ; x) \mathrm{d} \theta \geq M_n(\theta^*;x)+\frac{1}{2}\Delta_{n,\theta_0}^{\top}V_{\theta_0}\Delta_{n,\theta_0} - C_5\left(p+R_n\right)
\end{equation}
holds with probability at least $1-5e^{-y}-e^{-cr_0^2}$, based on \eqref{lower bound of renormalized q ddagger M_n final}--\eqref{control of log vn/v0}.
In summary, we have established both upper and lower bounds for $\int q^{\ddagger,\Theta_0\left(r_0\right)}(\theta) M_n(\theta ; x) \mathrm{d} \theta$.

{\textit{Step 2.}} Derive a non-asymptotic upper bound for $\int_{B^c\left(\theta^*, \eta\right)}q^{\ddagger,\Theta_0\left(r_0\right)}(\theta)\mathrm{d}\theta$.

Comparing \eqref{upper bound of renormalized q ddagger M_n final} and \eqref{Final lower bound of q ddagger Mn}, we have
\begin{equation}
\label{result for renormalized q ddagger}
\int_{B^c\left(\theta^*, \eta\right)} q^{\ddagger,\Theta_0\left(r_0\right)}(\theta)\mathrm{d} \theta \lesssim \frac{p+R_n}{\lambda_{\min}\left(V_{\theta_0}\right)n\eta^2+\Delta_{n,\theta_0}^{\top}V_{\theta_0}\Delta_{n,\theta_0}}
\end{equation}
holds with probability at least $1-5e^{-y}-e^{-cr_0^2}$.

{\textit{Step 3.}} Extend the result from $q^{\ddagger,\Theta_0\left(r_0\right)}(\theta)$ to $q^{\ddagger}(\theta)$. That is, we show that
\begin{equation}
\label{upper bound of the edge mass of the bridge}
\int_{B^c\left(\theta^*, \eta\right)} q^{\ddagger}(\theta) \mathrm{d} \theta\lesssim\frac{p+R_n}{\lambda_{\min}\left(V_{\theta_0}\right)n\eta^2+\Delta_{n,\theta_0}^{\top}V_{\theta_0}\Delta_{n,\theta_0}}
\end{equation}
holds with probability at least $1-5e^{-y}-e^{-cr_0^2}$.

Due to \eqref{result for renormalized q ddagger}, there exists a constant $C_6>0$ such that
\begin{equation*}
\int_{B\left(\theta^*, \eta\right)} q^{\ddagger,\Theta_0\left(r_0\right)}(\theta)\mathrm{d} \theta\geq 1-C_6\frac{p+ R_n}{\lambda_{\min}\left(V_{\theta_0}\right)n\eta^2+\Delta_{n,\theta_0}^{\top}V_{\theta_0}\Delta_{n,\theta_0}}
\end{equation*}
holds with probability at least $1-5e^{-y}-e^{-cr_0^2}$.
As $r_0$ grows with $n, p$, and  $\eta$ tends to zero, it suffices to consider the case where $B\left(\theta^*, \eta\right)\subset\Theta_0\left(r_0\right)$. Therefore, the above inequality implies that
\begin{equation}
\label{ratio of B/Theta_0 r_0}
\frac{\int_{B\left(\theta^*, \eta\right)} q^{\ddagger}(\theta) \mathrm{d} \theta}{\int_{\Theta_0\left(r_0\right)} q^{\ddagger}(\theta) \mathrm{d} \theta}
\geq 1-C_6\frac{p+R_n}{\lambda_{\min}\left(V_{\theta_0}\right)n\eta^2+\Delta_{n,\theta_0}^{\top}V_{\theta_0}\Delta_{n,\theta_0}}.
\end{equation}
Recall that $\int_{\Theta_0^c\left(r_0\right)} q^{\ddagger}(\theta) \mathrm{d} \theta$ is exponentially small by Lemma \ref{lemma: The Tail Moment Control for the p-Dimensional Normal Distribution} and Theorem \ref{theorem: the convergence rates of the VBBE and the covariance matrix of the VB bridge}. That is, there exists a constant $C_7>0$ such that
\begin{equation*}
\int_{\Theta_0^c\left(r_0\right)} q^{\ddagger}(\theta) \mathrm{d} \theta\lesssim\exp\left(-\frac{C_7r_0^2}{\rho\left(V_{\theta_0} V_{\ddagger}\right)}\right).
\end{equation*}
Consequently, there exists a constant $C_8>0$ such that
\begin{equation}
\label{lower bound of q ddagger Theta_0 r_0}
\int_{\Theta_0\left(r_0\right)} q^{\ddagger}(\theta) \mathrm{d} \theta\geq 1-C_8\exp\left(-\frac{C_7r_0^2}{\rho\left(V_{\theta_0} V_{\ddagger}\right)}\right).
\end{equation}
Based on \eqref{ratio of B/Theta_0 r_0} and \eqref{lower bound of q ddagger Theta_0 r_0}, we have
\begin{equation*}
\int_{B\left(\theta^*, \eta\right)} q^{\ddagger}(\theta) \mathrm{d} \theta
\geq \left(1-C_6\frac{p+R_n}{\lambda_{\min}\left(V_{\theta_0}\right)n\eta^2+\Delta_{n,\theta_0}^{\top}V_{\theta_0}\Delta_{n,\theta_0}}\right)\left[1-C_8\exp\left(-\frac{C_7r_0^2}{\rho\left(V_{\theta_0} V_{\ddagger}\right)}\right)\right]
\end{equation*}
holds with probability at least $1-5e^{-y}-e^{-cr_0^2}$, thereby arriving at conclusion \eqref{upper bound of the edge mass of the bridge}.

\section{Proof of Theorem \ref{theorem: the asymptotic normality of the VB bridge}}
\label{Proof of Theorem 3}

The Kullback-Leibler (KL) divergence between the variational distributions $q(\theta ;\mathbf{m},\mathbf{v})$ and $q(\theta ;\mathbf{m+\Delta m},\mathbf{v+\Delta v})$ is first computed as follows:
\begin{eqnarray*}
& & \mathrm{KL}(q(\theta ;\mathbf{m},\mathbf{v}) \| q(\theta ;\mathbf{m+\Delta m},\mathbf{v+\Delta v})\\
&=& \sum\limits_{i=1}^p\int N\left(\theta_i;m_i,\frac{v_i}{n}\right)\log N\left(\theta_i;m_i,\frac{v_i}{n}\right)\mathrm{d}\theta_i\\
& & -\sum\limits_{i=1}^p\int N\left(\theta_i;m_i,\frac{v_i}{n}\right)\log N\left(\theta_i;m_i+\Delta m_i,\frac{v_i+\Delta v_i}{n}\right)\mathrm{d}\theta_i\\
&=& -\frac{p}{2}-\frac{p}{2}\log\left(\frac{2\pi}{n}\right)-\frac{1}{2}\sum\limits_{i=1}^p\log v_i\\
& & +\frac{1}{2}\sum\limits_{i=1}^p\frac{v_i+n\left(\Delta m_i\right)^2}{v_i+\Delta v_i}+\frac{p}{2}\log\left(\frac{2\pi}{n}\right)+\frac{1}{2}\sum\limits_{i=1}^p\log\left(v_i+\Delta v_i\right)\\
&=& \frac{1}{2}\sum\limits_{i=1}^p\left(\frac{n\left(\Delta m_i\right)^2-\Delta v_i}{v_i+\Delta v_i}+\log\left(1+\frac{\Delta v_i}{v_i}\right)\right)\\ 
&=& \frac{n}{2}\sum\limits_{i=1}^p\frac{\left(\Delta m_i\right)^2}{v_i+\Delta v_i}+\frac{1}{2}\sum\limits_{i=1}^p\left(-\frac{\Delta v_i}{v_i+\Delta v_i}+\log\left(1+\frac{\Delta v_i}{v_i}\right)\right).
\end{eqnarray*}
Substituting $\mathbf{m},\mathbf{v},\mathbf{\Delta m},\mathbf{\Delta v}$ with $\mathbf{m}_n^*,\mathbf{v}_n^*,\mathbf{m}_0^*-\mathbf{m}_n^*,\mathbf{v}_0^*-\mathbf{v}_n^*$, we find
\begin{equation}
\label{KLn-KL0}
\mathrm{KL}(q(\theta ;\mathbf{m}_n^*,\mathbf{v}_n^*) \| q(\theta ;\mathbf{m}_0^*,\mathbf{v}_0^*))=\frac{n}{2}\sum\limits_{i=1}^p\frac{\left(m_{0i}-m_{ni}\right)^2}{v_{0i}}+\frac{1}{2}\sum\limits_{i=1}^p\left(\frac{v_{ni}}{v_{0i}}-\log\frac{v_{ni}}{v_{0i}}-1\right).
\end{equation}
By definition, we know that
\begin{eqnarray*}
 q(\theta ;\mathbf{m}_n^*,\mathbf{v}_n^*) &=& \underset{q \in \mathcal{Q}^p}{\arg \min } \operatorname{KL}\left(q(\theta) \| \pi^*(\theta \mid x)\right),  \\
q(\theta ;\mathbf{m}_0^*,\mathbf{v}_0^*) &=& \underset{q \in \mathcal{Q}^p}{\arg \min } \operatorname{KL}\left(q(\theta) \| N\left(\theta ; \theta^*+\frac{\Delta_{n, \theta_0}}{\sqrt{n}},\frac{V_{\theta_0}^{-1}}{n}\right)\right). 
\end{eqnarray*}
Then, taking the linear transformation $\tilde{\theta}=\sqrt{n}\left(\theta-\theta^*\right)$, we have
\begin{eqnarray*}
& &\mathrm{KL}\left(q(\theta ;\mathbf{m}_n^*,\mathbf{v}_n^*) \|q(\theta ;\mathbf{m}_0^*,\mathbf{v}_0^*) \right)\\
&=& \mathrm{KL}\left(\underset{q \in \mathcal{Q}^p}{\arg \min } \operatorname{KL}\left(q(\theta) \| \pi^*(\theta \mid x)\right)\|\underset{q \in \mathcal{Q}^p}{\arg \min } \operatorname{KL}\left(q(\theta) \| N\left(\theta ; \theta^*+\frac{\Delta_{n, \theta_0}}{\sqrt{n}},\frac{V_{\theta_0}^{-1}}{n}\right)\right) \right)\\
&=& \mathrm{KL}\left( \underset{q \in \mathcal{Q}^p}{\arg \min } \operatorname{KL}\left(q(\tilde{\theta}) \| \pi_{\tilde{\theta}}^*(\tilde{\theta}\mid x)\right)\|\underset{q \in \mathcal{Q}^p}{\arg \min } \operatorname{KL}\left(q(\tilde{\theta}) \| N\left(\tilde{\theta}; \Delta_{n, \theta_0}, V_{\theta_0}^{-1}\right)\right)\right).
\end{eqnarray*}
Combined with \eqref{final rate of ddagger} and \eqref{KLn-KL0}, we derive the following bound for the total variation distance between the two KL divergence minimizers:
\begin{eqnarray*}
& & \left\| \underset{q \in \mathcal{Q}^p}{\arg \min } \operatorname{KL}\left(q(\tilde{\theta}) \| \pi_{\tilde{\theta}}^*(\tilde{\theta}\mid x)\right)-\underset{q \in \mathcal{Q}^p}{\arg \min } \operatorname{KL}\left(q(\tilde{\theta}) \| N\left(\tilde{\theta}; \Delta_{n, \theta_0}, V_{\theta_0}^{-1}\right)\right)\right\|^2_{\mathrm{TV}}\\
&\leq& \frac{1}{2}\mathrm{KL}\left( \underset{q \in \mathcal{Q}^p}{\arg \min } \operatorname{KL}\left(q(\tilde{\theta}) \| \pi_{\tilde{\theta}}^*(\tilde{\theta}\mid x)\right)\|\underset{q \in \mathcal{Q}^p}{\arg \min } \operatorname{KL}\left(q(\tilde{\theta}) \| N\left(\tilde{\theta}; \Delta_{n, \theta_0}, V_{\theta_0}^{-1}\right)\right)\right) \\
&\lesssim& \frac{\max_{i}(v_{\theta_0,ii})R_n}{\lambda_{\min}\left(V_{\theta_0}\right)},
\end{eqnarray*}
where thee first inequality holds by the Pinsker inequality.

\section{Proof of Theorem \ref{theorem: the convergence rates of the distribution parameters of the VB posterior}}
\label{Proof of Theorem 4}

We follow an approach similar to that used in the proof of Theorem \ref{theorem: the convergence rates of the VBBE and the covariance matrix of the VB bridge} to derive the convergence rates of $\hat{\theta}_n^*$ and $\hat{V}^*_n$ around $\theta^*$ and $\mathrm{diag}^{-1}\left(V_{\theta_0}\right)$. 

According to the definition of $q^*(\theta)$, we first know that $\operatorname{ELBO}_p(q_0(\theta)) \leq \operatorname{ELBO}_p(q^*(\theta))$. Further, by Lemma \ref{lemma: the gap between two functionals}, we know that
\begin{equation*}
  \operatorname{ELBO}_p(q^*(\theta)) + \mathrm{KL}\left(q^*(\theta) \| \pi^*(\theta \mid x)\right) - \log \int p(\theta) \exp \left\{M_n(\theta ; x)\right\} \mathrm{d} \theta \leq 0,  
\end{equation*}
and
\begin{eqnarray*}
   0 &\leq& \log \int p(\theta) \exp \left\{M_n(\theta ; x)\right\} \mathrm{d} \theta-\mathrm{KL}\left(q_0(\theta) \| \pi^*(\theta \mid x)\right)-\operatorname{ELBO}_p(q_0(\theta))\\
   &\lesssim& \frac{p}{\min v_{\theta_0,ii}}\left[\left(\frac{p}{n\min v_{\theta_0,ii}}\right)^{\frac{s_1}{2}}+\left(\frac{p}{n\min v_{\theta_0,ii}}\right)^{\frac{s_2}{2}}\right]
\end{eqnarray*}
hold with probability at least $1-e^{-cp}$. Then, we can conclude that there exists a constant $C_1>0$ such that
\begin{equation}
\label{upper bound of KL q* - q0}
\mathrm{KL}\left(q^*(\theta) \| \pi^*(\theta \mid x)\right)-\mathrm{KL}\left(q_0(\theta) \| \pi^*(\theta \mid x)\right)\leq \frac{C_1p}{\min v_{\theta_0,ii}}\left[\left(\frac{p}{n\min v_{\theta_0,ii}}\right)^{\frac{s_1}{2}}+\left(\frac{p}{n\min v_{\theta_0,ii}}\right)^{\frac{s_2}{2}}\right]
\end{equation}
with probability at least $1-e^{-cp}$. Recall that $F_n(\mathbf{m},\mathbf{v}) =\operatorname{KL}\left(q(\theta ;\mathbf{m},\mathbf{v}) \| \pi^*(\theta \mid x)\right)$ and $\hat{V}^*_n = \mathrm{diag}\left(\left(\hat{v}^*_{n,11},\ldots,\hat{v}^*_{n,pp}\right)^{\top}\right)$. By \eqref{def of F0} and \eqref{Fn=F0+error}, we have
\begin{eqnarray}
\label{Fn* - Fn0}
& & \mathrm{KL}\left(q^*(\theta) \| \pi^*(\theta \mid x)\right)-\mathrm{KL}\left(q_0(\theta) \| \pi^*(\theta \mid x)\right)\nonumber\\
&=& \frac{1}{2}\sum\limits_{i=1}^p\left(\hat{v}^*_{n,ii}v_{\theta_0,ii}-\log \hat{v}^*_{n,ii}v_{\theta_0,ii}-1\right) +\frac{n}{2}\left(\hat{\theta}_n^*-\theta^*-\frac{\Delta_{n, \theta_0}}{\sqrt{n}}\right)^{\top}V_{\theta_0}\left(\hat{\theta}_n^*-\theta^*-\frac{\Delta_{n, \theta_0}}{\sqrt{n}}\right)\nonumber\\
& & +\log p\left(\theta^*+\frac{\Delta_{n, \theta_0}}{\sqrt{n}}\right)-\log p\left(\hat{\theta}_n^*\right) + O\left(R_n + \frac{1}{n}\mathrm{tr}\left(\hat{V}^*_n\right)\right).
\end{eqnarray}
As in Step 1 of Theorem \ref{theorem: the convergence rates of the VBBE and the covariance matrix of the VB bridge}, where $W_1$ is shown to be dominated by $W_2$, and in Step 2, where $\mathrm{tr}\left(\hat{V}^*_n\right)$ is absorbed into $R_n$, we have
\begin{eqnarray*}
 & & \sum\limits_{i=1}^p\left(\hat{v}^*_{n,ii}v_{\theta_0,ii}-\log \hat{v}^*_{n,ii}v_{\theta_0,ii}-1\right)+n\left(\hat{\theta}_n^*-\theta^*-\frac{\Delta_{n, \theta_0}}{\sqrt{n}}\right)^{\top}V_{\theta_0}\left(\hat{\theta}_n^*-\theta^*-\frac{\Delta_{n, \theta_0}}{\sqrt{n}}\right)\\
&\lesssim& R_n+\frac{p}{\min v_{\theta_0,ii}}\left[\left(\frac{p}{n\min v_{\theta_0,ii}}\right)^{\frac{s_1}{2}}+\left(\frac{p}{n\min v_{\theta_0,ii}}\right)^{\frac{s_2}{2}}\right], 
\end{eqnarray*}
with probability at least $1-5e^{-y}-e^{-cp}-e^{-cr_0^2}$, by \eqref{upper bound of KL q* - q0} and \eqref{Fn* - Fn0}.

\section{Proof of Corollary \ref{corollary of theorem 4}}
\label{Proof of Corollary 1}

From Theorem \ref{theorem: the convergence rates of the distribution parameters of the VB posterior}, it follows that
\begin{equation*}
\left\|\hat{\theta}_n^*-\theta^*\right\|^2 \lesssim \frac{1}{n\lambda_{\min}\left(V_{\theta_0}\right)} \left\{R_n+\frac{p}{\min v_{\theta_0,ii}}\left[\left(\frac{p}{n\min v_{\theta_0,ii}}\right)^{\frac{s_1}{2}}+\left(\frac{p}{n\min v_{\theta_0,ii}}\right)^{\frac{s_2}{2}}\right]\right\} +\frac{\left\|\Delta_{n, \theta_0}\right\|^2}{n}
\end{equation*}
holds with probability at least $1-5e^{-y}-e^{-cp}-e^{-cr_0^2}$. We now sequentially analyze the order of each term in this upper bound.

Recall that
\begin{equation*}
R_n=\Delta(r_0,y)+\frac{1}{\sqrt{n}} \left\{\left\|\Delta_{n,\theta_0}\right\|^2+\left\|\nabla\log p\left(\theta^*\right)\right\|^2\right\} + \frac{1}{n} \left\{\mathrm{tr}\left(V_{\theta_0}^{-1}\right)+\mathrm{tr}\left(\mathrm{diag}^{-1}\left(V_{\theta_0}\right)\right)\right\}.
\end{equation*}
As in Remark \ref{Remark 1}, we can chose $y\asymp p$. Then, combined with (i), we have
\begin{equation*}
\Delta(r_0,y)=O_p\left(\sqrt{\frac{p^3}{n}}\right).
\end{equation*}
For the second term in $R_n$, based on Conditions \ref{condition: local, exponential moment}, \ref{condition: identifiability}, (ii) and (iv), we can show that
\begin{equation*}
\frac{1}{\sqrt{n}} \left\{\left\|\Delta_{n,\theta_0}\right\|^2+\left\|\nabla\log p\left(\theta^*\right)\right\|^2\right\} = O_p\left(\frac{p+p^{\frac{3}{2}}}{\sqrt{n}}\right) = O_p\left(\sqrt{\frac{p^3}{n}}\right).
\end{equation*}
Under Conditions (ii) and (iii), it is straightforward to verify that
\begin{equation*}
\frac{1}{n} \left\{\mathrm{tr}\left(V_{\theta_0}^{-1}\right)+\mathrm{tr}\left(\mathrm{diag}^{-1}\left(V_{\theta_0}\right)\right)\right\} = O\left(\frac{p}{n}\right).
\end{equation*}
Therefore, in summary, $R_n=O_p\left(\sqrt{p^3/{n}}\right)$. 

Applying Condition (iii), one can easily verify that
\begin{equation*}
\frac{p}{\min v_{\theta_0,ii}}\left[\left(\frac{p}{n\min v_{\theta_0,ii}}\right)^{\frac{s_1}{2}}+\left(\frac{p}{n\min v_{\theta_0,ii}}\right)^{\frac{s_2}{2}}\right]=O\left(p\left[\left(\frac{p}{n}\right)^{\frac{s_1}{2}}+\left(\frac{p}{n}\right)^{\frac{s_2}{2}}\right]\right).
\end{equation*}
Based on Conditions \ref{condition: local, exponential moment}, \ref{condition: identifiability} and (ii), we can obtain that
\begin{equation*}
\frac{\left\|\Delta_{n, \theta_0}\right\|^2}{n}=O_p\left(\frac{p}{n}\right).
\end{equation*}
Consequently, under Condition (ii), we have
\begin{equation*}
\|\hat{\theta}_n^*-\theta^*\|^2 = O_p\left(\left(\frac{p}{n}\right)^{\frac{3}{2}}+\frac{p}{n}\left[\left(\frac{p}{n}\right)^{\frac{s_1}{2}}+\left(\frac{p}{n}\right)^{\frac{s_2}{2}}\right]+\frac{p}{n}\right),
\end{equation*}
which implies that
\begin{equation*}
\|\hat{\theta}_n^*-\theta^*\| = O_p\left(\sqrt{\frac{p}{n}}\right).
\end{equation*}

\section{Proof of Corollary \ref{Corollary 2}}
\label{Proof of Corollary 2}

Let $r_n = (\hat{\theta}_n^*-\theta^*) - \Delta_{n,\theta_0}/\sqrt{n}$. Then, from Theorem \ref{theorem: the convergence rates of the distribution parameters of the VB posterior}, we know that
\begin{equation*}
\left\|r_n\right\|^2 \lesssim \frac{1}{n\lambda_{\min}\left(V_{\theta_0}\right)} \left\{R_n+\frac{p}{\min v_{\theta_0,ii}}\left[\left(\frac{p}{n\min v_{\theta_0,ii}}\right)^{\frac{s_1}{2}}+\left(\frac{p}{n\min v_{\theta_0,ii}}\right)^{\frac{s_2}{2}}\right]\right\}
\end{equation*}
holds with probability at least $1-5e^{-y}-e^{-cp}-e^{-cr_0^2}$. According to the proof of Corollary \ref{corollary of theorem 4}, we can show that $\left\|r_n\right\|=o_P\left(n^{-1 / 2}\right)$ when $p^3=o(n)$ and $\min\left(s_1,s_2\right) \geq 1$. Thus, 
\begin{equation*}
\sqrt{n} \alpha^{\top}\left(\hat{\theta}_n^*-\theta^*\right)/\sigma_\alpha = \frac{1}{\sigma_\alpha} \left\{\alpha^{\top} V_{\theta_0}^{-1}\frac{1}{\sqrt{n}} \sum_{i=1}^n\nabla m(\theta^*;x_i)+ \sqrt{n} \alpha^{\top}r_n\right\}
\overset{\mathcal{L}}{\longrightarrow} N(0,1),
\end{equation*}
for any $\alpha \in S_p=\left\{\alpha \in \mathbb{R}^p:\|\alpha\|=1\right\}$, where $\sigma_\alpha^2=\alpha^{\top} V_{\theta_0}^{-1} \operatorname{Var}\left(\nabla m\left(\theta^* ; x\right)\right) V_{\theta_0}^{-1} \alpha$. 

\section{Proof of Theorem \ref{theorem: the non-asymptotic variational Bernstein--von Mises theorem}}
\label{Proof of Theorem 5}

To characterize the concentration of the VB posterior $q^*(\theta)$, we employ techniques similar to those used in the proof of Theorem \ref{theorem: the consistency of the VB bridge}. Specifically, the proof will be carried out in three main steps.

{\textit{Step 1.}} Renormalize $q^*(\theta)$ over $\Theta_0\left(r_0\right)$ as
\begin{equation*}
q^{*,\Theta_0(r_0)}(\theta) = \frac{q^*(\theta) \cdot I_{\Theta_0(r_0)}(\theta)}{\int q^*(\theta) \cdot I_{\Theta_0(r_0)}(\theta) \mathrm{d} \theta},
\end{equation*}
and establish both upper and lower bounds for $\int q^{*,\Theta_0\left(r_0\right)}(\theta) M_n(\theta; x) \mathrm{d} \theta$.

Similar to the derivation of the upper bound for $\int q^{\ddagger,\Theta_0\left(r_0\right)}(\theta) M_n(\theta ; x) \mathrm{d} \theta$ in \eqref{upper bound of renormalized q ddagger M_n final}, we derive the following upper bound for $\int q^{*,\Theta_0\left(r_0\right)}(\theta) M_n(\theta; x) \mathrm{d} \theta$:
\begin{eqnarray}
\label{upper bound of q*Mn}
& & \int q^{*,\Theta_0\left(r_0\right)}(\theta) M_n(\theta ; x) \mathrm{d} \theta \nonumber\\
&\leq& M_n(\theta^*;x)+\Delta(r_0,y)+\frac{1}{2}\Delta_{n, \theta_0}^{\top}V_{\theta_0}\Delta_{n, \theta_0}\left(1-\int_{B^c\left(\theta^*, \eta\right)} q^{*,\Theta_0\left(r_0\right)}(\theta)\mathrm{d} \theta\right)\nonumber\\
& & -C_1\lambda_{\min}\left(V_{\theta_0}\right)n\eta^2\int_{B^c\left(\theta^*, \eta\right)} q^{*,\Theta_0\left(r_0\right)}(\theta)\mathrm{d} \theta \nonumber\\
&=& M_n(\theta^*;x)+\Delta(r_0,y)+\frac{1}{2} \Delta_{n, \theta_0}^{\top}V_{\theta_0}\Delta_{n, \theta_0} -\left(\frac{1}{2} \Delta_{n, \theta_0}^{\top}V_{\theta_0}\Delta_{n, \theta_0}+C_1\lambda_{\min}\left(V_{\theta_0}\right)n\eta^2\right) \nonumber\\
& &\times\int_{B^c\left(\theta^*, \eta\right)} q^{*,\Theta_0\left(r_0\right)}(\theta)\mathrm{d} \theta
\end{eqnarray}
holds with probability at least $ 1 - e^{-y}$, where $C_1\in\left(0,1/2\right)$ is a constant.

We then derive a lower bound for $\int q^{*,\Theta_0\left(r_0\right)}(\theta) M_n(\theta; x) \mathrm{d} \theta$. First, according to the definition of KL divergence, $\int q^{*,\Theta_0\left(r_0\right)}(\theta) M_n(\theta; x) \mathrm{d} \theta$ can be decomposed as
\begin{eqnarray}
\label{decomposition of q*Mn}
 & & \int q^{*,\Theta_0\left(r_0\right)}(\theta) M_n(\theta ; x) \mathrm{d} \theta \nonumber\\
 &=& \int q^{*,\Theta_0\left(r_0\right)}(\theta)\log q^{*,\Theta_0\left(r_0\right)}(\theta)\mathrm{d}\theta +\log \int p(\theta) \exp \left\{M_n(\theta ; x)\right\} \mathrm{d} \theta -\int q^{*,\Theta_0\left(r_0\right)}(\theta)\log p(\theta) \mathrm{d} \theta \nonumber\\
 & & -\operatorname{KL}\left(q^{*,\Theta_0\left(r_0\right)}(\theta)\| \pi^*(\theta \mid x)\right) \nonumber\\
 &=& \int q^{*,\Theta_0\left(r_0\right)}(\theta)\log q^{*,\Theta_0\left(r_0\right)}(\theta)\mathrm{d}\theta +\log \int p(\theta) \exp \left\{M_n(\theta ; x)\right\} \mathrm{d} \theta -\int q^{*,\Theta_0\left(r_0\right)}(\theta)\log p(\theta) \mathrm{d} \theta \nonumber\\
 & & -\operatorname{KL}\left(q^*(\theta)\| \pi^*(\theta \mid x)\right) + \left\{\operatorname{KL}\left(q^*(\theta)\| \pi^*(\theta \mid x)\right)-\operatorname{KL}\left(q^{*,\Theta_0\left(r_0\right)}(\theta)\| \pi^*(\theta \mid x)\right)\right\}\nonumber\\
 &\triangleq& S_1 + S_2 - S_3 - S_4 + S_5.
\end{eqnarray}

We analyze these terms one by one. Using a similar analytical technique as that used to establish conclusion \eqref{Identity for the first term}, we can prove that 
\begin{equation}
\label{q*Theta0r0 log q*Theta0r0}
    S_1 = -\frac{p}{2}-\frac{p}{2}\log(2\pi)+\frac{p}{2}\log n-\frac{1}{2}\sum\limits_{i=1}^p\log \hat{v}^*_{n,ii}+S_{1n},
\end{equation}
where $S_{1n}$ is exponentially small. For the term $S_2$, note that it is identical to $T_2$,  and thus the conclusion in \eqref{Lower bound for the second term} follows. For the term $S_3$, using a similar analytical technique as that used to establish conclusion \eqref{Upper bound for the third term}, we can show that 
\begin{equation}
\label{q*Theta0r0 log p}
S_3 \leq \log p\left(\hat{\theta}_n^*\right) + \frac{\tau_5}{2n}\mathrm{tr}\left(\hat{V}^*_n\right) + S_{3n},
\end{equation}
where $S_{3n}$ is exponentially small. Substituting $\left(\mathbf{m},\mathbf{v}\right)=\left(\hat{\theta}_n^*,\left(\hat{v}^*_{n,11},\cdots,\hat{v}^*_{n,pp}\right)^{\top}\right)$ into \eqref{Fn=F0+error}, we obtain that
\begin{eqnarray}
\label{KL q*}
S_4 &=& -\frac{p}{2}-\frac{1}{2}\sum\limits_{i=1}^p\log \hat{v}^*_{n,ii}-\frac{1}{2}\log\det\left(V_{\theta_0}\right)+\frac{1}{2}\sum\limits_{i=1}^p\hat{v}^*_{n,ii}v_{\theta_0,ii} \nonumber\\
& & +\frac{n}{2}\left(\hat{\theta}_n^*-\theta^*-\frac{\Delta_{n, \theta_0}}{\sqrt{n}}\right)^{\top}V_{\theta_0}\left(\hat{\theta}_n^*-\theta^*-\frac{\Delta_{n, \theta_0}}{\sqrt{n}}\right) \nonumber\\
& & +\log p(\theta^*)-\log p(\hat{\theta}_n^*)\nonumber\\
& & +O\left(\Delta(r_0,y)+\frac{\left\|\Delta_{n,\theta_0}\right\|^2+\left\|\nabla\log p\left(\theta^*\right)\right\|^2}{\sqrt{n}}+\frac{\mathrm{tr}\left(V_{\theta_0}^{-1}\right)+\mathrm{tr}\left(\hat{V}^*_n\right)}{n}\right),
\end{eqnarray}
with probability at least $1-5e^{-y}-e^{-cr_0^2}$. Similar to $T_5$ defined in \eqref{Identity for the normalized q ddagger}, the discrepancy $S_5$ is exponentially negligible. Combined with the results \eqref{Lower bound for the second term} and \eqref{q*Theta0r0 log q*Theta0r0}--\eqref{KL q*}, we establish the following lower bound for $\int q^{*,\Theta_0\left(r_0\right)}(\theta) M_n(\theta ; x) \mathrm{d} \theta$: with probability at least $1-5e^{-y}-e^{-cr_0^2}$,
\begin{eqnarray}
\label{lower bound of q*Mn}
& & \int q^{*,\Theta_0\left(r_0\right)}(\theta) M_n(\theta ; x) \mathrm{d} \theta \nonumber\\
&=& M_n(\theta^*;x)+\frac{1}{2}\Delta_{n,\theta_0}^{\top}V_{\theta_0}\Delta_{n,\theta_0} - \frac{1}{2}\sum\limits_{i=1}^p\hat{v}^*_{n,ii}v_{\theta_0,ii}\nonumber\\
& & -\frac{n}{2}\left(\hat{\theta}_n^*-\theta^*-\frac{\Delta_{n, \theta_0}}{\sqrt{n}}\right)^{\top}V_{\theta_0}\left(\hat{\theta}_n^*-\theta^*-\frac{\Delta_{n, \theta_0}}{\sqrt{n}}\right)\nonumber\\
& & -C_2\left(\Delta(r_0,y)+\frac{\left\|\Delta_{n,\theta_0}\right\|^2+\left\|\nabla\log p\left(\theta^*\right)\right\|^2}{\sqrt{n}}+\frac{\mathrm{tr}\left(V_{\theta_0}^{-1}\right)+\mathrm{tr}\left(\hat{V}^*_n\right)}{n}\right).
\end{eqnarray}
Rewrite $-\sum\limits_{i=1}^p\hat{v}^*_{n,ii}v_{\theta_0,ii} / 2$ from \eqref{lower bound of q*Mn} as:
\begin{equation}
\label{Keq1-xpr}
-\frac{1}{2}\sum\limits_{i=1}^p\hat{v}^*_{n,ii}v_{\theta_0,ii}=-\frac{p}{2}-\frac{1}{2}\sum\limits_{i=1}^p\left(\hat{v}^*_{n,ii}v_{\theta_0,ii}-1\right).
\end{equation}
On the one hand, we observe that
\begin{equation}
\label{cauchy for v* v theta0}
\left[\sum\limits_{i=1}^p\left(\hat{v}^*_{n,ii}v_{\theta_0,ii}-1\right)\right]^2\leq p\sum\limits_{i=1}^p\left(\hat{v}^*_{n,ii}v_{\theta_0,ii}-1\right)^2
\end{equation}
by Cauchy-Schwarz inequality. On the other hand, Theorem \ref{theorem: the convergence rates of the distribution parameters of the VB posterior} ensures that
\begin{equation*}
\sum\limits_{i=1}^p\left(\hat{v}^*_{n,ii}v_{\theta_0,ii}-\log \hat{v}^*_{n,ii}v_{\theta_0,ii}-1\right)
\lesssim R_n +\frac{p}{\min v_{\theta_0,ii}}\left[\left(\frac{p}{n\min v_{\theta_0,ii}}\right)^{\frac{s_1}{2}}+\left(\frac{p}{n\min v_{\theta_0,ii}}\right)^{\frac{s_2}{2}}\right],
\end{equation*}
with probability at least $1-5e^{-y}-e^{-cp}-e^{-cr_0^2}$. Note that $x-\log x-1  \asymp (x-1)^2$ when $x \rightarrow 1$, we then obtain that 
\begin{equation}
\label{v*vtheta0-logv*vtheta0-1}
\sum\limits_{i=1}^p\left(\hat{v}^*_{n,ii}v_{\theta_0,ii}-1\right)^2 
\lesssim R_n + \frac{p}{\min v_{\theta_0,ii}}\left[\left(\frac{p}{n\min v_{\theta_0,ii}}\right)^{\frac{s_1}{2}}+\left(\frac{p}{n\min v_{\theta_0,ii}}\right)^{\frac{s_2}{2}}\right].
\end{equation}
By Theorem \ref{theorem: the convergence rates of the distribution parameters of the VB posterior}, combined with \eqref{Keq1-xpr} - \eqref{v*vtheta0-logv*vtheta0-1}, we can show that
\begin{eqnarray}
 \label{v*vtheta0 VBE}   
& & \sum\limits_{i=1}^p\hat{v}^*_{n,ii}v_{\theta_0,ii}+n\left(\hat{\theta}_n^*-\theta^*-\frac{\Delta_{n, \theta_0}}{\sqrt{n}}\right)^{\top}V_{\theta_0}\left(\hat{\theta}_n^*-\theta^*-\frac{\Delta_{n, \theta_0}}{\sqrt{n}}\right) \nonumber\\
&\lesssim& p+R_n +\frac{p}{\min v_{\theta_0,ii}}\left[\left(\frac{p}{n\min v_{\theta_0,ii}}\right)^{\frac{s_1}{2}}+\left(\frac{p}{n\min v_{\theta_0,ii}}\right)^{\frac{s_2}{2}}\right].
\end{eqnarray}
Further, using a similar contradiction-based argument as in Step 2 of Theorem \ref{theorem: the convergence rates of the VBBE and the covariance matrix of the VB bridge}, we can show that $R_n$ dominates $ n^{-1} \left\|\mathbf{v}_0^*\right\| \sqrt{\sum\limits_{i=1}^p\left(\hat{v}^*_{n,ii}v_{\theta_0,ii}-1\right)^2}$. Therefore, combined with \eqref{lower bound of q*Mn} and \eqref{v*vtheta0 VBE}, we have
\begin{eqnarray}
\label{final lower bound of q*Mn}
& & \int q^{*,\Theta_0\left(r_0\right)}(\theta) M_n(\theta ; x) \mathrm{d} \theta \nonumber\\
&\geq& M_n(\theta^*;x)+\frac{1}{2}\Delta_{n,\theta_0}^{\top}V_{\theta_0}\Delta_{n,\theta_0} \nonumber\\
& & -C_3\left(p+ R_n +\frac{p}{\min v_{\theta_0,ii}}\left[\left(\frac{p}{n\min v_{\theta_0,ii}}\right)^{\frac{s_1}{2}}+\left(\frac{p}{n\min v_{\theta_0,ii}}\right)^{\frac{s_2}{2}}\right]\right),
\end{eqnarray}
with probability at least $1-5e^{-y}-e^{-cp}-e^{-cr_0^2}$ for a certain constant $C_3>0$. In summary, we have established both upper and lower bounds for $\int q^{*,\Theta_0\left(r_0\right)}(\theta) M_n(\theta ; x) \mathrm{d} \theta$.

{\textit{Step 2.}} Derive a non-asymptotic upper bound for $\int_{B^c\left(\theta^*, \eta\right)} q^{*,\Theta_0\left(r_0\right)}(\theta)\mathrm{d} \theta$.

Combining \eqref{upper bound of q*Mn} and \eqref{final lower bound of q*Mn}, we have
\begin{eqnarray}
\label{consistency of q* renormalized}
& & \int_{B^c\left(\theta^*, \eta\right)} q^{*,\Theta_0\left(r_0\right)}(\theta)\mathrm{d} \theta  \nonumber\\
&\lesssim& \frac{p+R_n+\frac{p}{\min v_{\theta_0,ii}}\left[\left(\frac{p}{n\min v_{\theta_0,ii}}\right)^{\frac{s_1}{2}}+\left(\frac{p}{n\min v_{\theta_0,ii}}\right)^{\frac{s_2}{2}}\right]}{\lambda_{\min}\left(V_{\theta_0}\right)n\eta^2+\Delta_{n, \theta_0}^{\top}V_{\theta_0}\Delta_{n, \theta_0}}
\end{eqnarray}
holds with probability at least $1-5e^{-y}-e^{-cp}-e^{-cr_0^2}$. 

{\textit{Step 3.}} Extend the result from $q^{*,\Theta_0\left(r_0\right)}(\theta)$ to $q^*(\theta)$. 

Similar to the proof of \eqref{upper bound of the edge mass of the bridge}, we can neglect the impact of renormalization from $\Theta_0(r_0)$, demonstrating that
\begin{equation}
\label{consistency of q*}
\int_{B^c\left(\theta^*, \eta\right)} q^*(\theta)\mathrm{d} \theta\lesssim\frac{p+R_n+\frac{p}{\min v_{\theta_0,ii}}\left[\left(\frac{p}{n\min v_{\theta_0,ii}}\right)^{\frac{s_1}{2}}+\left(\frac{p}{n\min v_{\theta_0,ii}}\right)^{\frac{s_2}{2}}\right]}{\lambda_{\min}\left(V_{\theta_0}\right)n\eta^2+\Delta_{n, \theta_0}^{\top}V_{\theta_0}\Delta_{n, \theta_0}}
\end{equation}
holds with probability at least $1-5e^{-y}-e^{-cp}-e^{-cr_0^2}$.

We next use a similar analytical technique as that used in Theorem \ref{theorem: the asymptotic normality of the VB bridge} to derive a upper bound for the TV distance between $q_{\tilde{\theta}}^*(\tilde{\theta})$ and $N\left(\tilde{\theta};\Delta_{n, \theta_0},\mathrm{diag}^{-1}(V_{\theta_0})\right)$.

Recall that 
\begin{eqnarray*}
& & \mathrm{KL}(q(\theta ;\mathbf{m},\mathbf{v}) \| q(\theta ;\mathbf{m+\Delta m},\mathbf{v+\Delta v})\\
&=& \frac{n}{2}\sum\limits_{i=1}^p\frac{\left(\Delta m_i\right)^2}{v_i+\Delta v_i}+\frac{1}{2}\sum\limits_{i=1}^p\left\{-\frac{\Delta v_i}{v_i+\Delta v_i}+\log\left(1+\frac{\Delta v_i}{v_i}\right)\right\}.    
\end{eqnarray*}
Substituting $\mathbf{m},\mathbf{v},\mathbf{\Delta m},\mathbf{\Delta v}$ with $\hat{\theta}_n^*,\left(\hat{v}^*_{n,11},\cdots,\hat{v}^*_{n,pp}\right)^{\top},\mathbf{m}_0^*-\hat{\theta}_n^*,\mathbf{v}_0^*-\left(\hat{v}^*_{n,11},\cdots,\hat{v}^*_{n,pp}\right)^{\top}$, then
\begin{eqnarray*}
 \mathrm{KL}(q^*(\theta)\|q_0(\theta)) &=& \frac{1}{2}\sum\limits_{i=1}^p\left(\hat{v}^*_{n,ii}v_{\theta_0,ii}-\log \hat{v}^*_{n,ii}v_{\theta_0,ii}-1\right) \\
& & + \frac{n}{2}\left(\hat{\theta}_n^*-\theta^*-\frac{\Delta_{n, \theta_0}}{\sqrt{n}}\right)^{\top}\mathrm{diag}\left(V_{\theta_0}\right)\left(\hat{\theta}_n^*-\theta^*-\frac{\Delta_{n, \theta_0}}{\sqrt{n}}\right).
\end{eqnarray*}
By Theorem \ref{theorem: the convergence rates of the distribution parameters of the VB posterior}, this indicates that
\begin{eqnarray}
\label{KL q* q0}
& & \mathrm{KL}(q^*(\theta)\|q_0(\theta)) \nonumber\\
&\lesssim& \frac{\max_i(v_{\theta_0,ii})}{\lambda_{\min}\left(V_{\theta_0}\right)}\left(R_n+\frac{p}{\min v_{\theta_0,ii}}\left[\left(\frac{p}{n\min v_{\theta_0,ii}}\right)^{\frac{s_1}{2}}+\left(\frac{p}{n\min v_{\theta_0,ii}}\right)^{\frac{s_2}{2}}\right]\right)
\end{eqnarray}
holds with probability at least $1-5e^{-y}-e^{-cp}-e^{-cr_0^2}$. By utilizing the invariance of the KL divergence under the linear transformation $\tilde{\theta}=\sqrt{n}\left(\theta-\theta^*\right)$ in Gaussian distributions, we can show that
\begin{equation}
\label{KL invariance of q*}
\mathrm{KL}\left(q_{\tilde{\theta}}^*(\tilde{\theta})\|\underset{q \in \mathcal{Q}^p}{\arg \min } \operatorname{KL}\left(q(\tilde{\theta}) \| N\left(\tilde{\theta}; \Delta_{n, \theta_0}, V_{\theta_0}^{-1}\right)\right)\right)=\mathrm{KL}(q^*(\theta)\|q_0(\theta)).
\end{equation}
From the Pinsker inequality, it follows that
\begin{eqnarray}
\label{Pinsker of q*}
& & \left\|q_{\tilde{\theta}}^*(\tilde{\theta})-\underset{q \in \mathcal{Q}^p}{\arg \min } \operatorname{KL}\left(q(\tilde{\theta}) \| N\left(\tilde{\theta}; \Delta_{n, \theta_0}, V_{\theta_0}^{-1}\right)\right)\right\|_{\mathrm{TV}} \nonumber\\
&\lesssim& \sqrt{\mathrm{KL}\left(q_{\tilde{\theta}}^*(\tilde{\theta})\|\underset{q \in \mathcal{Q}^p}{\arg \min } \operatorname{KL}\left(q(\tilde{\theta}) \| N\left(\tilde{\theta}; \Delta_{n, \theta_0}, V_{\theta_0}^{-1}\right)\right)\right)}.
\end{eqnarray}
Therefore, combining \eqref{KL q* q0}--\eqref{Pinsker of q*}, we derive that
\begin{eqnarray}
\label{TV of q* - N}
& & \left\|q_{\tilde{\theta}}^*(\tilde{\theta})-N\left(\tilde{\theta};\Delta_{n, \theta_0},\mathrm{diag}^{-1}(V_{\theta_0})\right)\right\|^2_{\mathrm{TV}} \nonumber\\
&\lesssim& \frac{\max\nolimits_{i}v_{\theta_0,ii}}{\lambda_{\min}\left(V_{\theta_0}\right)}\left\{R_n+\frac{p}{\min v_{\theta_0,ii}}\left[\left(\frac{p}{n\min v_{\theta_0,ii}}\right)^{\frac{s_1}{2}}+\left(\frac{p}{n\min v_{\theta_0,ii}}\right)^{\frac{s_2}{2}}\right]\right\}
\end{eqnarray}
holds with probability at least $1-5e^{-y}-e^{-cp}-e^{-cr_0^2}$. In summary, we have established the non-asymptotic Variational Bernstein-von Mises Theorem by synthesizing the \eqref{consistency of q*} and \eqref{TV of q* - N}.

\section{Proof of Corollary \ref{corollary of theorem 5}}
\label{Proof of Corollary 3}

According to the proof of Corollary \ref{corollary of theorem 4}, we have $R_n=O_p\left(\sqrt{p^3/n}\right)$. Consequently, by Theorem \ref{theorem: the non-asymptotic variational Bernstein--von Mises theorem} (1), the following holds:
\begin{equation*}
\int_{B^c\left(\theta^*, \eta\right)} q^*(\theta)\mathrm{d} \theta=O_p\left(\frac{p+\sqrt{\frac{p^3}{n}}+p\left[\left(\frac{p}{n}\right)^{\frac{s_1}{2}}+\left(\frac{p}{n}\right)^{\frac{s_2}{2}}\right]}{n\eta^2}\right)=O_p\left(\frac{p}{n\eta^2}\right).
\end{equation*}
Furthermore, if $\max_i(v_{\theta_0, ii})\lesssim 1$ and $\min\left(s_1,s_2\right) \geq 1$, Theorem \ref{theorem: the non-asymptotic variational Bernstein--von Mises theorem} (2) ensures that
\begin{equation*}
\left\|q_{\tilde{\theta}}^*(\tilde{\theta})-N\left(\tilde{\theta};\Delta_{n, \theta_0},\mathrm{diag}^{-1}(V_{\theta_0})\right)\right\|^2_{\mathrm{TV}}\lesssim R_n+p\left[\left(\frac{p}{n}\right)^{\frac{s_1}{2}}+\left(\frac{p}{n}\right)^{\frac{s_2}{2}}\right]\lesssim R_n+p\left(\frac{p}{n}\right)^{\frac{1}{2}},
\end{equation*}
implying that:
\begin{equation*}
\left\|q_{\tilde{\theta}}^*(\tilde{\theta})-N\left(\tilde{\theta};\Delta_{n, \theta_0},\mathrm{diag}^{-1}(V_{\theta_0})\right)\right\|^2_{\mathrm{TV}}=O_P\left(\sqrt{\frac{p^3}{n}}\right).
\end{equation*}

\section{Proof of Proposition \ref{m=logp}}
\label{Proof of Proposition 1}

For the Gaussian mixture models (GMMs) defined in \eqref{abstract GMM}, we have 
\begin{eqnarray*}
m(\mu;x_i) &=& \sup\limits_{q(c_i) \in \mathcal{Q}} \int q(c_i) \log \frac{p(x_i,c_i\mid \mu)}{q(c_i)}\mathrm{d}c_i\\
&=& \sup_{\sum\limits_{k=1}^K q_k=1,q_k> 0}\sum_{k=1}^K q_k \log\left\{ \frac{K^{-1}(2\pi)^{-\frac{p}{2}} \exp \left(-\frac{\left\|x_i-\mu_k\right\|^2}{2}\right)}{q_k}\right\},
\end{eqnarray*}
where $c_i$ follows a categorical distribution in a variational framework over $K$ classes with assignments $\left(q_1, \ldots, q_K\right)$. Applying the Lagrange Multiplier method, we consider the function
\begin{equation*}
\begin{aligned}
&F\left(q_1,\ldots,q_K,\lambda\right)\\=&-\log K -\frac{p}{2} \log (2 \pi)-\frac{1}{2} \sum_{k=1}^K q_k\left\|x_i-\mu_k\right\|^2-\sum_{k=1}^K q_k \log q_k+\lambda\left(\sum_{k=1}^K q_k -1\right).
\end{aligned}
\end{equation*}
It indicates that the variational parameters are determined by
\begin{equation*}
\frac{\partial F\left(q_1,\ldots,q_K,\lambda\right)}{\partial q_k}=-\frac{1}{2}\left\|x_i-\mu_k\right\|^2-\left(1+\log q_k\right)+\lambda=0,
\end{equation*}
yielding the solutions:
\begin{equation*}
q_k^*=\frac{\exp \left(-\frac{1}{2}\left\|x_i-\mu_k\right\|^2\right)}{\sum\limits_{l=1}^K \exp \left(-\frac{1}{2}\left\|x_i-\mu_l\right\|^2\right)}
\end{equation*}
for $k=1,\ldots, K$. Substituting $q_k^*$ into $m(\mu;x_i)$, We obtain 
\begin{equation*}
m(\mu;x_i)=-\log K -\frac{p}{2} \log (2 \pi)+\log\left(\sum\limits_{k=1}^K \exp \left\{-\frac{1}{2}\left\|x_i-\mu_k\right\|^2\right\}\right),
\end{equation*}
 which turns out to coincide with $\log p(x_i\mid\mu)$, showing that the variational log-likelihood simplifies to the standard log-likelihood.

\section{Proof of Proposition \ref{mustar=mu0}}
\label{Proof of Proposition 2}

By Proposition \ref{m=logp}, any maximizer of $\mathbb{E}_{\mu^0}m(\mu;x_i)$ is also a maximizer of $\mathbb{E}_{\mu^0}\log p(x_i;\mu)$. To show that $\mu^*=\mu^0$ up to a permutation, we first establish the integrability of both  $\mathbb{E}_{\mu^0}\nabla m(\mu;x_i)$ and $\mathbb{E}_{\mu^0}\nabla^2 m\left(\mu ; x_i\right)$. Then, we obtain the result by applying the dominated convergence theorem.

Let
\begin{equation*}
w_k(\mu; x_i)=\frac{\exp\left(-\frac{1}{2}\left\|x_i-\mu_k\right\|^2\right)}{\sum_{l=1}^K\exp\left(-\frac{1}{2}\left\|x_i-\mu_l\right\|^2\right)},
\end{equation*}
$w(\mu ; x_ i) = \left(w_1(\mu; x_i),\ldots,w_K(\mu; x_i)\right)^\top$, and $Z(\mu; x_i) = \left(x_i-\mu_1,\ldots,x_i-\mu_K\right)$. Then, the gradient vector is given by
\begin{equation*}
\nabla m(\mu;x_i)=w(\mu ; x_ i)\otimes 1_p \circ \operatorname{Vec}\left(Z(\mu; x_i)\right),
\end{equation*}
where $1_p$ is a $p$-dimensional vector of ones, $\otimes$ denotes the Kronecker product, $\circ$ denotes the Hadamard product, and $\operatorname{Vec}(\cdot)$ is the vectorization operator. Since $0 < w_k(\mu; x_i)< 1$, by Cauchy-Schwarz inequality, we have
\begin{equation*}
 \mathbb{E}_{\mu^0} \|\nabla m(\mu;x_i)\| \leq \sqrt{\mathbb{E}_{\mu^0} \left\| \nabla m(\mu; x_i) \right\|^2} < \sqrt{\mathbb{E}_{\mu^0} \left\| \mathrm{Vec}(Z(\mu; x_i)) \right\|^2}
= \sqrt{\sum_{k=1}^{K} \mathbb{E}_{\mu^0} \| x_i - \mu_k \|^2 }.
\end{equation*}
Note that
\begin{equation*}
\mathbb{E}_{\mu^0} \|x_i - \mu_k\|^2 = \frac{1}{K} \sum^K_{k=1} \sum_{k=1}^{K} \mathrm{tr} \left( \mathbb{E}_{x_i \sim \mathcal{N}(\mu_k^0, I_p)} \left\{(x_i - \mu_k)(x_i - \mu_k)^\top\right\} \right) < \infty.
\end{equation*}
Hence, the integrability of $\mathbb{E}_{\mu^0}\nabla m(\mu;x_i)$ is guaranteed. Consequently, we obtain the interchangeability of differentiation and expectation via the dominated convergence theorem. Thus,
\begin{equation}
\label{Peq1-xpr}
\nabla{\mathbb{E}_{\mu^0}} m(\mu^0; x_i) = \mathbb{E}_{\mu^0}\nabla m(\mu^0;x_i) =  {\mathbb{E}_{\mu^0}} \nabla\log p(x_i; \mu^0) = \mathbb{E}_{\mu^0} \frac{\nabla p(x_i; \mu^0)}{p(x_i; \mu^0)} = 0.
\end{equation}

The Hessian matrix is given by
\begin{eqnarray*}
\nabla^2 m\left(\mu ; x_i\right) &=& -  \operatorname{diag}\left(w\left(\mu ; x_i\right)\right) \otimes I_p +\sum_{k_1<k_2} w_{k_1}\left(\mu ; x_i\right) w_{k_2}\left(\mu ; x_i\right) \\
& & \times \operatorname{Vec}\left\{Z\left(\mu ; x_i\right) \operatorname{diag}\left(e_{k_1}-e_{k_2}\right)\right\}\operatorname{Vec}\left\{Z\left(\mu ; x_i\right) \operatorname{diag}\left(e_{k_1}-e_{k_2}\right)\right\}^{\top}\\
&\triangleq& -T_1 + T_2,
\end{eqnarray*}
where $e_k$ is a $K$-dimensional standard basis vector with a 1 in the $k$-th position and 0s elsewhere. On the one hand, it is easy to verify that $0_{Kp \times Kp} < T_1 <  I_{Kp}$. On the other hand, since 
\begin{equation*}
    \mathbb{E}_{\mu^0} \left\{\mathrm{Vec} \left(Z(\mu; x_i) \, \mathrm{diag}(e_{k_1} - e_{k_2})\right) 
\mathrm{Vec} \left(Z(\mu; x_i) \, \mathrm{diag}(e_{k_1} - e_{k_2})\right)^\top\right\} < \infty,
\end{equation*}
we have $0_{Kp \times Kp} < \mathbb{E}_{\mu^0} (T_2) < \infty$. Thus, the integrability of $\mathbb{E}_{\mu^0}\nabla^2 m\left(\mu ; x_i\right)$ is guaranteed. Again, by the dominated convergence theorem, we have
\begin{equation}
\label{Peq2-xpr}
\nabla^2\mathbb{E}_{\mu^0}m(\mu^0;x_i)= \mathbb{E}_{\mu^0}\nabla^2m(\mu^0;x_i) = \mathbb{E}_{\mu^0}\nabla^2 \log p(x_i;\mu^0)=-\operatorname{Var}_{\mu^0}\left(\nabla \log p(x_i;\mu^0)\right).
\end{equation}
Given that $V_{\mu^0}=\operatorname{Var}_{\mu^0}\left(\nabla \log p(x_i;\mu^0)\right)>0$, it follows from \eqref{Peq1-xpr} and \eqref{Peq2-xpr} that $\mu^*$ coincides with $\mu_0$ up to a permutation.

\section{Proof of Proposition \ref{Proposition: the (variational) information matrix}}
\label{Proof of Proposition 3}

By the dominated convergence theorem, we have
\begin{eqnarray*}
V_{\mu^0} &=& -\nabla^2 \mathbb{E}_{\mu^0} m\left(\mu^0 ; x_i\right) \\
&=& \mathbb{E}_{\mu^0} \operatorname{diag}\left(w\left(\mu^0 ; x_i\right)\right) \otimes I_p-\sum_{k_1<k_2} \mathbb{E}_{\mu^0}\left[w_{k_1}\left(\mu^0 ; x_i\right) w_{k_2}\left(\mu^0 ; x_i\right) \right. \\
& & \left.\times\operatorname{Vec}\left\{Z\left(\mu^0 ; x_i\right) \operatorname{diag}\left(e_{k_1}-e_{k_2}\right)\right\}\operatorname{Vec}\left\{Z\left(\mu^0 ; x_i\right) \operatorname{diag}\left(e_{k_1}-e_{k_2}\right)\right\}^{\top}\right].
\end{eqnarray*}
Note that $\mathbb{E}_{\mu^0}w_k(\mu^0; x_i)=1/K$. Hence, we obtain
\begin{equation*}
\mathbb{E}_{\mu^0}\operatorname{diag}\left(w(\mu^0;x_ i)\right)\otimes I_p=\frac{1}{K}I_{K}\otimes I_p=\frac{1}{K}I_{Kp}.
\end{equation*}
Consequently, for any eigenvalue $\lambda\left(V_{\mu^0}\right)$ of $V_{\mu^0}$, we have
\begin{eqnarray*}
0 &\leq& \frac{1}{K}  -\lambda\left(V_{\mu^0}\right) \\
&\leq& \operatorname{tr}  \left(\sum _ { k _ { 1 } < k _ { 2 } } \mathbb { E } _ { \mu ^ { 0 } } \left[w_{k_1}\left(\mu^0 ; x_i\right) w_{k_2}\left(\mu^0 ; x_i\right)\right.\right.\\
& & \left.\left.\times \operatorname{Vec}\left\{Z\left(\mu^0 ; x_i\right) \operatorname{diag}\left(e_{k_1}-e_{k_2}\right)\right\}\operatorname{Vec}\left\{Z\left(\mu^0 ; x_i\right) \operatorname{diag}\left(e_{k_1}-e_{k_2}\right)\right\}^{\top}\right]\right) \\
&=& \sum_{k_1<k_2}\operatorname{tr}\left(\int \frac{\exp \left(-\frac{1}{2}\left\|x_i-\mu_{k_1}^0\right\|^2-\frac{1}{2}\left\|x_i-\mu_{k_2}^0\right\|^2\right)}{K(2 \pi)^{\frac{p}{2}}\left[\sum_{k=1}^K \exp \left(-\frac{1}{2}\left\|x_i-\mu_k^0\right\|^2\right)\right]}\right.\\
& &\left.\times\operatorname{Vec}\left\{Z\left(\mu^0 ; x_i\right) \operatorname{diag}\left(e_{k_1}-e_{k_2}\right)\right\}\operatorname{Vec}\left\{Z\left(\mu^0 ; x_i\right) \operatorname{diag}\left(e_{k_1}-e_{k_2}\right)\right\}^{\top} \mathrm{d}x_i\right).
\end{eqnarray*}
Note that
\begin{eqnarray*}
& & \frac{\exp \left(-\frac{1}{2}\left\|x_i-\mu_{k_1}^0\right\|^2-\frac{1}{2}\left\|x_i-\mu_{k_2}^0\right\|^2\right)}{\sum_{l=1}^K \exp \left(-\frac{1}{2}\left\|x_i-\mu_l^0\right\|^2\right)}\\
&\leq& \frac{\exp \left(-\frac{1}{2}\left\|x_i-\mu_{k_1}^0\right\|^2-\frac{1}{2}\left\|x_i-\mu_{k_2}^0\right\|^2\right)}{\exp \left(-\frac{1}{2}\left\|x_i-\mu_{k_1}^0\right\|^2\right)+\exp \left(-\frac{1}{2}\left\|x_i-\mu_{k_2}^0\right\|^2\right)} \leq \frac{\exp \left(-\frac{1}{2}\left\|x_i-\frac{\mu_{k_1}^0+\mu_{k_2}^0}{2}\right\|^2\right)}{2\exp \left(\frac{1}{8}\left\|\mu_{k_1}^0-\mu_{k_2}^0\right\|^2\right)}.
\end{eqnarray*}
Hence, we have 
\begin{eqnarray*}
0 &\leq& \frac{1}{K}  -\lambda\left(V_{\mu^0}\right) 
\\
&\leq& \sum_{k_1<k_2} \frac{1}{2K}\exp \left(-\frac{1}{8}\left\|\mu_{k_1}^0-\mu_{k_2}^0\right\|^2\right) \\
& & \times \operatorname{tr}\left(\mathbb{E}_{x_i \sim N\left(\frac{1}{2}(\mu_{k_1}^0+\mu_{k_2}^0), I_p\right)}\left[\operatorname{Vec}\left\{Z\left(\mu^0 ; x_i\right) \operatorname{diag}\left(e_{k_1}-e_{k_2}\right)\right\}\operatorname{Vec}\left\{Z\left(\mu^0 ; x_i\right) \operatorname{diag}\left(e_{k_1}-e_{k_2}\right)\right\}^{\top}\right]\right)\\
&=& \sum_{k_1<k_2} \frac{\left\|\mu_{k_1}^0-\mu_{k_2}^0\right\|^2+4 p}{4 K \exp \left(\frac{1}{8}\left\|\mu_{k_1}^0-\mu_{k_2}^0\right\|^2\right)}\leq\frac{C}{K}.
\end{eqnarray*}
Additionally, for any $r \in \{1,\ldots, Kp\}$,
\begin{equation*}
\lambda_{\min}\left(V_{\mu^0}\right)=\lambda_{\min}\left(V_{\mu^0}\right)\left\| e_r\right\|^2\leq e_r^\top V_{\mu^0}e_r=\left(V_{\mu^0}\right)_{r,r}\leq\lambda_{\max}\left(V_{\mu^0}\right)\left\| e_r\right\|^2=\lambda_{\max}\left(V_{\mu^0}\right),
\end{equation*}
where $e_r$ is a $Kp$-dimensional standard basis vector with a 1 in the $r$-th position and 0s elsewhere. Therefore, we conclude that 
\begin{equation*}
\frac{1-C}{K}\leq\lambda_{\min}\left(V_{\mu^0}\right) \leq \left(V_{\mu^0}\right)_{r,r} \leq \lambda_{\max}\left(V_{\mu^0}\right) \leq \frac{1}{K}.
\end{equation*}

\section{Some Verifications of Design Conditions}
\label{Some Verifications of Design Conditions}

We first verify Condition \ref{condition: local, twice continuously differentiable}. 
According to the proof of Proposition \ref{Proposition: the (variational) information matrix}, we know that
\begin{equation*}
\frac{1-C}{K}\leq\lambda\left(V_{\mu^0}\right)\leq\frac{1}{K},
\end{equation*}
for any eigenvalue $\lambda\left(V_{\mu^0}\right)$ of $V_{\mu^0}$. Consequently, $\lambda\left(D^2_0\right)=n\lambda\left(V_{\mu^0}\right)$, which leads to
\begin{equation*}
\Theta_0(r_0)\subseteq\left\{\mu\in\mathbb{R}^{p\times K}:\left\|\operatorname{Vec}\left(\mu-\mu^0\right)\right\| \leq \sqrt{\frac{K}{1-C}}\frac{r_0}{\sqrt{n}}\right\}.
\end{equation*}
It follows that for $\mu\in\Theta_0(r_0)$, $\left\|\operatorname{Vec}\left(\mu-\mu^0\right)\right\|=o(1)$. Thus, for $\mu\in\Theta_0(r_0)$, we have
\begin{eqnarray*}
\rho\left(D_0^{-1} D^2(\mu) D_0^{-1}-I_p\right) &=& \rho\left(D_0^{-1}\left(D^2(\mu)-D_0^2\right)D_0^{-1}\right)\\
&\leq& \frac{1}{n}\left(\frac{K}{1-C}\right)^2\rho\left(D^2(\mu)-D_0^2\right)\\
&=& \left(\frac{K}{1-C}\right)^2\rho\left(\nabla^2 \mathbb{E}_{\mu^0}m\left(\mu;x_i\right)-\nabla^2 \mathbb{E}_{\mu^0}m\left(\mu^0;x_i\right)\right)=o(1),
\end{eqnarray*}
because the variational log-likelihood process is three times stochastically differentiable. For related discussions on three times stochastic differentiability, see Section 5.1 of \cite{Spokoiny2012} and Section 2.4.3 of \cite{spokoiny2014bernsteinvonmises}.

As for Condition \ref{condition: local, exponential moment}, given that $\tau_2$ can be sufficiently small, for all $\gamma \in \mathbb{R}^{Kp}$ and $|\lambda| \leq \tau_2$, we use the Taylor expansion to derive the following inequality:
\begin{eqnarray*}
& & \mathbb{E}_{\mu^0}\exp \left[\lambda \frac{\gamma^{\top} \nabla M_n\left(\mu^0 ; x\right)}{\left\|\Sigma_0 \gamma\right\|}\right]\\
&\leq& 1+\mathbb{E}_{\mu^0}\left[\lambda \frac{\gamma^{\top} \nabla M_n\left(\mu^0 ; x\right)}{\left\|\Sigma_0 \gamma\right\|}\right]+ \tau_1\mathbb{E}_{\mu^0}\left[\frac{\lambda^2}{2} \frac{\gamma^{\top} \nabla M_n\left(\mu^0 ; x\right)\nabla M_n\left(\mu^0 ; x\right)^\top\gamma}{\left\|\Sigma_0 \gamma\right\|^2}\right]\\
&=& 1  +\frac{\tau_1\lambda^2}{2}\frac{\gamma^\top \mathbb{E}_{\mu^0}\left[\nabla M_n\left(\mu^0 ; x\right)\nabla M_n\left(\mu^0 ; x\right)^\top\right]\gamma}{\gamma^\top\Sigma_0^2\gamma}= 1+\frac{C\lambda^2}{2} \leq \exp\left(\frac{C}{2}\lambda^2\right)
\end{eqnarray*}
for some constant $\tau_1>1$. 

For the global exponential moment Condition \ref{condition: global, exponential moment, variational}, similar to Condition \ref{condition: local, exponential moment}, we need to verify the following inequality: for all $|\lambda| \leq \delta_2(r)$, 
\begin{equation}
\label{PC3}
\sup _{\mu\in \Theta_0(r)}\sup _{\gamma_1, \gamma_2 \in \mathbb{R}^{Kp}} \mathbb{E}_{\mu^0} \exp \left[\frac{\lambda}{\tau_3} \frac{\gamma_1^{\top} \nabla^2\left\{M_n\left(\mu; x\right)-\mathbb{E}_{\theta_0} M_n\left(\mu; x\right)\right\}\gamma_2}{\left\|D_0 \gamma_1\right\| \cdot\left\|D_0 \gamma_2\right\|}\right] \leq \exp\left(\frac{\tau_2^2 \lambda^2}{2}\right),
\end{equation}
where $\tau_3$ is a positive constant, and $\delta_2(r)\geq Cp$ for all $r$ with a fixed constant $C>0$. The constant $\tau_3$ can be chosen sufficiently large to meet our requirements. We note that
\begin{eqnarray*}
& & \mathbb{E}_{\mu^0} \exp \left[\frac{\lambda}{\tau_3} \frac{\gamma_1^{\top} \nabla^2\left\{M_n\left(\mu; x\right)-\mathbb{E}_{\theta_0} M_n\left(\mu; x\right)\right\}\gamma_2}{\left\|D_0 \gamma_1\right\| \cdot\left\|D_0 \gamma_2\right\|}\right]\\
&=& \prod_{i=1}^n\mathbb{E}_{\mu^0} \exp \left[\frac{\lambda}{\tau_3} \frac{\gamma_1^{\top} \nabla^2\left\{m\left(\mu; x_i\right)-\mathbb{E}_{\theta_0} m\left(\mu; x_i\right)\right\}\gamma_2}{\left\|D_0 \gamma_1\right\| \cdot\left\|D_0 \gamma_2\right\|}\right]\\
&=& \prod_{i=1}^n \frac{1}{K} \sum_{k=1}^K\mathbb{E}_{x_i \sim N\left(\mu_k^0, I_p\right)}\exp \left[\frac{\lambda}{\tau_3} \frac{\gamma_1^{\top} \nabla^2\left\{m\left(\mu; x_i\right)-\mathbb{E}_{\theta_0} m\left(\mu; x_i\right)\right\}\gamma_2}{\left\|D_0 \gamma_1\right\| \cdot\left\|D_0 \gamma_2\right\|}\right].
\end{eqnarray*}
Since $\lambda\left(D^2_0\right)=n\lambda\left(V_{\mu^0}\right)$, by Proposition \ref{Proposition: the (variational) information matrix} we know that $\left\|D_0 \gamma\right\|\asymp\sqrt{n}\left\|\gamma\right\|$ and thus $\left\|D_0 \gamma_1\right\| \cdot\left\|D_0 \gamma_2\right\|\asymp n\left\|\gamma_1\right\|\cdot\left\|\gamma_2\right\|$. On the other hand, by tedious calculations, we can  decompose the quadratic form $\gamma_1^{\top} \nabla^2\left\{m\left(\mu; x_i\right)-\mathbb{E}_{\theta_0} m\left(\mu; x_i\right)\right\}\gamma_2$ based on the corresponding blocks for the $K$ components, and thus show that 
\begin{equation*}
    \frac{\gamma_1^{\top} \nabla^2\left\{m\left(\mu; x_i\right)-\mathbb{E}_{\theta_0} m\left(\mu; x_i\right)\right\}\gamma_2}{\left\|D_0 \gamma_1\right\| \cdot\left\|D_0 \gamma_2\right\|}
\end{equation*}
is bounded. Hence, by applying the mean value theorem and completing the square, we can formally verify \eqref{PC3}. We also recommend employing advanced mathematical analysis techniques and numerical simulations for a more thorough verification of Condition \ref{condition: global, exponential moment, variational}.

Next, we examine Condition \ref{condition: global identification property}: for each $r\geq r_0$, there exists a constant $C>0$ such that, for $\mu$ satisfying
\begin{equation*}
\inf _{\pi \in S_K}\left\|D_0 \operatorname{Vec}\left(\mu \pi-\mu^0\right)\right\|=r,
\end{equation*}
we have
\begin{equation}
\label{GMM KL lower bound}
\mathbb{E}_{\mu^0}m\left(\mu^0;x_i\right)-\mathbb{E}_{\mu^0}m\left(\mu;x_i\right)\geq\frac{Cr}{n},
\end{equation}
where $\pi$ is a permutation matrix, and $S_K$ is the symmetric group of all permutations of $K$ elements. The inclusion of $\left(\pi, S_K\right)$ addresses potential identifiability issues. We now discuss this condition in two cases. 

{\textit{Case 1.}} When the distance between $\mu$ and $\mu^0$ is relatively small, considering the Taylor expansion, we obtain
\begin{eqnarray*}
\mathbb{E}_{\mu^0}m\left(\mu^0;x_i\right)-\mathbb{E}_{\mu^0}m\left(\mu;x_i\right) &\gtrsim& \operatorname{Vec}\left(\mu-\mu^0\right)^\top V_{\mu^0}\operatorname{Vec}\left(\mu-\mu^0\right)\\
&=& \frac{1}{n}\left\|D_0 \operatorname{Vec}\left(\mu-\mu^0\right)\right\|^2 \geq\frac{r^2}{n}\gtrsim\frac{r}{n}
\end{eqnarray*}
due to $r\geq r_0$. 

{\textit{Case 2.}} When $\mu$ and $\mu^0$ are far apart, under Proposition \ref{m=logp}, we have
\begin{eqnarray*}
\mathbb{E}_{\mu^0}m\left(\mu^0;x_i\right)-\mathbb{E}_{\mu^0}m\left(\mu;x_i\right) &=& \mathbb{E}_{\mu^0}\log p\left(x_i;\mu^0\right)-\mathbb{E}_{\mu^0}\log p\left(\mu;x_i\right)\\
&=& \operatorname{KL}\left(p\left(x_i;\mu^0\right) \| p(x_i;\mu)\right).
\end{eqnarray*}
Note that $p(x_i;\mu)=p(x_i;\mu\pi)$ due to cyclic symmetry. Hence, the left side of \eqref{GMM KL lower bound} becomes the KL divergence between two Gaussian mixtures, which lacks an analytic form. This is a key topic in information theory, with numerous studies focusing on its bounds and approximations. For example, Section 3 of \cite{durrieu2012lower} provided a strict but relaxed lower bound for $\operatorname{KL}\left(p\left(x_i;\mu^0\right) \| p(x_i;\mu\pi)\right)$ as follows:
\begin{eqnarray*}
& & \operatorname{KL}\left(p\left(x_i;\mu^0\right) \| p(x_i;\mu\pi)\right)\\
&\geq& -\frac{p}{2}\left(1-\log2\right)\\
& & +\sum_{k=1}^K\frac{1}{K} \log\left(\sum_{k^\prime=1}^K\exp\left\{-\frac{1}{2}\left\|\mu_{k^\prime}^0-\mu_k^0\right\|^2\right\}\right) -\log\left(\sum_{k^\prime=1}^K\exp\left\{-\frac{1}{4}\left\|\mu_{\pi\left(k^\prime\right)}-\mu_k^0\right\|^2\right\}\right),
\end{eqnarray*}
and Proposition 3 of \cite{2017On} improved the lower bound by mediant inequality such that
\begin{eqnarray*}
& & \operatorname{KL}\left(p\left(x_i;\mu^0\right) \| p(x_i;\mu)\right)\\
&\geq& \frac{1}{K}\sum_{k=1}^K \mathbb{E}_{x_i \sim N\left(\mu^0_k,I_p\right)}\log\left(\max_\pi \min _{k^\prime\in \{1,\ldots,K\}} \frac{\exp\left\{\mu^{0\top}_{k^\prime}x_i-\frac{1}{2}\left\|\mu^{0}_{k^\prime}\right\|^2\right\}}{\exp\left\{\mu^\top_{\pi\left(k^\prime\right)}x_i-\frac{1}{2}\left\|\mu_{\pi\left(k^\prime\right)}\right\|^2\right\}}\right).
\end{eqnarray*}
It requires extensive computation and is difficult to analyze in high-dimensional cases, as indicated in Section 7 of \cite{2016Guaranteed}. For more details and other bounds or approximations, readers are referred to Sections 2.3.3 and 5 of \cite{2016Guaranteed}, Section 5 of \cite{hershey2007approximating}, and Sections 2.2 and 3 of \cite{durrieu2012lower}, including the product of Gaussians approximation and the variational approximation. We only need to reasonably leverage these related studies to complete our verification. Additionally, from another perspective, we can define $y_i=x_i-\mu_k^0$ and derive that
\begin{eqnarray*}
& & \mathbb{E}_{x_i \sim N\left(\mu^0_k,I_p\right)} \log\left(\frac{\sum_{k^\prime=1}^K \exp\left\{\mu^{0\top}_{k^\prime}x_i-\frac{1}{2}\left\|\mu^{0}_{k^\prime}\right\|^2\right\}}{\sum_{k^\prime=1}^K\exp\left\{\mu^\top_{k^\prime}x_i - \frac{1}{2}\left\|\mu_{k^\prime}\right\|^2\right\}}\right)\\&=& \mathbb{E}_{y_i \sim N\left(0_p,I_p\right)}\log\left(\frac{\sum_{k^\prime=1}^K \exp\left\{\mu^{0\top}_{k^\prime}y_i+\mu^{0\top}_{k^\prime}\mu^0_k-\frac{1}{2}\left\|\mu^{0}_{k^\prime}\right\|^2\right\}}{\sum_{k^\prime=1}^K \exp\left\{\mu^\top_{k^\prime}y_i+\mu^\top_{k^\prime}\mu^0_k -\frac{1}{2}\left\|\mu_{k^\prime}\right\|^2\right\}}\right).
\end{eqnarray*}
Then, Proposition 3 of \cite{2017On} implies that 
\begin{equation*}
\operatorname{KL}\left(p\left(x_i;\mu^0\right) \| p(x_i;\mu)\right) \gtrsim \mathbb{E}_{y_i \sim N\left(0_p,I_p\right)}\sum_{k^\prime=1}^K\log\left(\frac{\sum_{k^\prime=1}^K\exp\left\{\mu^{0\top}_{k^\prime}y_i+\mu^{0\top}_{k^\prime}\mu^0_k-\frac{1}{2}\left\|\mu^{0}_{k^\prime}\right\|^2\right\}}{\sum_{k^\prime=1}^K\exp\left\{\mu^\top_{k^\prime}y_i+\mu^\top_{k^\prime}\mu^0_k-\frac{1}{2}\left\|\mu_{k^\prime}\right\|^2\right\}}\right).
\end{equation*}
We thus can utilize Monte Carlo simulations, combined with studies on the bounds and approximations of the Kullback-Leibler divergence in GMMs, to help verify this design condition.

Conditions \ref{condition: identifiability} and \ref{condition: prior} are straightforward due to Proposition \ref{m=logp} and the flexibility in selecting our prior. Proposition \ref{m=logp} allows us to set $\tau_4=1$ for Condition \ref{condition: identifiability}. Using a standard multivariate Gaussian prior lets us assume $\tau_5=1$ for Condition \ref{condition: prior} as a simple example. In practice, however, we hope to choose the most suitable prior tailored to the specific problem to facilitate optimization and inference.

As for Condition \ref{condition: global, exponential moment, complete}, we first note that 
\begin{equation}
\label{hessian of log mu}
\nabla^2\log p(x_i, c_i;\mu)=\nabla^2\left[\log p(x_i;c_i,\mu)+\log p(c_i;\mu)\right]=\nabla^2\log p(x_i;c_i,\mu),
\end{equation}
which is a diagonal matrix whose entries from the $(c_i-1)p+1$-th to the $c_ip$-th positions are -1, with all other elements being zero. Therefore, we can readily deduce that
\begin{equation*}
\nabla^2\log p(x, c;\mu)=\mathbb{E}_{\mu^0}\nabla^2\log p(x, c;\mu),
\end{equation*}
and hence confirm Condition \ref{condition: global, exponential moment, complete}.

To examine the first inequality in Condition \ref{condition: lipschitz}, we first know that
\begin{eqnarray*}
\nabla^2\mathbb{E}_{\mu^0}m\left(\mu ; x_i\right)&=& - \mathbb{E}_{\mu^0}\operatorname{diag}\left(w\left(\mu ; x_i\right)\right) \otimes I_p \\
 & & + \sum_{k_1<k_2}\mathbb{E}_{\mu^0}\left[w_{k_1}\left(\mu ; x_i\right) w_{k_2}\left(\mu ; x_i\right) \right.\\
& & \left.\times \operatorname{Vec}\left\{Z\left(\mu; x_i\right) \operatorname{diag}\left(e_{k_1}-e_{k_2}\right)\right\}\operatorname{Vec}\left\{Z\left(\mu ; x_i\right) \operatorname{diag}\left(e_{k_1}-e_{k_2}\right)\right\}^{\top}\right]\\
&\triangleq& -T_1(\mu; x_i) + T_2(\mu; x_i),
\end{eqnarray*}
based on the proof of Proposition \ref{mustar=mu0}. Then,
\begin{equation*}
\nabla^2\mathbb{E}_{\mu^0}m\left(\mu; x_i\right) - \nabla^2\mathbb{E}_{\mu^0}m\left(\mu^\prime; x_i\right) = \left[T_1(\mu^\prime; x_i) - T_1(\mu; x_i)\right]  + \left[T_2(\mu; x_i) - T_2(\mu^\prime; x_i) \right].
\end{equation*}
Since
\begin{equation*}
0_{Kp\times Kp}<\mathbb{E}_{\mu^0}\operatorname{diag}\left(w\left(\mu ; x_i\right)\right) \otimes I_p<I_{Kp}
\end{equation*}
for all $\mu$, we have 
\begin{equation}
\label{Peq3-xpr}
-I_{Kp}< T_1(\mu^\prime; x_i) - T_1(\mu; x_i) < I_{Kp}    
\end{equation}
We next analyze the term $T_2(\mu; x_i) - T_2(\mu^\prime; x_i)$. For notational simplicity, we define 
\begin{equation*}
    S(\mu; x_i) = \operatorname{Vec}\left\{Z\left(\mu ; x_i\right) \operatorname{diag}\left(e_{k_1}-e_{k_2}\right)\right\}.
\end{equation*}
Then, $T_2(\mu; x_i) - T_2(\mu^\prime; x_i)$ can be decomposed as follows:
\begin{eqnarray*}
& & T_2(\mu; x_i) - T_2(\mu^\prime; x_i) \\ &=& \mathbb{E}_{\mu^0}\left[w_{k_1}\left(\mu ; x_i\right) w_{k_2}\left(\mu ; x_i\right) S(\mu; x_i) S(\mu; x_i)^\top\right] -\mathbb{E}_{\mu^0}\left[w_{k_1}\left(\mu ; x_i\right) w_{k_2}\left(\mu ; x_i\right) S(\mu^\prime; x_i) S(\mu^\prime; x_i)^{\top}\right]\\
 & & + \mathbb{E}_{\mu^0}\left[w_{k_1}\left(\mu ; x_i\right) w_{k_2}\left(\mu ; x_i\right) S(\mu^\prime; x_i) S(\mu^\prime; x_i)^{\top} \right] -\mathbb{E}_{\mu^0}\left[w_{k_1}\left(\mu^\prime ; x_i\right) w_{k_2}\left(\mu^\prime ; x_i\right) S(\mu^\prime; x_i) S(\mu^\prime; x_i)^{\top}\right]\\
 &=& \mathbb{E}_{\mu^0}\left\{w_{k_1}\left(\mu ; x_i\right) w_{k_2}\left(\mu ; x_i\right) \left[S(\mu; x_i) S(\mu; x_i)^\top -  S(\mu^\prime; x_i) S(\mu^\prime; x_i)^{\top}\right]\right\}\\
 & & + \mathbb{E}_{\mu^0}\left\{ \left[w_{k_1}\left(\mu ; x_i\right) w_{k_2}\left(\mu ; x_i\right) - w_{k_1}\left(\mu^\prime ; x_i\right) w_{k_2}\left(\mu^\prime ; x_i\right)\right] S(\mu^\prime; x_i) S(\mu^\prime; x_i)^{\top}\right\}
\end{eqnarray*}
For the first part, $S(\mu; x_i) S(\mu; x_i)^\top -  S(\mu^\prime; x_i) S(\mu^\prime; x_i)^{\top}$ can be decomposed as follows:
\begin{eqnarray*}
   & & S(\mu; x_i) S(\mu; x_i)^\top -  S(\mu^\prime; x_i) S(\mu^\prime; x_i)^{\top}\\ 
   &=& \operatorname{Vec}\left\{\left(\mu^\prime-\mu\right)\operatorname{diag}\left(e_{k_1}-e_{k_2}\right)\right\}\operatorname{Vec}\left\{\left(\mu^\prime-\mu\right)\operatorname{diag}\left(e_{k_1}-e_{k_2}\right)\right\}^\top\\
     & & + \left\{\operatorname{Vec}\left\{\left(\mu^\prime-\mu\right)\operatorname{diag}\left(e_{k_1}-e_{k_2}\right)\right\}S(\mu^\prime; x_i)^{\top} + S(\mu^\prime; x_i) \operatorname{Vec}\left\{\left(\mu^\prime-\mu\right)\operatorname{diag}\left(e_{k_1}-e_{k_2}\right)\right\}^{\top}\right\}\\
     &\triangleq& W(\mu, \mu^\prime) W(\mu, \mu^\prime)^\top + \left\{W(\mu, \mu^\prime) S(\mu^\prime; x_i)^\top + S(\mu^\prime; x_i) W(\mu, \mu^\prime)^\top\right\}.
\end{eqnarray*}
Since $0 < w_k\left(\mu ; x_i\right) < 1$ for all $k=1,\ldots,K$, we have
\begin{equation*}
\rho\left(\mathbb{E}_{\mu^0} \left\{w_{k_1}\left(\mu ; x_i\right) w_{k_2}\left(\mu ; x_i\right) W(\mu, \mu^\prime) W(\mu, \mu^\prime)^\top\right\}\right) \leq \frac{1}{4}\left\|W(\mu, \mu^\prime) \right\|^2 \leq\frac{1}{4}\left\|\operatorname{Vec}\left(\mu^\prime-\mu\right)\right\|^2
\end{equation*}
and
\begin{eqnarray*}
& & \rho\left(\mathbb{E}_{\mu^0}\left\{w_{k_1}\left(\mu; x_i\right) w_{k_2}\left(\mu; x_i\right) \left[W(\mu, \mu^\prime) S(\mu^\prime; x_i)^\top + S(\mu^\prime; x_i) W(\mu, \mu^\prime)^\top\right]\right\}\right)\\
&\leq& \frac{1}{2} \mathbb{E}_{\mu^0} \left\{\left\|W(\mu, \mu^\prime)\right\| \left\|S(\mu^\prime; x_i)\right\|\right\}\\
&\leq& \frac{1}{4}\left\|\operatorname{Vec}\left(\mu^\prime-\mu\right)\right\|^2 + \frac{1}{4} \left[\mathbb{E}_{\mu^0} \left\|S(\mu^\prime; x_i)\right\|\right]^2.
\end{eqnarray*}
And thus, 
\begin{eqnarray*}
& & \rho\left( \mathbb{E}_{\mu^0}\left\{w_{k_1}\left(\mu; x_i\right) w_{k_2}\left(\mu; x_i\right) \left[S(\mu; x_i) S(\mu; x_i)^\top -  S(\mu^\prime; x_i) S(\mu^\prime; x_i)^{\top}\right] \right\} \right) \\
&\leq& \frac{1}{2}\left\|\operatorname{Vec}\left(\mu^\prime-\mu\right)\right\|^2 + \frac{1}{4} \left[\mathbb{E}_{\mu^0} \left\|S(\mu^\prime; x_i)\right\|\right]^2.
\end{eqnarray*}
For the second part, we have
\begin{equation*}
  \rho\left(\mathbb{E}_{\mu^0}\left\{\left[w_{k_1}\left(\mu ; x_i\right) w_{k_2}\left(\mu ; x_i\right)-w_{k_1}\left(\mu^\prime ; x_i\right) w_{k_2}\left(\mu^\prime ; x_i\right)\right] S(\mu^\prime; x_i) S(\mu^\prime; x_i)^\top\right\}\right) \leq   \frac{1}{4}\mathbb{E}_{\mu^0}\left\|S(\mu^\prime; x_i)\right\|^2.
\end{equation*}
Therefore,
\begin{equation*}
    \rho\left(T_2(\mu; x_i) - T_2(\mu^\prime; x_i)\right) \lesssim \left\|\operatorname{Vec}\left(\mu^\prime-\mu\right)\right\|^2 + \mathbb{E}_{\mu^0}\left\|S(\mu^\prime; x_i)\right\|^2.
\end{equation*}
Note that
\begin{eqnarray*}
& & \mathbb{E}_{\mu^0}\left\|S(\mu^\prime; x_i)\right\|^2\\
&=& \mathbb{E}_{\mu^0}\left\|x_i-\mu^\prime_{k_1}\right\|^2+\mathbb{E}_{\mu^0}\left\|x_i-\mu^\prime_{k_2}\right\|^2\\
&=& 2p + \frac{1}{K}\sum_{k=1}^K\left(\left\|\mu_{k_1}^\prime-\mu^0_k\right\|^2+\left\|\mu_{k_2}^\prime-\mu^0_k\right\|^2\right)\\
&\leq& 2\left\{p+\frac{1}{K}\sum_{k=1}^K\left(\left\|\mu_{k_1}^\prime-\mu^0_{k_1}\right\|^2+\left\|\mu^0_{k_1}-\mu^0_k\right\|^2+\left\|\mu_{k_2}^\prime-\mu^0_{k_2}\right\|^2+\left\|\mu^0_{k_2}-\mu^0_k\right\|^2\right)\right\}\\
&=& 2\left\{p+\left\|\mu_{k_1}^\prime-\mu^0_{k_1}\right\|^2+\left\|\mu_{k_2}^\prime-\mu^0_{k_2}\right\|^2 + \frac{1}{K}\sum_{k=1}^K\left(\left\|\mu^0_{k_1}-\mu^0_k\right\|^2+\left\|\mu^0_{k_2}-\mu^0_k\right\|^2\right)\right\}.
\end{eqnarray*}
If $\mu^0$ satisfies
\begin{equation*}
\sum_{k=1}^K\left(\left\|\mu^0_{k_1}-\mu^0_k\right\|^2+\left\|\mu^0_{k_2}-\mu^0_k\right\|^2\right)\lesssim p,
\end{equation*}
then combined with \eqref{Peq3-xpr}, we can conclude that
\begin{equation*}
\rho\left(\nabla^2\mathbb{E}_{\mu^0}m\left(\mu; x_i\right)-\nabla^2\mathbb{E}_{\mu^0}m\left(\mu^\prime; x_i\right)\right) \lesssim p+\left\|\operatorname{Vec}\left(\mu^\prime-\mu^0\right)\right\|^2+\left\|\operatorname{Vec}\left(\mu-\mu^\prime\right)\right\|^2.
\end{equation*}
As noted in Remark \ref{remark: the gap between two functionals}, it is sufficient to consider scenarios where $\mu^{\prime}$ is very close to $\mu^0$, then this inequality supports the first requirement of Condition \ref{condition: lipschitz}.
As for the second inequality of Condition \ref{condition: lipschitz}, from \eqref{hessian of log mu}, we can show that
\begin{equation*}
\nabla^2 \mathbb{E}_{\mu^0} \log p\left(x_i, c_i;\mu\right)=\nabla^2 \mathbb{E}_{\mu^0} \log p\left(x_i, c_i;\mu^\prime\right),
\end{equation*}
indicating that $s_2$ can be selected as any positive number. In summary, under specific assumptions, we have verified that Condition \ref{condition: lipschitz} is satisfied, with $\max \left(s_1, s_2\right) \geq 2$.

Now, all design conditions \ref{condition: local, twice continuously differentiable}-\ref{condition: lipschitz} have been checked. Denote $q^*\left(\operatorname{Vec}\left(\mu\right)\right)$as the VB posterior of $\mu$ after vectorizing. Similarly, let $\operatorname{Vec}\left(\hat{\mu}_n^*\right)$ represents the VB estimator after vectorization. If $\left\|\nabla\log p\left(\operatorname{Vec}\left(\mu^0\right)\right)\right\|^2\lesssim p^{\frac{3}{2}}$ and $r_0^2/p$ can be sufficiently large, Corollary \ref{corollary of theorem 4} demonstrates that 
\begin{equation*}
\|\operatorname{Vec}\left(\hat{\mu}_n^*-\mu^0\right)\| = O_p\left(\sqrt{\frac{p}{n}}\right).
\end{equation*}
Here, the condition $\left\|\nabla \log p\left(\operatorname{Vec}\left(\mu^0\right)\right)\right\|^2 \lesssim p^{\frac{3}{2}}$ is generally straightforward to satisfy, especially when assuming a multivariate standard normal prior and that each element of $\mu^0$ is uniformly bounded. 
Furthermore, if $p^3=o(n)$, Corollary \ref{Corollary 2} demonstrates that
\begin{equation*}
\sqrt{n} \alpha^{\top}\operatorname{Vec}\left(\hat{\mu}_n^*-\mu^0\right) / \sigma_\alpha \overset{\mathcal{L}}{\longrightarrow} N(0,1),
\end{equation*}
for any $\alpha \in S_{Kp}$, where $\sigma_\alpha^2=\alpha^{\top}V_{\mu^0}^{-1}\operatorname{Var}\left\{\nabla m\left(\mu^0; x\right)\right\}V_{\mu^0}^{-1}\alpha$. 

In summary, we have systematically verified the design conditions of our theoretical framework based on this proposed model. 

\bibliographystyle{apalike}
\bibliography{ref1}